\RequirePackage{fix-cm}
\documentclass[smallextended]{svjour3}       % onecolumn (second format)
\smartqed  % flush right qed marks, e.g. at end of proof
\usepackage{graphicx}
\RequirePackage[numbers]{natbib}
\usepackage{multirow}
\usepackage[ruled,vlined]{algorithm2e}
\usepackage{graphicx}
\usepackage{setspace}
\usepackage{float}
\usepackage{breakcites}
\usepackage{epsfig}
\usepackage{amssymb}
\usepackage{amsmath}
\usepackage{hyperref}
\usepackage{verbatim}
\usepackage{subcaption}
\usepackage{tikz}
\usepackage{color}
\newcommand{\AZ}{\textcolor{black}}

\usepackage{bbm}
\usetikzlibrary{arrows}
\usepackage[margin=1.2in]{geometry}

\tikzset{cross/.style={cross out, draw=black, minimum size=2*(#1-\pgflinewidth), inner sep=0pt, outer sep=0pt},
cross/.default={5pt}}
\usetikzlibrary{shapes.misc}
\usepackage{array}
\newcommand{\be}{\begin{eqnarray}}
\newcommand{\ee}{\end{eqnarray}}
\newcommand{\bea}{\begin{eqnarray*}}
\newcommand{\eea}{\end{eqnarray*}}

\newcommand{\ra}{\rightarrow}

\usepackage{nicefrac}
\numberwithin{equation}{section}
\usepackage{tikz}
 \def\firstcircle{(90:1.75cm) circle (2cm)}
 \def\secondcircle{(210:1.75cm) circle (2.5cm)}
 \def\thirdcircle{(330:1.75cm) circle (2.5cm)}
\newcolumntype{P}[1]{>{\centering\arraybackslash}m{#1}}
\usepackage{wrapfig}
\usepackage{caption}
 \newcolumntype{P}[1]{>{\centering\arraybackslash}m{#1}}

\def\a#1#2{\left ( {#1}\atop{#2} \right )}

\newcommand{\dd}{{b}}

\newcommand{\PP}{{\cal P}}
\newcommand{\RR}{{\mathbb{ R}}}
\newcommand{\QQ}{{\mathbb{ Q}}}
\newcommand{\T}{{\cal T}}
\newcommand{\D}{{ D}}
\newcommand{\KK}{{h}}
\newcommand{\VV}{{\mathbb{ S}}}
\newcommand{\UU}{{\mathbb{ B}}}
\newcommand{\X}{{\cal D}}
\newcommand{\Y}{{\cal Y}}

\newcommand{\Xb}{{\cal A}}
\newcommand{\suml}{\sum\limits}
\newcommand{\M}{{\cal X}}

\SetKwInput{KwInput}{Input}

\begin{document}

\title{Random and quasi-random designs  in group testing%\thanks{Grants or other notes
}

\titlerunning{Random and quasi-random designs  in group testing}        % if too long for running head

\author{ Jack Noonan   \and Anatoly Zhigljavsky  (Corresponding Author)
       %etc.
}


\institute{A. Zhigljavsky  \at
              School of Mathematics, Cardiff University, Cardiff, CF24 4AG, UK \\
              \email{ZhigljavskyAA@cardiff.ac.uk}           %  \\
           \and
           J. Noonan \at
              School of Mathematics, Cardiff University, Cardiff, CF24 4AG, UK \\
               \email{Noonanj1@cardiff.ac.uk}
}

\date{Received: date / Accepted: date}

\maketitle

\begin{abstract}
For  large classes of group testing problems,
  we derive lower bounds for the probability  that  all significant items are uniquely identified using specially constructed random designs.
  These bounds allow us to optimize  parameters of the randomization schemes.
We also suggest and numerically justify a procedure of constructing designs with better  separability properties than pure random designs. We illustrate theoretical considerations with a large simulation-based study. This study indicates, in particular, that in the case of the common binary group testing, the suggested families of designs have better separability than the popular designs constructed from disjunct matrices. We also derive several asymptotic expansions and discuss the situations when the resulting approximations achieve high  accuracy.

\end{abstract}

%\KEYWORDS{f} % Separate items with ;
%\AMSSUBJ{60G50;60G35} % Edit. Separate items with ;
%\AMSSUBJSECONDARY{60G70;94C12 } % Optional, separate items with ;

%\tableofcontents

\raggedbottom

\section{Introduction}

Assume that there are $n$ items (units, elements,  variables, factors, etc.)
                $a_1,\ldots,a_n$ with
 some of them  {\it defective}
(significant, important, etc.).
The problem of group testing (also known as ``pooling'' or ``factor screening'') is to determine  the defective items
  by testing a certain number of {\it test groups}~$X_j$.
 A design $\D_N=\{X_1, \ldots, X_N\}$  is a collection of $N$ test groups.
We assume that all test groups $X_j \in \D_N$ belong to some set $\X$ containing certain subsets
of the set
$\Xb
 =  \{a_1,\ldots,a_n \}.$ The set $\X \subseteq 2^\Xb$ will be called {\it design set.}

The group testing problems differ in the following aspects:
\vspace{-0.1cm}

 \begin{enumerate}
   \item[\,(i)\;] assumptions concerning the occurrence of defective items;
   \item [(ii)\,] assumptions on admissible designs;
   \item [\,(iii)] forms of the test function which provides observation results;
   \item [\,(iv)] assumptions on the number of allowed wrong answers (lies);
   \item [(v)\,] definitions of the problem solution.

 \end{enumerate}

\vspace{-0.1cm}
%
%
%{\bf (i)}
%assumptions concerning the occurrence of defective items;
%
%{\bf (ii)}
%assumptions on admissible designs;
%
%{\bf (iii)}
%forms of the test function which provides observation results;
%
%{\bf (iv)}
%definitions of the problem solution.

%\AZ{(v) lies ?}

 \noindent
The group testing  problems  considered in this paper are specified
 by the following properties.

\vspace{-0.1cm}

 \begin{enumerate}
   \item[(i)] As the main special case, we assume that
 there are exactly $d$ defective items with \mbox{$0\!<\!d\! \leq \! n$}.  Many statements, however, are formulated for the very general models defined by
 prior distributions for the number of defective items, see Section~\ref{sec:typ_ass}.
Moreover, a few results (e.g. Theorem~\ref{th:sec6.2} and  points three of Corollary~\ref{Main_binary_corollaries} and Corollary~\ref{simple_corollary}) cover the
problem of finding defectives the so-called binomial sample, where
the
 events ``item $a_i$ is defective'' are independent and have the same
 prior probability.

% Such setting is  frequently considered  in group testing literature \cite{Johnson,d2014lectures}.

   \item [(ii)] We only consider non-adaptive designs $\D_N= \{X_1,\ldots, X_N\} \subset \X$.
As the principal case,  we consider the design sets $\X$,
which contain the test groups  $X$ consisting of
exactly $s$  items with suitable $s$, see Section~\ref{sec:typ_ass}; { such designs are normally called constant-row-weight designs, see Section~\ref{sec:improve}. For brevity, we will call these designs simply   {\it constant-weight designs}}. The constant-weight random designs  seem to be  marginally more efficient than Bernoulli designs, where
 in order to build every test group  $X_j \in \D_N$,  each item is included into $X_j$  with  given probability; see Section~\ref{general_x} and Section~\ref{sec:4.1}.
%In addition, in Section~\ref{sec:improve} we suggest an alternative scheme of construction of nearly constant-weight designs that demonstrate very high practical efficiency (that is, providing high values of the probability $1-\gamma$ discussed in (v) below).
   \item [(iii)]
   Let $T \subset \Xb$ denote an unknown collection of defective items and
$X \subset \Xb$ be a test group.
We consider group testing models where the observation result for given $X$ and $T$ is
\be
\label{eq:f(X,T)}
f_{\KK}(X,T){=}
%f_K(X,T){=}
\min\{\KK,|X\cap T|\}\, ,
\ee
where
$|\cdot|$ stands for the number of elements in a discrete set and
 $\KK$ is a positive  integer.
In the most important special case of {\it binary} (or {\it disjunctive}) model, $\KK=1$. In this model, by inspecting a group $X\subset \Xb$
we receive 1 if there is at least one defective
item in $X$ and 0 otherwise.
In the {\it additive} (or ``adder'', in the terminology of \cite{d2014lectures}) model, $f_\infty(X,T) = | X \cap T |$ so that we choose $\KK=\infty$; in fact, any number between $n$ and
 $\infty$ can be chosen as $\KK$.
(In the additive model, after inspecting a group $X$
we receive the number of defectives in $X$.)
In the so-called {\it multiaccess channel} model, $\KK=2$.

   \item [(iv)] In the main body of the paper, we assume that the test results are {\it noiseless} (or {\it error-free}). In Section~\ref{lies_section} we show how most of our results can be extended  to the case of noisy testing, where
   up to $L$ {\it lies} (wrong answers, errors) are allowed. Moreover, in Section~\ref{sec:4.6} some specific results are specialized for the important case of
   binary group testing with lies.
    \item [(v)] As a rule, we are not interested in the designs that provide 100\% guarantee that  all defective items are correctly identified (in the group testing literature, this criterion is often referred to as ``zero-error probability criterion'' or ``exact recovery''). Instead,
         we are interested in studying the probability $1-\gamma$ that all defective items are discovered (for random designs) with the main theoretical contribution of this paper being the derivation of the lower bounds $1-\gamma^*$ for this probability; when it suffices to recover the defective set with high probability we are considering the small error probability criterion. Moreover, in Section~\ref{sec:improve}  we propose designs  that seem to provide very high values of $1-\gamma$,  even in comparison to the designs constructed from suitable disjunct matrices,  see
         Tables~\ref{disjunct_table} and \ref{disjunct_table2} in Section~\ref{sec:disj}.
 \end{enumerate}

Group testing is a well established area and has
attracted significant attention of  specialists in
optimum design, combinatorics, information theory and discrete search.
The origins of group testings can be traced back to the paper
\cite{Dorfman43}
devoted to adaptive
procedures of blood testing for detection of syphilitic men. Since then, the field of group testing has seen significant developments with extensive literature and numerous books dedicated to the field. The textbooks
\cite{DuH93,du2006pooling} and lecture notes \cite{d2014lectures} provide a background on group testing especially for zero-error non-adaptive problems.  An excellent introduction and summary of recent developments in  group testing and its connection to information theory can be found in \cite{Johnson}. The group testing problem in the binomial sample is especially  popular in the
group testing literature, see \cite{Johnson,SobelG59,Bruno}.

{Research in group testing often concentrates around the following important areas:
\begin{description}
  \item[(a)] construction of efficient designs (both, adaptive and non-adaptive);
  \item[(b)] studying properties of different families of designs;
  \item[(c)] derivation of upper and lower bounds for the lengths of designs providing either exact or weak recovery of the defective items;
  \item[(d)] extension of results in the noiseless setting for the case of noisy group testing;
  \item[(e)] construction of efficient decoding procedures to locate the defective items (given a design).
\end{description}
}
%Note that many available theoretical results are asymptotic established under the assumption that $n \to \infty$.

{In this paper, we touch upon all the above areas. In particular:
\begin{description}
  \item[(a)] in Section~\ref{sec:improve}    we develop a procedure of construction of a sequence
  of nested nearly doubly regular designs $D_1, D_2, \ldots$
  which, for all $N$, have large Hamming distances between all pairs $X_i,X_j \in D_N$ $(i \neq j)$  and, as a consequence, excellent separability properties (this is confirmed by a numerical study described in  Sections~\ref{sec:improve} and~\ref{sec:disj});
    \item[(b)] one of the main purposes of the paper is an extensive study of the probability of recovery of defective items for constant-weight  random designs (both, in non-asymptotic and asymptotic regimes);
  \item[(c)] as explained  in Remark 1 of Section~\ref{sec:exist},
most results on the probability of recovery of defective items can be reformulated as existence theorems of deterministic designs providing weak recovery; moreover,  in Sections  \ref{sec:add_exact} and \ref{sec:bin_exact} we derive asymptotic upper bounds for the lengths of deterministic designs providing exact recovery;
  \item[(d)] in Sections \ref{lies_section}, \ref{sec:4.6} and \ref{sec:noisy_as} we show how  most important results  obtained in  the noiseless setting can be extended for  the noisy group testing when up to $L$ lies are allowed;
  \item[(e)] in Section~\ref{sec:4.8} we numerically demonstrate that the so-called \AZ{Combinatorial Orthogonal Matching Pursuit (COMP)} decoding procedure alone could be  very inefficient; \AZ{see Section~\ref{sec:disj} for the definition of the COMP procedure.}
\end{description}
}

Existence theorems for group testing problems were extensively studied in Russian literature by M.B. Malutov, A.G.
Dyachkov, V.V. Rykov and other representatives of the Moscow
probability school, see e.g.  \cite{dr83,Tsyb}. The construction of upper bounds for the
length of optimal zero-error designs in the binary group
testing model has attracted significant attention; see \cite{DuH93} for a good survey.
In the papers
\cite{KatonaS83,Macula2,MaculaReuter},
the  construction schemes of group testing designs in important specific
cases, including the case of the binary model with two and, more generally,
 $d$ defectives, are studied.
Using probabilisitic arguments, existence theorems for designs under the zero-error criterion for the additive model
have been thoroughly studied in
\cite{ZhigljavskyZ95}. Motivated by the results of \cite{ZhigljavskyZ95}, in \cite{zhigljavsky2003probabilistic} expressions for the binary model were derived under the zero-error and small-error criterions. The results of \cite{zhigljavsky2003probabilistic} provided the inspiration for this paper. Note that there is a limited number of results on construction of optimal algorithms for
finding one, two or three defectives in search with lies, see e.g.
\cite{DeBonis,HillKarim,Macula97}. Some asymptotic expansions in existence theorems for general group testing problems have been
derived in \cite{zhigljavsky2010nonadaptive}.

{In the majority of papers devoted to construction of designs for the non-adaptive binary group testing problem, the designs are built from the so-called disjunct matrices, these are defined in Section~\ref{sec:disj}. Moreover, the COMP decoding procedure (according to COMP, all items in a negative test are identified as non-defective whereas all remaining items are identified as potentially defective, see Section~\ref{sec:disj})
is often used for identification of the set of defective items; see
e.g. a popular paper \cite{chan2014non} and
  a survey on non-adaptive group testing algorithms through the point of view  of decoding of test results \cite{chen2008survey}.
   Despite common claims, as explained in Sections \ref{sec:disj} and \ref{sec:4.8}, the designs based on the use of disjunct matrices are
inefficient  and the COMP decoding procedure alone leads to poor decoding.
}

{
In the asymptotic considerations, we assume that the
number of defective items is small relative  to the total number
of items $n$; that is, we consider a very sparse regime. Many results can be generalized to a sparse regime when $d $ slowly increases with $n$ but $d/n \to 0$ as $n \to \infty$.  There is  a big difference between the
asymptotic results in the sparse regime and results in the case when $d/n \to {\rm const}>0$ as $n \to \infty$. In particular, in view of  \cite{cantor1966determination,ErdosR63,lindstrom1964combinatory,Lindstrem75}, where the non-adaptive  group
testing problem for the additive model is considered with no
constraints on both the test groups and the number of defective
items, $N \sim
 {2n}/{\log_2 n},$ $n \ra \infty,$ for the minimal length of the
non-adaptive strategies that guarantee detection of all defective items. For fixed $d$, the best known explicit constructions of designs come from number theory \cite{bose1962theorems,lindstrom1969determination} and are closely related to the concept of Bose-Chaudhuri-Hocquenghem codes. For these constructions it is shown that $N \leq d \log_2 n(1 + o(1)) $ tests are required. For $d\geq  3$, the best currently known construction is with $N \leq 4d \log_2 n/ \log_2 d(1 + o(1))$ and can be obtained from results of \cite{d1981coding,poltyrev1987improved}. This result is constructed using random coding and is shown to be order-optimal.
%It is well known that non-adaptive designs are doomed to fail with high probability $d = \alpha n$ if $n \rightarrow \infty$ and $\alpha$ is constant.
}

{
In the very sparse regime with $d$ constant and $n\ra\infty$, the best
known upper bound for the length of zero-error designs
in the binary group testing problem has been derived in
 \cite{drr},
see also Theorem~7.2.15 in
\cite{DuH93}:
$ N\leq \frac1{2}{dc_d}(1+o(1))\log_2 n$,
where
$$1/{c_d}=\max\limits_{0\leq  q\leq 1} \max\limits_{0\leq Q\leq 1}
\left\{ -(1-Q)\log_2(1-q^d)+d\left[Q\log_2\frac{q}{Q}
+(1-Q)\log_2\frac{1-q}{1-Q}\right]\right\}$$
and
$c_d={d\log_2 e}(1+o(1))\;\mbox{as}\;
d\ra\infty.$
Asymptotically, when both $n$ and $d$ are large,
this is a marginally
%\marginpar{check it in the library}
 better bound
than  the asymptotic bound
\bea
N \leq N_*(n,d) \sim \frac{e}{2} d^2 \log n \, , \;\;
n \ra \infty,\; d=d(n) \ra \infty,\; d(n)/n \ra 0\, ,
\eea
which has been  derived
in
\cite{dr83} by  the probabilistic method based on the use of the Bernoulli design.
Exactly the same upper bound can be obtained using
random constant-weight designs, see Corollary~5.2 in
\cite{zhigljavsky2003probabilistic}. 
\AZ{Development of existential (upper) bounds for group testing designs for  binary group testing has 
has been complemented by establishing various lower bounds; for  comparison of the lower and upper bounds, see  the well-written Section 7.2 of \cite{DuH93}. } }

{Primarily for the binary model, notable contributions in recent years are as follows.
In \cite{aldridge2014group}, the authors consider the problem of nonadaptive noiseless group testing problem using Bernoulli designs and describe a number of algorithms used to locate the defective set after the design has been constructed; one of these is the COMP procedure which will be discussed in Section~\ref{sec:disj}. For bounds on the number of tests when using Bernoulli designs, also see \cite{scarlett2016limits,scarlett2016phase}. In \cite{aldridge2016improved}, instead of Bernoulli designs  the authors consider designs where each item is placed in a constant number of tests. The tests are chosen uniformly at random with replacement so the test matrix has (almost) constant column weights, these terms will be fully explained in Section~\ref{sec:improve}. The authors show that application of the COMP detection algorithm with these constant column-weight-designs significantly increases detection of the defective items in all sparsity regimes. This (almost) constant-column-weight property will be discussed further in Section~\ref{sec:improve} where it will be combined with a Hamming distance constraint to improve the probability of separation. In  \cite{coja2020information}, for the randomised design construction discussed in \cite{aldridge2016improved}, the authors provide a sharp bound on the number of tests required to locate the defective items. In \cite{coja2020optimal},
the authors consider existence bounds for both a test design and an efficient algorithm that solve the group testing problem with high probability. In \cite{mezard2011group}, the authors consider the binomial sample group testing problem where each item is defective with probability $q$. The authors construct a class of two-stage algorithms that reach the asymptotically optimal value of $nq|\log(q)|$. The asymptotic bounds for the one-stage (nonadaptive) setting for the binomial sample problem are studied in \cite{mezard2008group}.
}

{
This paper differs from the aforementioned papers in the following aspects: (a) the majority of known theoretical results require  large $n$ and only numerical evidence is presented  when $n$ is small; this paper, however, provides rigorous results for any $n$ where many asymptotic results do not apply; (b) the asymptotic expansions in this paper provide constants that have crucial significance  when $n$ is only moderately  large (this additional constant term is not present in many asymptotic results for group testing); (c) many of the previously cited papers use decoding procedures that do not guarantee identification of the defective set even if it is possible to locate it. Procedures like COMP are fast to execute, and as previously mentioned, with certain design constructions can in a large number of cases locate the defective set. However, in this paper we will use decoding procedures that will guarantee the location of the defective set if this is possible given the design.
}

{
By requiring a given design to satisfy the constraint of being able to find the defective items, we are considering an example of a (random) constraint satisfaction problem (CSP). Many of the main advances of this paper can be viewed as the careful counting of satisfying assignments for a CSP, where the satisfying assignments can correspond to tests that are able to differentiate between different subsets of $ \Xb$.
 The techniques used in this paper are related to approaches used in the random CSPs literature, see for instance \cite{zdeborova2016statistical}. However group testing problems are very specific  and cannot be simply considered as specific application of the general CSP methodology.
 }

The rest of the paper is organized as follows.
In Section~\ref{sec:2} we  develop a general methodology of derivation the lower bounds for $1-\gamma$,
the probability that all defective items are uniquely identifiable from test results taken according to constant-weight random designs  and establish several important auxiliary results. In Section~\ref{sec:3} we derive  lower bounds for $1-\gamma$ in a general group testing problem and consider the case of additive model for discussing examples and numerical results. The more practically important case of the binary model is treated in Section~\ref{sec:4}. Section~\ref{sec:2} is devoted to asymptotic existence bounds and construction of accurate approximations. In Appendix~A we provide some proofs and in Appendix B we formally describe the algorithm of Section~\ref{sec:improve}. Let us consider the content of Sections~\ref{sec:2},~\ref{sec:3}, \ref{sec:4} and~\ref{sec:5} in more detail.

In Section~\ref{sec:2.1}   we discuss general discrete search problems.
In Section~\ref{sec:2.2} we develop the general framework for derivation of the upper bounds $\gamma^*$ for $\gamma$, the probability that for a random design all defective items cannot be recovered; the main result is formulated as Theorem~\ref{th:sec2.2}. Theorem~\ref{th:22} of  Section~\ref{lies_section} extends Theorem~\ref{th:sec2.2} to the case when some of $N$ test results are allowed to be wrong (the case of lies).
In Section~\ref{sec:exist} we show how many of our results can be reformulated in terms of
existence bounds in the cases of weak and exact recovery.
In Section~\ref{sec:typ_ass} we consider different assumptions  on the occurrence of defective items and
the  randomisation schemes used for the construction of the randomized designs. In Sections~\ref{sec:2.5} and ~\ref{sec:2.6} we
formulate two  important combinatorial results, Lemmas~\ref{th:q} and \ref{th:R}.

In Section~\ref{sec:3.1} we derive upper bounds $\gamma^*$ for $\gamma$ for a general test function  \eqref{eq:f(X,T)} in  the most important case $\X=\PP_n^{s}$; that is,  when all $X_i \in\D_N$ have exactly $s$ items (see \eqref{eq:G_k} for the formal definition of $\PP_n^{s}$).
In Section~\ref{sec:3.2} we specialize the general results of Section~\ref{sec:3.1} to a relatively easy case of the additive model and consider special instances of the information about the defective items
 including the case of the binomial sample case, see Corollary~\ref{Main_binary_corollaries}. In Section~\ref{sec:3.3} we provide some results of simulation studies for the additive model. In Section~\ref{general_x} we show how to extend the results established for the case $\X=\PP_n^{s}$ to cover other randomization schemes for choosing the groups of items $X_i$ including the case of Bernoulli designs.

In Section~\ref{sec:4.1} we provide a collection of upper bounds $\gamma^*$ for $\gamma$ for different instances of
 the binary model. All  results formulated in this section follow from general results and specific considerations  of Sections~\ref{sec:3.1} and~\ref{general_x}. In Section~\ref{sec:4.2}, we illustrate  some of the theoretical results formulated in Section~\ref{sec:4.1} by results of simulation studies.  In Section~\ref{sec:improve} we develop
                 a procedure for construction of a sequence
  of nested nearly doubly regular designs $D_1, D_2, \ldots$
  which, for all $N$, have large Hamming distances between all pairs $X_i,X_j \in D_N$ $(i \neq j)$.  With the help of numerical studies we
   also demonstrate  excellent separability properties of the resulting designs.
 In Section~\ref{quasi_sim} we apply the technique of Section~\ref{sec:improve} and numerically demonstrate that indeed the resulting designs provide a superior separability relative to random designs.
In Section~\ref{sec:disj} we numerically compare random, improved random of Section~\ref{sec:improve} and the very popular designs constructed from the disjunct matrices. In particular, we find that improved random designs have a better separability than the designs constructed from the disjunct matrices, see
Tables~\ref{disjunct_table} and \ref{disjunct_table2}. In Section~\ref{sec:4.8} we discuss the (in)efficiency of the COMP decoding procedure.
In Section~\ref{sec:4.6} some specific upper bounds  are specialized for the
   binary group testing with lies; simulation results are provided to illustrate theoretical bounds.

{In Section~\ref{sec:5.1} we describe the technique used to transform finite-$n$ results into the asymptotic expansions. A very important feature of the developed expansions is that in the very-sparse  regime we have explicit expressions for the constant term, additionally to the main term involving $\log n$. Sections~\ref{sec:add_exact}, \ref{sec:bin_exact} and \ref{sec:as_b_weak} we apply results of
Section~\ref{sec:5.1} respectively to the cases  of additive model (both exact and weak recoveries), binary model with exact recovery and the
binary model with weak recovery. Results of these sections clearly demonstrate the following: (a) weak recovery is much simpler  than exact
recovery, (b) the constant terms in the asymptotic expansions play an absolutely crucial role if these expansions are used as approximations, and (c) the resulting approximations have rather simple form and are very accurate already for moderate values of $n$.
Finally, in Section~\ref{sec:noisy_as} we discuss   a technique of transforming the asymptotic upper bounds for $N$ for noise-free group testing problems into upper bounds for $N$ in the same model  when up to $L$ lies are allowed.
}

\section{General discrete search problem, random designs and the probability of solving the problem }
\label{sec:2}

\subsection{Problem statement}
\label{sec:2.1}

We consider the group testing problems
from the  general point of view of
discrete search.
Following~\cite{geran}
a discrete {\em search problem} can often be
determined as a quadruple
$\{\T,\X,f, \Y \}$,
where $\T=\{T\}$ is a {\it target set},
which is an ordered collection
of all possible {\em targets}  $T$, $\X=\{X\}$ is a {\it design set}, a collection of all allowed test groups $X$,
 and
\mbox{$f:\X\times
\T\rightarrow\Y$} is a
{\em test function} mapping $\X\times\T$ to $\Y$, the set of all possible outcomes of a single test.
In group testing, the targets $T$ are allowed  collections of defective items and a value $f(X,T)$ for fixed $X\in\X$ and $T\in\T$ is
a test result at the test group  $X$ under the assumption that the  unknown target is $T$.
 For a pair of targets
$T_i, T_j \in {\cal T}$,
we say that
$X \in \X$
{ separates} $T_i$ and
$T_j$ if
$f(X,T_i) \neq f(X,T_j)$.
We say that a design
$\D_N=\{ X_1, \ldots, X_N\}$
{separates} $T\in{\cal T}$ if
for any
$T' \in {\cal T}$, such that $T'\neq T$,
there exists a test group
$X \in \D_N$ separating the pair  $(T,T')$. We only consider {\it solvable}  search problems where
each $T \in {\cal T}$ can be uniquely identified from  test results at all $X \in \X$.

In this paper, we are interested in studying properties of random designs for solving group testing problems.
Let $\RR$ and $\QQ$ be distributions on $\X$ and $\T$ respectively.
Let $\D_N= \{X_1,\ldots, X_N\} $ be a random $N$-point design with mutually independent and $\RR$-distributed test groups $X_i \,\,(i=1,\ldots,N)$ and let  $T\in \T$ be a $\QQ$-distributed random target. For a random $N$-point design $\D_N$,  we are interested in estimating the value of $\gamma=\gamma(\QQ, \RR,N) $ such that
\be \label{weak_solution}
{{\rm Pr}_{\QQ,\RR}\{ T \textrm{ is separated by } \D_{N}  \} } = 1-\gamma \,.
\ee
The intractable nature of the l.h.s in \eqref{weak_solution} makes it (unless the problem is very easy and hence impractical) impossible to
explicitly compute  $\gamma$. One of the main aims of this paper is the derivation  of explicit upper bounds
$\gamma^*=\gamma^*(\QQ, \RR,N)$ for $\gamma$ so that
\be\label{weak_solution_bound}
{{\rm Pr}_{\QQ,\RR}\{ T \textrm{ is separated by } \D_{N}  \} } \geq 1-\gamma^* \,.
\ee
This will allow us to state that a random design $\D_N$ solves the group testing problem with probability at least $1-\gamma^*$.

Another way of interpreting the results of the form \eqref{weak_solution_bound} is as follows.  For a  given search problem  $\{\T,\X,f, \Y \}$, an algorithm  of generating the test groups  $X_1, X_2, \ldots$ and $\gamma \in (0,1)$, define
$N_\gamma$ to be the smallest integer $N$ such that
\be \label{inverse problem}
{{\rm Pr}_{\QQ,\RR}\{ T \textrm{ is separated by } \D_{N}  \} } \ge 1-\gamma \,,
\ee
where the probability is taken over randomness  in $T$ and   $X_1, X_2, \ldots$
Computation of the exact value of $\gamma$ is a very difficult problem. However, as formulated in the following lemma, the  ability of  computing any upper bound $\gamma^*=\gamma^*(N)$ for $\gamma$ in \eqref{inverse problem} implies the possibility of derivation of the corresponding upper bound for $N_\gamma$.
\begin{lemma}
%{ \it
 Let
$\{\T,\X,f, \Y\}$
be a solvable discrete search problem with random $T$,   $X_1, X_2, \ldots$ be a sequence of test groups $X_i \in \X$ and
$\gamma^*=\gamma^*(N)$ be an upper bound for $\gamma$ in \eqref{inverse problem} for a design $D_N=\{X_1, \ldots, X_N\}$.
 Then for any
 $0\!<\!\gamma\!<\!1$, \eqref{inverse problem} is satisfied for any $N\geq N_{\gamma} $ where
\be
\label{eq:N_gamma}
N_{\gamma} := \min\Biggl\{\! N =1,2,\ldots:  \gamma^*(N)<  \gamma\Biggl\} \,.
\ee
\end{lemma}

\begin{remark} \label{rem1}
{Even if the test groups $X_1, X_2, \ldots$ leading to \eqref{inverse problem} are random, from formula \eqref{inverse problem} with $N=N_\gamma$ we deduce that there exists a deterministic design
 $\D_N= \{X_1,\ldots, X_N\} $ with $N \leq N_\gamma$ such that \eqref{inverse problem} holds, where the probability in \eqref{inverse problem} is taken over $\QQ$ (random $T$) only. This follows from the discreteness of the space of all $N$-point designs and that the expectation of the event ``$T \textrm{ is separated by } \D_{N} $'' with respect to  random designs is the l.h.s. in \eqref{inverse problem}.}
\end{remark}

\subsection{A general technique for derivation of  upper bounds $\gamma^*=\gamma^*(\QQ, \RR,N)$ for $\gamma$}
\label{sec:2.2}

For fixed $T_i$ and $T_j \in \T$, let
\be \label{eq:pij}
p_{ij}= \mbox{Pr}_{\RR}\{f(X,T_i) = f(X,T_j)\}\,
\ee
be the probability that the targets
$T_i$ and $T_j$
are not separated by one random test  $X \in \X$, which is
distributed according to $ \RR$.
The following theorem is a straightforward application of the union bound.

%{\bf Theorem 2.2.}
\begin{theorem}
 \label{th:sec2.2}
%{ \it
 Let
$\{\T,\X,f, \Y\}$
be a solvable discrete search problem with $\RR$ and $\QQ$ being any distributions on $\X$ and $\T$ respectively. For a fixed $N\geq1$, let $\D_N= \{X_1,\ldots, X_N\} $ be a random $N$-point design with each $X_i \in \D_N $ chosen independently and $\RR$-distributed. Then for $\gamma=\gamma(\QQ, \RR,N) $  of  \eqref{weak_solution}, we have $\gamma(\QQ, \RR,N) \leq \gamma^*(\QQ, \RR,N)$ with
\be
\label{new_weak}
 \gamma^*(\QQ, \RR,N) =\sum_{i=1}^{|\T|}{\rm Pr}_\QQ\{ T=T_i\}\sum_{j\neq i}p_{ij}^N  \,.
\ee
\end{theorem}
{\bf Proof.}
By applying the union bound,  the probability  that $T_i$ is  not
separated from at least one $T_j\in\T$ $ (T_j\neq T_i) $
after $N$ random tests is less
than or equal to
$
\sum_{j \neq i} \left(p_{ij}\right)^N
$
and we thus have
$
1 - \sum_{j \neq i} \left(p_{ij} \right)^N
$
as a lower bound for the probability that $T_i$ is separated from
all other $T_j \in \T$. Averaging over $T_i$ we obtain
\bea
{{\rm Pr}_{\QQ,\RR}\{ T \textrm{ is separated by } \D_N  \} }
&=& \sum_{i=1}^{|\T|}{\rm Pr}_{\RR}\{ T_i \textrm{ is separated by } \D_N  \} {\rm Pr}_\QQ\{ T=T_i\} \nonumber \\
&\ge& 1- \sum_{i=1}^{|\T|}{\rm Pr}_\QQ\{ T=T_i\}\sum_{j\neq i}p_{ij}^N = 1- \gamma^*(\QQ, \RR,N)\,. \label{key_step}
\eea
The statement of the theorem follows. \hfill $\Box$

For the very  common scenario when $\QQ$ is uniform  on $\T$, that is  ${\rm Pr}_\QQ\{ T=T_i\}=1/|\T|$ for all $i=1,\ldots |\T|$,
the formula \eqref{new_weak} for $\gamma^*(\QQ, \RR,N)$ simplifies to
\be
\label{eq:missing}
\gamma^*(\QQ, \RR,N) =\frac{2}{|\T|}\sum_{i=1}^{|\T|}\sum_{j=1}^{i-1}p_{ij}^N  \,.
\ee

Note also that the in order to apply  the upper bound \eqref{new_weak}, the test function $f(X,T)$ does not have to be of the form \eqref{eq:f(X,T)}. Indeed,  this bound can be used for many discrete search problems of different nature from group testing; in particular, for solving the  ``Mastermind'' game \cite{o1993mastermind}.

\subsection{Extension to the case when several lies (errors) are allowed}\label{lies_section}

 Assume the so-called {\it $L$-lie search problem}, where
 up to $L$  test results $Y(X_j,T)$ at some $X_j \in \D_N= \{X_1,\ldots, X_N\} $ may differ from $f(X_j,T)$.
For a random $N$-point design $\D_N$  we are interested in bounding the value of $\gamma$, $0<\gamma<1$, such that
\bea
{\rm Pr}_{\QQ,\RR}\{ T \textrm{ can be uniquely identified by } \D_{N}  \textrm{ with at most $L$ lies} \} = 1-\gamma \,.
\eea

An important observation is that if
a non-adaptive  design
$\D_N=\{X_1,\ldots,X_N\}$
is applied in a general $L$-lie search problem, then one can
guarantee that the target can be uniquely identified
if and only if the two vectors
$F_{T}=(f(X_1,T),\ldots,f(X_N,T))$ and
$F_{T'}=(f(X_1,T'),\ldots,f(X_N,T'))$ differ in at least
$2L+1$ components where $(T,T')$ is any pair of different targets in
$\T$.
That is, a target $T\in \T$ can be uniquely identified if and only if for all
$T'\in\T\setminus\{ T\}$
\be
\label{eq:min-dH}
d_H(F_T,F_{T'})\geq 2 L+1\, ,
\ee
where $d_H(a,a')$ is the Hamming distance between two $n$-vectors $a$ and $a'$; that is, the number of components of $a$
and $a'$ that are different.

The following statement is a generalization of Theorem~\ref{th:sec2.2}   to the
case of $L$-lie search problem.

\begin{theorem} \label{th:22}
Let
$\{\T,\X,f,\Y\}$ be a solvable $L$-lie search problem with $\RR$ and $\QQ$ being any distributions on $\X$ and $\T$ respectively.  For a fixed $N\geq1$, let $\D_N= \{X_1,\ldots, X_N\} $ be a random $N$-point design with each $X_i \in \D_N $ chosen independently and $\RR$-distributed. Then
\be
\label{eq:exist_lies}
\gamma^*(\QQ, \RR,N)  =
\sum_{i=1}^{|\T|}{\rm Pr}_{\QQ}\{ T=T_i \} \ \sum_{j\neq i} \sum_{l=0}^{2L} {{N}\choose {l}}
\left(p_{ij}\right)^{N-l}
\left(1-p_{ij}\right)^l  \, .
 \ee
\end{theorem}

Proof of Theorem~\ref{th:22}  can be found in Appendix A.
Theorem~\ref{th:22} can be seen as  a generalisation of Theorem 9 of \cite{zhigljavsky2003probabilistic}.
Note that in the most important case when  $\QQ$ is uniform on $\T$, \eqref{eq:exist_lies} becomes
\be
\label{eq:exist_lies0}
\gamma^*(\QQ, \RR,N)  = \frac{2}{|\T|}
\sum_{i=2}^{|\T|} \sum_{j=1}^{i-1} \sum_{l=0}^{2L} {{N}\choose {l}}
\left(p_{ij}\right)^{N-l}
\left(1-p_{ij}\right)^l  \, .
 \ee

One can consider a version of the $L$-lie search problem where all wrong
answers  are the same; that is, the wrong results  are equal to some
$y\in\Y$, and this value $y$  can be obtained by
correct answers as well. This problem
is a little simpler than the general $L$-lie problem and in this problem
it is enough
to ensure
$
d_H(F_{T},F_{T^{'}}) \geq L+1
$
rather than $(\ref{eq:min-dH})$, to guarantee the unique identification of the defective set. For this problem the upper bound is:
\be
\label{eq:lie_red0}
\gamma^*(\QQ, \RR,N)  =
\sum_{i=1}^{|\T|}{\rm Pr}_{\QQ}\{ T=T_i \} \ \sum_{j\neq i} \sum_{l=0}^{L} {{N}\choose {l}}
\left(p_{ij}\right)^{N-l}
\left(1-p_{ij}\right)^l  \, .
\ee

For several setups of the group testing problem, we will derive
 closed-form expressions for $p_{ij}$;
we therefore can easily compute the upper bounds (\ref{eq:exist_lies}) and
(\ref{eq:lie_red0}) for the corresponding $L$-lie group testing problems as well.
These bounds will be very similar to the ones formulated for problems with no lies but with an extra summation in the right-hand side.

{

\subsection{Existence bounds in the cases of weak and exact recovery}
\label{sec:exist}

%\begin{remark} \label{rem2}
As was noted in Remark~\ref{rem1}, $N_\gamma$ of \eqref{eq:N_gamma} has the following interpretation as an existence  bound in the case of weak recovery: for a given $\gamma \in (0,1)$ and  any $N\geq N_\gamma$, there exist deterministic designs $D_N$ such that ${{\rm Pr}_{\QQ}\{ T \textrm{ is separated by } \D_{N}  \} } \ge 1-\gamma$. In the most important case when  $\QQ$ is uniform on $\T$,  in view of \eqref{eq:missing}, we can write the existence bound $N_\gamma$ of \eqref{eq:N_gamma} as

\be
\label{eq:exist_g}
N_{\gamma} = \min\Biggl\{\! N =1,2,\ldots:
\sum_{i=2}^{|\T|} \sum_{j=1}^{i-1}
p_{ij}^{N} < \frac{ \gamma |\T|}2 \Biggl\} \,.
\ee

In case of exact recovery, we need to separate all possible pairs $(T,T')\in \T\times \T$.
Let, as in Theorem~\ref{th:sec2.2}, $\D_N= \{X_1,\ldots, X_N\} $ be a random $N$-point design with independent $\RR$-distributed test groups $X_i  $.  By the union bound, similarly to the proof of Theorem~\ref{th:sec2.2}, the probability that at least one pair $(T,T')\in \T\times \T$ is not separated by $D_N$, is not larger than
$\sum_{i=2}^{|\T|}\sum_{j=1}^{i-1}p_{ij}^N $. If this expression is smaller than 1, then, by the discreteness of $\T$, there is at least one deterministic design $\D_N= \{X_1,\ldots, X_N\} $ separating all $(T,T')\in \T\times \T$. The smallest $N$ when this happens is
\be
\label{eq:N_0}
N_{0} := \min\Biggl\{\! N =1,2,\ldots:  \sum_{i=2}^{|\T|}\sum_{j=1}^{i-1}p_{ij}^N  < 1\Biggl\} \,
\ee
and  for all $N\geq N_0$
there exist deterministic designs $D_N$ guaranteeing unique identification of the unknown target $T\in \T$.

By comparing \eqref{eq:exist_g} and \eqref{eq:N_0} we observe that if we set $\gamma= 2/ |\T|$ then $N_{\gamma}$ and $N_0$ coincide so we might suggest that $N_0$ is the limit of $N_{\gamma}$ as $\gamma \to 0$. This intuition rarely works, however, as in typical group testing problems values of $|\T|$ are astronomically large but values of $\gamma$ are simply small. As we demonstrate in several subsections of Section~\ref{sec:5}, weak recovery is indeed a much simpler problem than exact recovery, at least in the case of fixed $\gamma>0$.

Assume now that  up to $L$ lies are allowed. Similarly to \eqref{eq:N_0} and using  \eqref{eq:exist_lies0}, we deduce that there are deterministic designs $D_N$ guaranteeing unique identification of the unknown target $T\in \T$ if $N\geq N_{0,L}$ where
\bea
\label{eq:N_0L}
N_{0,L} := \min\Biggl\{\! N =2L,2L+1,\ldots:  \sum_{i=2}^{|\T|} \sum_{j=1}^{i-1} \sum_{l=0}^{2L} {{N}\choose {l}}
\left(p_{ij}\right)^{N-l}
\left(1-p_{ij}\right)^l < 1\Biggl\} \,.
\eea

}

\subsection{Typical target and design sets  and assumptions on the randomisation schemes $\QQ$ and $\RR$ in  group testing}
\label{sec:typ_ass}

{
In group testing problems, the target set $\T$   has, as a rule, a very particular structure considered below.
Denote the  collection of all subsets of $\Xb=\{a_1, \ldots, a_n\}$ of length $k$  by $\PP_n^k$:
\be
\label{eq:G_k}
\PP_n^k= \left\{   (a_{i_1},\dots,a_{i_k} ), \; 1\leq i_1<\dots i_k\leq n
 \right\} .
\ee
The collection of groups of items containing  $k$ items or less will be denoted by
$
\label{eq:G_leq-k}
\PP_n^{\leq k}=\bigcup_{j=0}^k\PP_n^j ,
$ where $\PP_n^{0}= \emptyset$.
All target sets $\T$ considered in this paper will have the form
$\label{eq:25}
\T =  \cup_{j \in B } \PP_n^j\, ,
$
where $B  $  is a subset of $ \{0,1,\ldots,n \}.$
The main choices of $B $ are $B  =\{d \}$ and $B =\{0,1,\ldots,d \}$  for $1\leq d\leq n$; this corresponds to $\T = \PP_n^d$ and $\T= \PP_n^{\leq d}$ respectively.

The distribution $\QQ$  for $T \in \T$ defines  the assumptions on  the occurrence of defective items.
In a typical group testing setup, $\QQ$ has  the property of exchangeability; that is, symmetry with respect to  re-numeration of the items.
We express this as follows.
Let $\UU$ be a probability distribution on $ \{0,1,\ldots,n \}$ and $\xi$ be a $\UU$-distributed random variable.
Then for a $\QQ$-distributed  random target $T\in \T$ and any $ j \in  \{0,1,\ldots,n \}$:
 \be\label{Q_restrictions}
 \mbox{Pr}_{\QQ}\{ |T| = j  \} =  \mbox{Pr}_\UU\{ \xi =j  \}\; {\rm and} \;
 \mbox{Pr}_{\QQ}\{ T = \textsf{T}    \, | \,|T| ={j} \} = \left\{
                                                   \begin{array}{cc}
                                                       1/{{n \choose {j} }} & {\rm if}\;\; \textsf{T}  \in  \PP_n^{j} \\
                                                     0 & {\rm otherwise}
                                                   \end{array}
                                                 \right.
\,,
 \ee
 where the term   ${{n \choose {j} }} $
 %$\left(       \begin{array}{c}                 n \\            j \\                 \end{array}                                                  \right)$
                                                   is the number of elements in  $\PP_n^j.$
                     In the main two particular cases, when $\T=\PP_n^{d}$ and  $\T=\PP_n^{\leq d}$, the measure $\UU$ is concentrated on the one-point set $\{d\}$ and on $\{0,1, \ldots,d\}$, respectively.
The assumption  \eqref{Q_restrictions} can also be expressed as follows: $\forall j \text{ and } \forall \textsf{T}  \in  \PP_n^{j}  $
\bea
 \mbox{Pr}_{\QQ}\{ T= \textsf{T}    \} ={\mbox{Pr}_\UU\{ \xi=j \}}/{ {n \choose j}}  \,.
 \eea

The main objective of choosing the design set $\X$ (as well as the randomization scheme $\cal R$) is the efficiency of the resulting group testing procedure. Bearing this in mind, we mostly use  $\X=\PP_n^s$ with suitable $s$.
As a rule, in this case we  achieve better bounds than, say,
in the case $\X=\PP_n^{\leq s}$,
with  optimal  $s$ as well as in the case of  Bernoulli designs, when each item is included into a test  group with  probability  $p$, with optimal $p$; see Table~\ref{compare_binom2}.

For the main choice ${\X}  = \PP_n^{s}$, we choose the distribution $\RR$ to be the uniform  on ${\X} $ so that
$
 \mbox{Pr}_{\RR}\{ X= \textsf{X}  \} ={1}/{ {n \choose s}}
$  for all $\textsf{X} \in  \PP_n^{s}  $.
For this choice of    $\RR$, we can rewrite the probabilities $p_{ij}$ of \eqref{eq:pij} as
$
p_{ij} ={k_{ij}}/{|{\X} |} ={k_{ij}}/{{n \choose s}}  \, ,
$
where
\bea
k_{ij} = \left| \{ X\in {\X }:\; f(X,T_i) = f(X,T_j)\} \right| \;\;\;\;\;
\mbox{\rm for $\;\;T_i,T_j  \in {\T}$}\, .
\eea
In accordance with  \cite{geran}, these coefficients will be called {\it R\'{e}nyi coefficients}.
As shown below, computation of these coefficients involves  some counting only.

}

\subsection{An important auxiliary result} \label{sec:2.5}
%\section{Computation of the R\'{e}nyi coefficients}

Consider integers  $m, l$ and $p$  satisfying the conditions
$0\leq p\leq m\leq l\leq n$ and $p<l$.
Denote
\be
\label{eq:T}
\T(n,l,m,p) = \{ ( T , T') \in
  \PP_n^{\leq n}
\times
  \PP_n^{\leq n}: \;
 | T | = l,\; | T' | = m,\; | T \cap T' | = p \} \subset \PP_n^{l}
\times
  \PP_n^{m}\, .\;\;\;
\ee
Note that the condition $p<l$ guarantees that $ T \neq  T'$
for all pairs $( T , T' ) \in \T(n,l,m,p).$ {${\cal T}(n,l,m,p)$ is simply the collection of pairs of assignments $(T, T')$ of defective items such that $T$ contains $l$ defective items, $T'$ contains $m$ defective items and they have exactly $p$ defective items in common.
Interpretation for the numbers $l,m$ and $p$ is given on Figure~\ref{XandRfigure} (left).}

The following lemma
 allows computing the number of elements in the
sets
(\ref{eq:T}).

\begin{lemma}
\label{th:q}
The number of different non-ordered pairs in $\T(n,l,m,p)$ equals
\be
\label{eq:Q}
Q(n,l,m,p)=\left\{\begin{array}{ll}
            {{n} \choose {\,p\;m-p\;l-p\;n-l-m+p\,} }&  {\rm if} \; m<l  \\
                                        &  \\
            \frac{1}{2} { {n}\choose {\,p\;m-p\;m-p\;n-2m+p\,}} & {\rm if} \; m=l \,,
                   \end{array}\right. \,
\ee
where
\bea
{{n} \choose {n_1\;n_2 \dots n_k}} = \left\{
                                       \begin{array}{cc}
                                         \frac{n!}{n_1!n_2!\dots n_k!} & \mbox{\rm if }
n_r\geq 0, \; \sum_{r=1}^k n_r=n \\ &\\
 0 & \mbox{\rm if }\,  \min\{n_1,\ldots , n_k \} <0 \\
                                       \end{array}
                                     \right.
\eea
is the multinomial coefficient.
\end{lemma}
{For the proof of \eqref{eq:Q}, which only involves simple counting arguments, see Theorem~4.1 in
\cite{ZhigljavskyZ95}. Note the coefficient $\frac12$ in \eqref{eq:Q} for the case $l=m$; it is related to the fact that $Q(n,l,m,p)$
is the number of {\it non-ordered} pairs $(T,T')$ in $\T(n,l,m,p)$.}

\subsection{Balanced design sets} \label{sec:2.6}

\begin{figure}[h]
\centering
\begin{minipage}{.49\textwidth}
  \centering
\vspace{-1.25cm}
\begin{tikzpicture}[scale=0.6]
    \begin{scope}[shift={(3cm,-5cm)}, fill opacity=0.7]
        \fill[white] \firstcircle;
        \fill[white] \secondcircle;
        \fill[white] \thirdcircle;
      %  \draw \firstcircle ;
        \draw \secondcircle ;
        \draw \thirdcircle;
%               \node [scale=3]   at (0,4.5) {$X$};
                              \node [scale=1.5]   at (-2.5,-1) {$l\!-\!p$};
                            \node [scale=1.5]   at (2.5,-1) {$m\!-\!p$};
                              \node [scale=2]   at (-4,-3) {$T$};
                                     \node [scale=2]  at (4,-3) {$T'$};
 %        \node  [scale=2]  at (0,0.3) {$r$};
   %               \node  [scale=2]  at (-1,0.8) {$u$};
   %                \node  [scale=2]  at (1,0.8) {$v$};
                  \node  [scale=1.5]  at (0,-1.1) {$p$};
                 % \foreach \Point in {(-1,1), (-1,1), (-2.5,-1),(-0.5,-0.5), (-1,2.5), (1,3)}{
 %   \node at \Point {\textbullet};
%}

    \end{scope}

\end{tikzpicture}
\end{minipage}%
\hspace{-0.5cm}
\begin{minipage}{.49\textwidth}
  \centering

\begin{tikzpicture}[scale=0.6]
    \begin{scope}[shift={(3cm,-5cm)}, fill opacity=0.7]
        \fill[white] \firstcircle;
        \fill[white] \secondcircle;
        \fill[white] \thirdcircle;
        \draw \firstcircle ;
        \draw \secondcircle ;
        \draw \thirdcircle;
               \node [scale=2]   at (-2.5,3.0) {$X$};
              %                \node [scale=2]   at (-2,-1.1) {$l$};
            %                \node [scale=2]   at (2,-1.1) {$m$};
                              \node [scale=2]   at (-4,-3) {$T$};
                                     \node [scale=2]  at (4,-3) {$T'$};
         \node  [scale=1.5]  at (0,0.3) {$r$};
                  \node  [scale=1.5]  at (-1,0.8) {$u$};
                   \node  [scale=1.5]  at (1,0.8) {$v$};
                   \node  [scale=1.3]  at (0,2.2) {$s\!-\!u\!-\!v\!-\!r$};
                           \node  [scale=1.3]  at (-2.4,-1) {$l\!-\!p\!-\!u$};
   \node  [scale=1.3]  at (2.5,-1) {$m\!-\!p\!-\!v$};
                     \node  [scale=1.3]  at (0,-1.3) {$p\!-\!r$};
        %          \node  [scale=2]  at (0,-1.5) {$p$};
                 % \foreach \Point in {(-1,1), (-1,1), (-2.5,-1),(-0.5,-0.5), (-1,2.5), (1,3)}{
 %   \node at \Point {\textbullet};
%}

    \end{scope}

\end{tikzpicture}
\end{minipage}
\caption{Depiction of the sets $T,T' $ with $(T,T')\in\T(n,l,m,p)$, $X\in{\PP_n^{s}}$ and their intersections.}
\label{XandRfigure}
\end{figure}
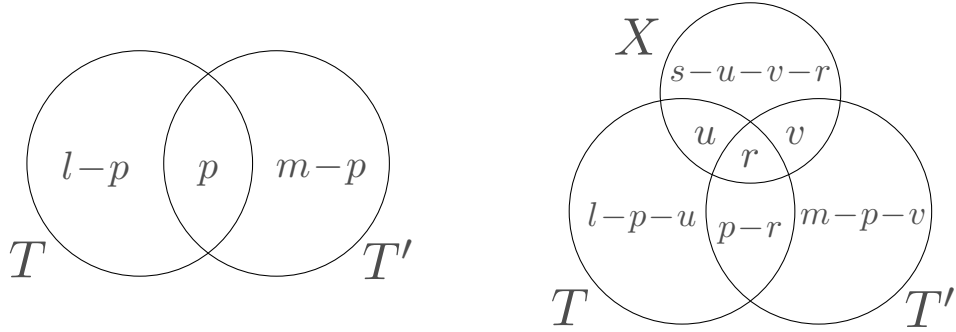

Let the design set $\X$ be $\X=\PP_n^{s}$
  and
 $(T,T')\in\T(n,l,m,p)\,$ both fixed such that
$T\neq T'$ and
$l,m,p$ satisfy
$0\leq p\leq m\leq l\leq n$ and $p<l$.
 Define the quantity
\be
\label{def_R}
\!\!\!\!\!\!\!\!\!\!\!R(n,l,m,p,u,v,r) =
| \left\{X\in\X:\, |X\cap (T\backslash T')|=u,\,
|X\cap (T'\backslash T)|=v,\,
 |X\cap T \cap T'|=r \right\}
 | \,,
\ee
where $u, v, r$ are some nonnegative integers. { $R(n,l,m,p,u,v,r)$ is the number of tests in $\X$ that contain $u$ defective items from $T \setminus T'$, $v$ defective items from $T'\setminus T$ and $r$ defective items from $T \cap T'$. }  Interpretation for the numbers $u, v$ and $ r$  is given on Figure~\ref{XandRfigure} (right).

Observe  that the number $R(n,l,m,p,u,v,r) $
is non-zero only if
\bea
%\label{eq:ineq7}
 0\leq u\leq l-p, \; 0\leq v\leq m-p, \;
0\leq r \leq p\,.
\eea
Joining these restrictions on the parameters $u,v,r$ with the restrictions on
$p,m$ and $l$ in the definition of
the sets
 $\T(n,l,m,p)$, we
obtain  the combined parameter restriction
\be
\label{ae:param_g}
0\leq p\leq m\leq l\leq n, \; p<l, \; 0\leq u\leq l-p, \; 0\leq v\leq m-p, \;
0\leq r \leq p\,.
\ee

%Another consequence of
%(\ref{eq:ineq7})
%is the fact that for any pair $(T,T')$, such that $T\neq T'$, the design set
%$\X$ can be partitioned as follows
%\be
%\label{summationX0}
%\X =
%\bigcup_{r{=}0}^p
%\bigcup_{u{=}0}^{l{-}p}\bigcup_{v{=}0}^{m{-}p}
%\X_{uvr}(T,T')\, .
%\ee

As discussed and proved in Theorem~3.2 in \cite{zhigljavsky2003probabilistic},  formally the design set $\X=\PP_n^{s}$ is balanced. This means the number
$R(n,l,m,p,u,v,r) $ does not depend on the choice of the pair $(T,T')\in
\T(n,l,m,p)$ for any set of integers
$u,v,r,p,m,l$
satisfying
(\ref{ae:param_g}). Moreover, as shown in the next lemma, the number $R(n,l,m,p,u,v,r) $ can be explicitly computed.

\begin{lemma}
\label{th:R}
The design set
$\X\!=\! \PP_n^s$ is balanced  for any  $s \leq n$. For this design set,
and for any set of integers
$u,v,r, p, m, l$
satisfying
{\rm (\ref{ae:param_g})},
we have
\be
\label{eq:RR}
R(n,l,m,p,u,v,r)={{p} \choose {r}} {{l-p} \choose {u}} { {m-p} \choose {v}}  { {n-l-m+p} \choose {s-r-u-v}}
\ee
where the convention
$
\left( \begin{array}{c}
       b \\ a
       \end{array} \right)
=0\; \mbox{{\rm  for $ a<0 \; $ and $a>b$} }
$
may be used for certain values of parameters.
\end{lemma}
For the proof of Lemma~\ref{th:R}, see Theorem~3.2 in \cite{zhigljavsky2003probabilistic}. Lemma~\ref{th:R}   implies, in particular, that the design sets $\X=\PP_n^{\leq s}$ are also balanced for all $1\leq s \leq n$: clearly, a union of disjoint balanced design sets is also a balanced design set.

\section{Derivation of an upper bound for $\gamma$ in a general group testing problem}
\label{sec:3}

\subsection{General test function (\ref{eq:f(X,T)}) and $\X=\PP_n^{s}$}
\label{sec:3.1}

In this section, we consider test functions  $f(\cdot,\cdot)=f_\KK(\cdot,\cdot)$ of the form (\ref{eq:f(X,T)}).  The
following theorem provides a closed-form expression for the
R\'{e}nyi coefficients in this case and represents the major
input into the non-asymptotic expressions of the upper bounds in
specific cases.

\begin{theorem}
\label{th:k_ij}
Let
the test function be defined by
{\rm (\ref{eq:f(X,T)}),}
$0\leq p\leq m\leq l\leq n$, $p<l$, $\X=\PP_n^s$
 and $(T_i,T_j)\in \T(n,l,m,p)$. Then
the value of the R\'{e}nyi coefficient $k_{ij}$ does not depend
on the choice of the pair $(T_i,T_j)\in\T(n,l,m,p)$ and equals
$k_{ij}=
K(\PP_n^s,n,l,m,p),$ where
\be
\label{eq:k_ij-2}
{ K(\PP_n^s,n,l,m,p)}\, &=& \,
 \sum_{r{=}0}^p \sum_{u{=}0}^{m{-}p}
R(n,l,m,p,u,u,r)
\nonumber \\
&{+}&
\sum_{r{=}0}^p \sum_
{u{=}w}^{l{-}p}
\sum_{v{=}u{+}1}^{m{-}p}
R(n,l,m,p,u,v,r)
+ \sum_{r{=}0}^p \sum_{v{=}w}^{m{-}p} \sum_{u{=}v{+}1}^{l{-}p}
R(n,l,m,p,u,v,r)\,  .\;\;\;\;\;
\ee
Here
$ w=\max\{0,\KK-r\} $ and the terms $R(n,l,m,p,u,v,r)$   are as  in \eqref{eq:RR}.
\end{theorem}
The proof of Theorem~\ref{th:k_ij} can be found in Appendix A; it also follows from Theorem~3.3 in \cite{zhigljavsky2003probabilistic}.
Set
\be \label{eq:qq}
q_{\X,n,l,m,p}= {K(\PP_n^s,n,l',m',p)}/{{n \choose s}} \;\;{\rm with}\;\ l'=\max(l,m),m'=\min(l,m)\; {\rm and}\; \X=\PP_n^s\, ,\;\;\;
\ee
where $K(\PP_n^s,n,l,m,p)$ are the
R\'{e}nyi coefficients of \eqref{eq:k_ij-2}; note that using the convention of Lemma~\ref{th:R}, for all $d=0,\ldots,n$ we have
$K(\X,n,d,d,d)=0$ and hence
 $q_{\X,n,d,d,d}= 0$.
Then we have the following theorem.

\begin{theorem}  \label{th:sec6.2}
%{ \it
Let $\T=\PP_n^{\leq d}$ and $\X =\PP_n^{ s}$ where $n\geq 2$, $1\leq d\leq n$, $1\leq s\leq n$. Let $\QQ$ be a distribution satisfying \eqref{Q_restrictions} and let $\RR$ be the uniform distribution on $\X$.  For a fixed $N\geq1$, let $\D_N= \{X_1,\ldots, X_N\} $ be a random $N$-point design with each $X_i \in \D_N $ chosen independently and $\RR$-distributed. Then
\be
\label{eq:binomial_sample1}
\!\!\!\gamma^*(\QQ, \RR,N) =\sum_{\dd =0}^d {\rm Pr}_\UU \{\xi=\dd  \} \min \left\{
1,
\frac{1}{ {{n}\choose {\dd }}} \sum_{m=0}^{d}\, \sum_{p=0}^{\min\{\dd ,m\}}\!
{\textstyle
 {{n}\choose {p\;m-p \; \dd -p\; n-\dd -m+p}}
}
 q^N_{\X,n,\dd ,m,p}
\right\} \,.\;\;\;
\ee
\end{theorem}

The proof of Theorem~\ref{th:sec6.2} is included in the Appendix A; it  is a generalisation of Theorem 6.2 in \cite{zhigljavsky2003probabilistic}.
The following corollary follows from Theorem~\ref{th:sec6.2} and its proof. More specifically, the only adjustment needed in the proof of Theorem~\ref{th:sec6.2} is to set $Q_{N,n,\dd }(\X)=\min \{1,  S_2 \}$, where $S_2$ is defined in the proof.

\begin{corollary} \label{cor:2}
%{ \it
Let $\T=\PP_n^d$ and $\X=\PP_n^s$, where $n\geq 2$, $1\leq d<n$, $1\leq s<n$.
Let $\QQ$ and $\RR$ be uniform distributions on $\T$ and $\X$ respectively.
 For a fixed $N\geq1$, let $\D_N= \{X_1,\ldots, X_N\} $ be a random $N$-point design with each $X_i \in \D_N $ chosen independently and $\RR$-distributed. Then
\be \label{main_cor1}
\gamma^*(\QQ, \RR,N) =\min \left\{
1,
\frac{1}{ {{n}\choose {d}}}\, \sum_{p=0}^{d-1}\!
{\textstyle
 {{n}\choose {p\;d-p \; d-p\; n-2d+p}}
}
 \left( {K(\PP_n^s,n,d,d,p)}/ {{n \choose s}} \right)^N  \
\right\} \,.
\ee
\end{corollary}

\subsection{Additive model} \label{sec:3.2}

In this section we specialize general results of Section~\ref{sec:3.1} to the case of additive model,  where
$
f(X,T) = | X \cap T |
$
so that we can set $\KK=\infty$ in
(\ref{eq:f(X,T)}) and
(\ref{eq:k_ij-2}).
This removes two terms in
(\ref{eq:k_ij-2})
hence
 simplifying  this expression. Furthermore, using \eqref{eq:RR} and the Vandermonde convolution formula, we obtain the following
 statement.

\begin{lemma}
\label{th:ak_ij}
Let $
f(X,T) = | X \cap T |
$, $\X=\PP_n^{s}$ and
$0 \leq p \leq m \leq l \leq n,\; p < l$.
Then $k_{ij}=K(\PP_n^s,n,l,m,p)$ with
%\be
%\label{eq:ak_ij2}
%K(\X,n,l,m,p) =
% \sum_{r=0}^p \sum_{u=0}^{m-p}
%R(n,l,m,p,u,u,r)\, .
%\ee
%For $\X=\PP_n^s$  the formula
%{\rm (\ref{eq:ak_ij2})}
%can be simplified to
\bea
%\label{eq:k_ija0}
K(\PP_n^s,n,l,m,p) =
\sum_{u=0}^{m-p}{{l-p}\choose{u}}
\, {{m-p}\choose {u}} \, {{n-l-m+2p}\choose{s-2u}}\, .
\eea
\end{lemma}

By considering Lemma~\ref{th:ak_ij}, Corollary~\ref{cor:2} and specialising  Theorem~\ref{th:sec6.2} to some specific cases,  we obtain the following corollary.

{

\begin{corollary}\label{Main_binary_corollaries}
Let $
f(X,T) = | X \cap T |
$ and set $n\geq 2$, $1\leq d<n$, $1\leq s<n$ and $\X=\PP_n^s$. For a fixed $N\geq1$, let $\D_N= \{X_1,\ldots, X_N\} $ be a random $N$-point design with each $X_i \in \D_N $ chosen independently and $\RR$-distributed, where $\RR$ is the uniform distribution on  $\X$. We consider the following cases for $\T=\PP_n^d$ and $\QQ$:
\begin{enumerate}
\item Let $\T=\PP_n^d$ and $\QQ$ be the uniform distribution on $\T$.  Then $\gamma^*(\QQ, \RR,N) $ can be obtained from \eqref{main_cor1}
with
\bea
{K(\PP_n^s,n,d,d,p)}=
\sum_{u=0}^{d-p}{{d-p}\choose{u}}^2
\,  \, {{n-2d+2p}\choose{s-2u}}\, \,.
\eea

\item Let $\T=\PP_n^{\leq d}$ and $\QQ$ be a distribution satisfying \eqref{Q_restrictions}. Then $\gamma^*(\QQ, \RR,N) $ can be obtained from \eqref{eq:binomial_sample1} with
\be\label{q_n}
q_{\X,n,\dd ,m,p}=\frac{1}{{n \choose s}} \sum_{u=0}^{m-p}{{\dd -p}\choose{u}} \, {{m-p}\choose {u}} \, {{n-\dd -m+2p}\choose{s-2u}}\  \,.
\ee
%where $\dd '=\max(\dd ,m)$ and $m'=\min(\dd ,m)$.

\item  Let
 $\T=\PP_n^{\leq n}$, $\QQ$ satisfy \eqref{Q_restrictions} and suppose $\UU$ is the $Bin(n,q)$ distribution on $\{0,1,\ldots n \}$. Then from Theorem~\ref{th:sec6.2} we obtain
\bea
\gamma^*(\QQ, \RR,N) =\sum_{\dd =0}^n {n\choose \dd }q^\dd (1-q)^{n-\dd }\min \left\{
1,
\frac{1}{ {{n}\choose {\dd }}} \sum_{m=0}^{n}\, \sum_{p=0}^{\min\{\dd ,m\}}\!
{\textstyle
 {{n}\choose {p\;m-p \; \dd -p\; n-\dd -m+p}}
}
 q^N_{\X,n,\dd ,m,p}
\right\}
\eea
with  $q_{\X,n,\dd ,m,p}$ given in \eqref{q_n}.
\end{enumerate}
\end{corollary}

}

\subsection{Simulation study for the additive model}
\label{sec:3.3}

\begin{figure}[h]
\centering
\begin{minipage}{.5\textwidth}
  \centering
  \includegraphics[width=1\textwidth]{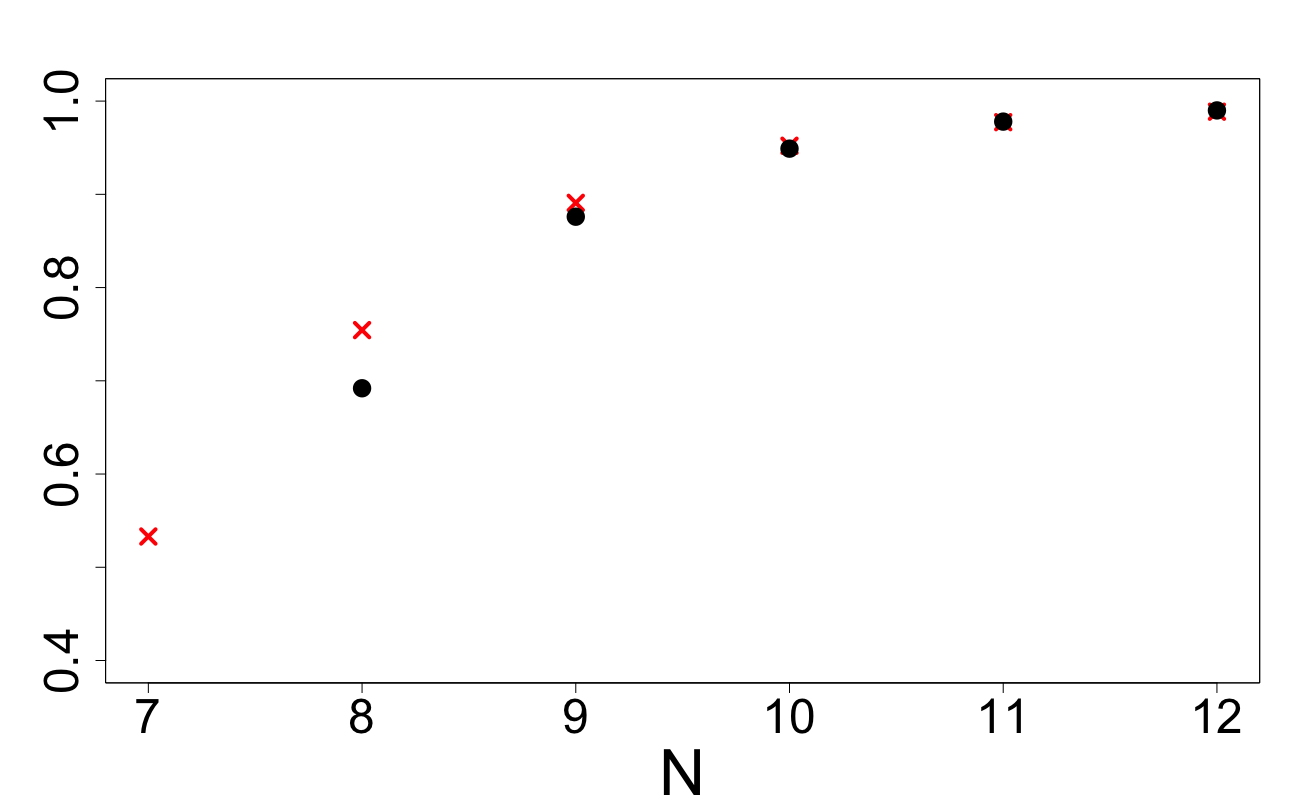}
\caption{Additive model;  $n=20, s=10$.}
\label{alg1_effect1}
\end{minipage}%
\begin{minipage}{.5\textwidth}
  \centering
 \includegraphics[width=1\textwidth]{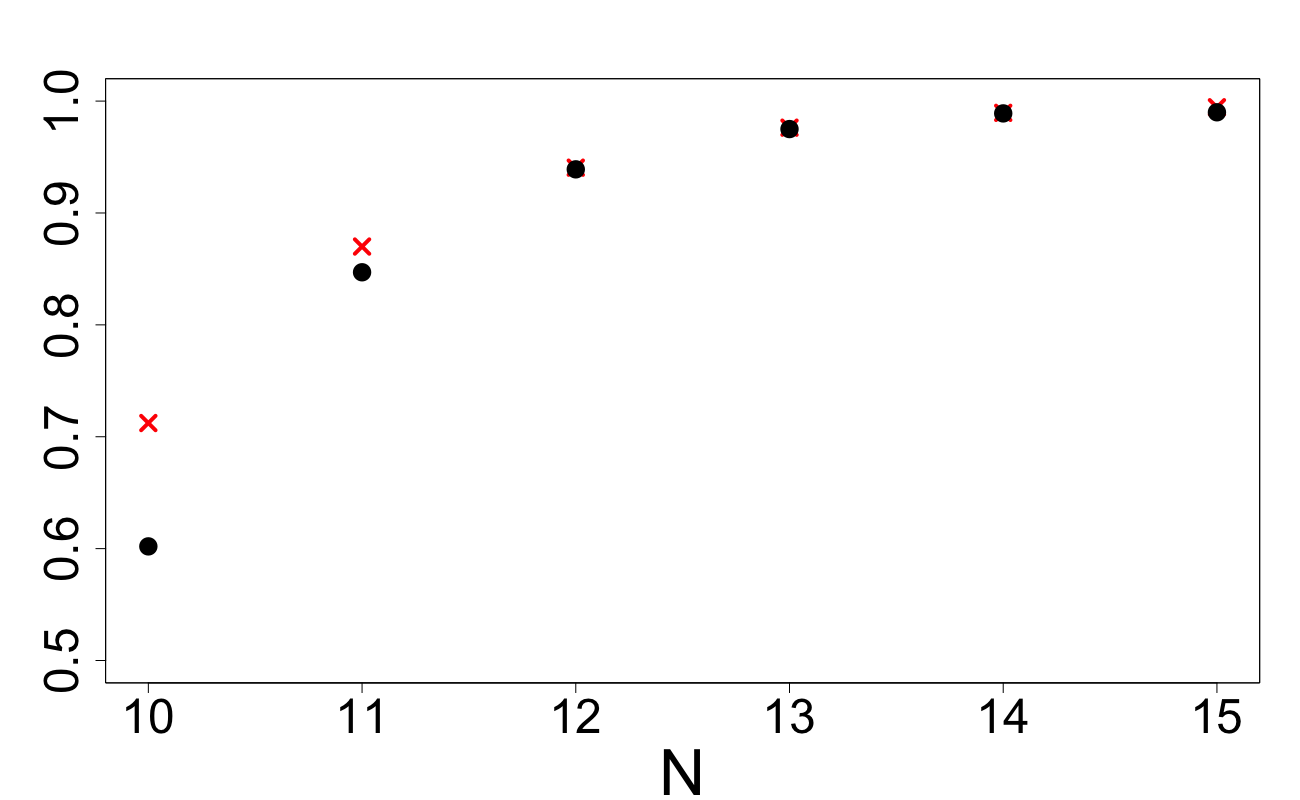}
\caption{Additive model;  $n=50, s=25$.}
\label{alg1_effect2}
\end{minipage}
\end{figure}

In Figures~\ref{alg1_effect1}--\ref{alg1_effect2}, using red crosses we depict the probability ${{\rm Pr}_{\QQ,\RR}\{ T \textrm{ is separated by } \D_{N}  \} }$ as a function of $N$. These values have been obtained via Monte Carlo simulations with $50,000$ repetitions. With the black dots we plot the value of $1-\gamma$ as a function of $N_\gamma$. For these figures, we have set $\T=\PP_n^3$ and chosen $s=n/2 $ based on the asymptotic considerations discussed in the beginning of Section~\ref{sec:as_b_weak}.

In Tables~\ref{table1_add}--\ref{table2_add}, for a given value of $1-\gamma^*$ we tabulate the value of $1-\gamma$ for the additive group testing model, where $\T=\PP_n^3$, $\X=\PP_n^s$ and $\QQ$ and $\RR$ are uniform on $\T$ and $\X$ respectively. The values have been obtained via Monte Carlo simulations. When considering the inverse problem discussed in \eqref{inverse problem}, we also include the explicit upper bounds $N_\gamma$ and the value of $N_{\gamma^*}$ obtained via Monte Carlo for different values of $n, s$ and $\gamma^*$. In all Monte Carlo simulations, we have used $50,000$ repetitions. Tables~\ref{table1_add}--\ref{table2_add} and Figures~\ref{alg1_effect1}--\ref{alg1_effect2} demonstrate that when $\gamma^*$ is small,  the union bound used in the proof of Theorem~\ref{th:sec2.2} appears very sharp since the values of $1-\gamma $ and $1-\gamma^*$ almost coincide.

\begin{table}[h]

\begin{center}
\begin{tabular}{|c||c|c|c|c|c|c|c|c|c|c|c|c|}
\cline{2-13}
\multicolumn{1}{c|}{} &\multicolumn{3}{c|}{$n=20$}&\multicolumn{3}{c|}{$n=50$}&\multicolumn{3}{c|}{$n=100$}&\multicolumn{3}{c|}{$n=150$}  \\
\cline{1-13}
$\lambda$&$N_{\gamma^*} $ & $N_\gamma$ & $1-\gamma $ &$N_{\gamma^*}$ & $N_\gamma$  & $1-\gamma$&$N_{\gamma^*}$ & $N_\gamma$  & $1-\gamma$& $N_{\gamma^*}$ &$N_\gamma$  & $1-\gamma$ \\
\hline
\hline
0.10  & 31& 34  & 0.96&    38 &  40& 0.96 &42  & 44  & 0.97  & 42  & 46 &0.96  \\
0.20 & 16 & 17   &0.96   & 19 &  21 & 0.97 &  21 &23   &  0.96 & 23 &  24 & 0.96  \\
 0.30&11 & 12   & 0.97   & 14 &    15 &  0.97 & 14 & 16   & 0.97  & 17 & 18 & 0.97  \\
0.40 & 9 & 11   &  0.98 &  11 &   13 &  0.98 & 13 & 15 &  0.98 & 14 &   16 & 0.98 \\
0.50 &8 & 11 &   0.98 & 11 &   13 &  0.98 &  12 & 14 &   0.98  & 13 &  15 &0.98 \\
\hline
% etc. ...
\end{tabular}
\end{center}
\caption{Additive model with $\gamma^*=0.05$, $d=3$
$s =\lceil \lambda n \rceil$,
various  $n$ and $\lambda$.}
\label{table1_add}
\end{table}

\begin{table}[h]

\begin{center}
\begin{tabular}{|c||c|c|c|c|c|c|c|c|c|c|c|c|}
\cline{2-13}
\multicolumn{1}{c|}{} &\multicolumn{3}{c|}{$n=20$}&\multicolumn{3}{c|}{$n=50$}&\multicolumn{3}{c|}{$n=100$}&\multicolumn{3}{c|}{$n=150$}  \\
\cline{1-13}
$\lambda$&$N_{\gamma^*} $ & $N_\gamma$ & $1-\gamma $ &$N_{\gamma^*}$ & $N_\gamma$  & $1-\gamma$&$N_{\gamma^*}$ & $N_\gamma$  & $1-\gamma$& $N_{\gamma^*}$ &$N_\gamma$  & $1-\gamma$ \\
\hline
\hline
0.10  & 28  & 30 & 0.92 & 34 &  36 & 0.93  &  38 &  40 & 0.93 &41 & 43&0.94 \\
0.20 &  15 &  16 & 0.93 &  17 & 19 & 0.93  & 20 &   21 &  0.93 & 21& 22 & 0.93  \\
 0.30& 9 &   11  &  0.94 &12  & 14 & 0.94  &14  &  15 & 0.95  & 15&  17 &0.95 \\
0.40 &  8 &  10  &  0.95 &10 &12 & 0.95  &  12&  14 &  0.95 & 13& 15 & 0.96 \\
0.50 &  8 &  10  &  0.96 &  10 & 12 &  0.97   & 12 & 14  & 0.96 &12  &14  &0.97 \\
\hline
% etc. ...
\end{tabular}
\end{center}
\caption{ Additive model with $\gamma^*=0.1$, $d=3$
$s =\lceil \lambda n \rceil$,
various  $n$ and $\lambda$.}
\label{table2_add}
\end{table}

\subsection{Extension for $\X\neq \PP_n^{s}$}\label{general_x}

{

In this  section we demonstrate how the key results of the Sections~\ref{sec:3.1} and \ref{sec:3.2} can be easily modified  for  the case when $\X = \cup_s \PP_n^{s}$, where the union is taken over any subset of $\{0,1, \ldots, n\}$,  and for a distribution $\RR$  that is not necessarily uniform on $\X$.

Let $\X =\PP_n^{ \leq n}$, $\VV$ be a probability distribution on $\{0,1, \ldots, n\}$ and $\zeta$ be a $\VV$-distributed random variable on
$\{0,1, \ldots, n\}$. The distribution $\RR$ depends on $\VV$ in the following way:
for a $\RR$-distributed  random test $X \in \X$  we have
 \be\label{R_restrictions}
 \mbox{Pr}_{\RR}\{ |X| = {s}  \} =  \mbox{Pr}_\VV\{ \zeta={s}  \}, \,\,\  \mbox{Pr}_{\RR}\{ X=x \,|\, |X|={s}\} = 1/{{n \choose s }} \,\,\,\, \mbox{$\forall x\in  \PP_n^{s} $, else 0} \,.
 \ee
These two requirements mean that for all $s \in \{0,1, \ldots, n\}$  and $\textsf{X} \in  \PP_n^{s}  $ we have
\bea
 \mbox{Pr}_{\RR}\{ X= \textsf{X}
    \} ={\mbox{Pr}_\VV\{ \zeta=s  \}}/{ {n \choose s}}  \,.
 \eea
Note that in the case of Bernoulli design, when each item is included into a group of items with  probability  $p$, $\VV$ is Bin($n,p$), the Binomial distribution  with parameters $n$ and $p$.

For a general test function $f(X,T)$  we introduce the probability
$$
p_{ijs}  = \mbox{Pr}\{f(X,T_i) = f(X,T_j) \,|\,|X|=s  \}.
$$
By conditioning on $s$, we obtain
$
p_{ij} = \sum_{s=0}^n p_{ijs}\mbox{Pr}_\VV\{ \zeta={s}  \}.
$
In view of the conditional uniformity of $\RR$, which is the second condition in  \eqref{R_restrictions},    the
probabilities $p_{ijs}$
can be written as
\bea
p_{ijs} = {k_{ijs}}/{| \PP_n^{s} |} = {k_{ijs}}/{{{n\choose s}}}
\eea
where
$k_{ijs} = k(T_i,T_j,s)\,$
is the number of $X\in  \PP_n^{s} $ such that
$
f(X,T_i)\, =\, f(X,T_j);\;$
that is,
 \bea
%\label{eq:3}
k_{ijs} = \left| \{ X\in { \PP_n^{s} }:\; f(X,T_i) = f(X,T_j)\} \right| \;\;\;\;\;
\mbox{\rm for $\;\;T_i,T_j  \in {\T}$}\, .
\eea
From this, we obtain
\bea
p_{ij} = \sum_{s=0}^n{k_{ijs}}\mbox{Pr}_\VV\{ \zeta={s}\} / {{{n\choose s}}} \, .
\eea
Set
\bea
q_{\X,n,l,m,p;\,\VV}=  \sum_{s =0}^n \frac{K(\PP_n^s,n,l',m',p)} {{n \choose s}}\mbox{Pr}_\VV\{ \zeta={s}  \} \;\;{\rm with}\;\; l'=\max(l,m),m'=\min(l,m)\,.
\eea
Then all results of the previous sections established for the case $\X=\PP_n^{ s}$ can be can be extended for the group testing problems with $\X =\cup_{s}\PP_n^{ s}$ by replacing  $q_{\X,n,l,m,p} $ of \eqref{eq:qq} with $q_{\X,n,l,m,p;\,\VV}$.

}

\section{Group testing for the binary model}
\label{sec:4}

\subsection{A general result and its specialization to particular cases }
\label{sec:4.1}

In the binary group testing,  we have  $\KK=1$ in \eqref{eq:f(X,T)} and thus the test function is
\be
\label{eq:f(X,T)b}
f(X,T)=f_1(X,T)=\left\{\begin{array}{ll}
                 0  & \; \mbox{ if } \; |X\cap T|=\emptyset, \\
                 1  &  \mbox{ otherwise.}
              \end{array}\right.
\ee

%\begin{theorem}
%\label{th:k_ijb}
%Let
%the test function be
%{\rm (\ref{eq:f(X,T)b}), }
%$0\leq p\leq m\leq l\leq n$, $p<l$ and $\X$ be a balanced design set.
%Then $k_{ij}=K(\X,n,l,m,p)$ with
%\be
%\label{eq:k_ijb}
%K(\X,n,l,m,p) =
%|\X|-\left(\sum_{u=1}^{l-p}R(n,l,m,p,u,0,0)+\sum_{v=1}^{m-p}
%R(n,l,m,p,0,v,0)\right)\, .
%\ee
%\end{theorem}
%
%
%
%
%\begin{corollary}
% \label{cor:upper1}
%Assume either
%$\T=\PP_n^d$ or
%$\T=\PP_n^{\leq d}$, $d \leq n,\;$
%and let the design set $\X$ be balanced.
% For a fixed $N\geq1$, let $\D_N= \{X_1,\ldots, X_N\} $ be a random $N$-point design with each $X_i \in \D_N $ chosen independently and $\RR$ distributed. Then
%\be
%\gamma^*(\QQ, \RR,N)&= & \frac{2}{|\T|}\;\suml_{l,m}\suml_{{0\leq p\leq m}\atop{p<l}}
%Q(n,l,m,p) \times \Biggr.
%\nonumber\\
%\label{eq:N*gam}
%&&\Biggl. \times
%\left( 1{-}\frac{1}{|\X|} \left(\suml_{u=1}^{l-p}R(n,l,m,p,u,0,0){+}\suml_{v=1}^{m-p}
%R(n,l,m,p,0,v,0)\right) \right)^N \,,
%\ee
%where the first summation in the above formulae
%is over $0\leq m\leq l\leq d$ for
%the case $\T=\PP_n^{\leq d}$ and $l=m=d$; that is, the first summation
%disappears in
%{\rm
%(\ref{eq:N*gam}),} in the case $\T=\PP_n^d$.
%\end{corollary}
%
%The binary group testing model is the most popular in literature.
%In view of its importance we  provide below
%several specific existence theorems
%for this model.
%Of course, all these theorems are easy corollaries of
%(\ref{eq:k_ijb0}), (\ref{eq:N*gam})
%and Theorems~\ref{th:q}, \ref{th:R}.
%

\begin{theorem}
\label{th:k_ijb2}
Let
the test function be
{\rm (\ref{eq:f(X,T)b}),}
$0\leq p\leq m\leq l\leq n$, $p<l$, $\X=\PP_n^s$
 and $(T_i,T_j)\in \T(n,l,m,p)$. Then
the value of the R\'{e}nyi coefficient $k_{ij}$ does not depend
on the choice of the pair $(T_i,T_j)\in\T(n,l,m,p)$ and equals
$k_{ij}=
K(\PP_n^s,n,l,m,p),$ where
\be
\label{eq:k_ijb0}
K(\PP_n^s,n,l,m,p) =
{{n} \choose {s}}
-{{n-l} \choose {s}}
-{{n-m} \choose {s}}+2 {{n-l-m+p} \choose {s}}\, .
\ee
\end{theorem}
The proof of Theorem~\ref{th:k_ijb2} can be obtained from \cite{zhigljavsky2003probabilistic} and is included in Appendix A for completeness.

%\boldmath

{

\begin{corollary}\label{simple_corollary}
Let
the test function be
{\rm (\ref{eq:f(X,T)b})} and set $n\geq 2$, $1\leq d<n$, $1\leq s<n$ and $\X=\PP_n^s$. For a fixed $N\geq1$, let $\D_N= \{X_1,\ldots, X_N\} $ be a random $N$-point design with each $X_i \in \D_N $ chosen independently and $\RR$-distributed, where $\RR$ is the uniform distribution on  $\X$. We consider the following cases for $\T=\PP_n^d$ and $\QQ$:
\begin{enumerate}
\item Let $\T=\PP_n^d$ and $\QQ$ be the uniform distribution on $\T$.  For a fixed $N\geq1$, let $\D_N= \{X_1,\ldots, X_N\} $ be a random $N$-point design with each $X_i \in \D_N $  independent and $\RR$-distributed. Then $\gamma^*(\QQ, \RR,N) $ can be obtained from \eqref{main_cor1}
with
\bea
{K(\PP_n^s,n,d,d,p)} =
{{n \choose s}}-2{{{n-d}\choose {s}} +2
{{n-2d+p} \choose {s}}} \,.
\eea

\item Let $\T=\PP_n^{\leq d}$ and $\QQ$ be a distribution satisfying \eqref{Q_restrictions}. Then $\gamma^*(\QQ, \RR,N) $ can be obtained from \eqref{eq:binomial_sample1} with
\be\label{q_n2}
q_{\X,n,\dd ,m,p}=1{-}\left[{\a{n{-}\dd }{s}
+\a{n{-}m}{s}{-}2\a{n{-}\dd {-}m{+}p}{s}}\right] \big/ {\a{n}{s}} \,.
\ee

\item Let $\T=\PP_n^{\leq n}$,  $\QQ$ be a distribution satisfying \eqref{Q_restrictions} and suppose $\UU$ is the $Bin(n,q)$ distribution on $\{0,1,\ldots n \}$ for some $q>0$. Then application of Theorem~\ref{th:sec6.2} provides
\bea
\gamma^*(\QQ, \RR,N) =\sum_{\dd =0}^n {n\choose \dd }q^\dd (1-q)^{n-\dd }\min \left\{
1,
\frac{1}{ {{n}\choose {\dd }}} \sum_{m=0}^{n}\, \sum_{p=0}^{\min\{\dd ,m\}}\!
{\textstyle
 {{n}\choose {p\;m-p \; \dd -p\; n-\dd -m+p}}
}
 q^N_{\X,n,\dd ,m,p}
\right\}
\eea
with  $q_{\X,n,\dd ,m,p}$ obtained from \eqref{q_n2}.

\end{enumerate}
\end{corollary}

}

In Table~\ref{exactly_3}, using the results of part one of Corollary~\ref{simple_corollary} we consider the inverse problem discussed in \eqref{inverse problem} and tabulate the value of $N_\gamma$ supposing $\T=\PP_n^3$  for different values of $s$ and $n$. In Table~\ref{binomial_3}, using the results of part three of Corollary~\ref{simple_corollary} we tabulate the value of $N_\gamma$ supposing $\UU$ is the $Bin(n,3/n)$ distribution. In distribution $\UU$, the probability of success has been set to $3/n$ so that each target $T\in \T$ will have three elements on average to compare with the results of Table~\ref{exactly_3}. We see the binomial sample problem requires significantly more tests to locate the defective items with high probability than the case of exactly $d$ defectives.

\begin{table}[!h]
\begin{center}
\begin{tabular}{|P{7mm}|}
  \multicolumn{1}{c}{  } \\
 \hline
  \multicolumn{1}{|c|}{ $\lambda$ } \\ \hline
  0.10 \\
   0.15 \\
  0.20 \\
  0.25   \\
    0.30   \\
  \hline
\end{tabular}
\begin{tabular}{ |P{12mm} |P{12mm}|P{13mm}| }
  \hline
  \multicolumn{3}{|c|}{ $\gamma=0.01$ } \\ \hline
$n=20$ & $n=50$& $n=100$   \\ \hline
47&58 &64 \\
37& 47& 50\\
33& 40& 44\\
32& 39& 43\\
34& 40& 44\\
  \hline
\end{tabular}
\begin{tabular}{ |P{12mm} |P{12mm}|P{13mm}|}
  \hline
  \multicolumn{3}{|c|}{ $\gamma=0.05$ } \\ \hline
$n=20$ & $n=50$& $n=100$   \\ \hline
38& 48&54 \\
30& 39& 42\\
27& 33& 37\\
26& 33& 36\\
28& 34& 38\\
  \hline\end{tabular}
\end{center}
\caption{Values of $N_\gamma$  for binary model with $d=3$, $s =\lceil \lambda n \rceil$ for
various  $n$ and $\lambda$.}
\label{exactly_3}
\end{table}

\begin{table}[!h]
\begin{center}
\begin{tabular}{|P{7mm}|}
  \multicolumn{1}{c}{  } \\
 \hline
  \multicolumn{1}{|c|}{ $\lambda$ } \\ \hline
  0.10 \\
   0.15 \\
  0.20 \\
  0.25   \\
  0.30  \\
  \hline
\end{tabular}
\begin{tabular}{ |P{12mm} |P{12mm}|P{13mm}| }
  \hline
  \multicolumn{3}{|c|}{ $\gamma=0.01$ } \\ \hline
$n=20$ & $n=50$& $n=100$   \\ \hline
 90& 119&142 \\
 84& 117&184 \\
 105& 187& 410\\
 166&283 & 731\\
 316   &547   &1334  \\
  \hline
\end{tabular}
\begin{tabular}{ |P{12mm} |P{12mm}|P{13mm}|}
  \hline
  \multicolumn{3}{|c|}{ $\gamma=0.05$ } \\ \hline
$n=20$ & $n=50$& $n=100$   \\ \hline
71 & 95&113 \\
 63&91 &154 \\
 70&129 & 242\\
 101& 186& 380\\
   170  & 330  & 604 \\
  \hline\end{tabular}
\end{center}
\caption{Values of $N_\gamma$ for binary model with $\UU$ the $Bin(n,3/n)$ distribution, $s =\lceil \lambda n \rceil$ for
various  $n$ and $\lambda$.}
\label{binomial_3}
\end{table}

The results below will address the scenario of Bernoulli designs.
In the following corollaries we set $\X =\PP_n^{ \leq n}$ and $\VV$ is the $Bin(n,\kappa)$ distribution for some $0<\kappa<1$.
The discussion of Section~\ref{general_x} results in the following.

{

\begin{corollary}\label{cor:N*1}
Let
the test function be
{\rm (\ref{eq:f(X,T)b})} and set $n\geq 2$, $1\leq d<n$, $1\leq s<n$. Let $\X =\PP_n^{ \leq n}$, $\RR$ be a distribution satisfying the constraints \eqref{R_restrictions} and suppose $\VV$ is the $Bin(n,\kappa)$ distribution on $ \{0,1,\ldots n \}$. Let $\D_N= \{X_1,\ldots, X_N\} $ be a random  design with each $X_i \in \D_N $ chosen independently and $\RR$-distributed. We consider the following cases for $\T=\PP_n^d$ and $\QQ$:
\begin{enumerate}

\item Let $\T=\PP_n^{d}$ and $\QQ$ be the uniform distribution on $\T$.  Then $\gamma^*(\QQ, \RR,N) $ can be obtained from \eqref{main_cor1} by replacing
$
{K(\PP_n^s,n,d,d,p)}/ {{n \choose s}}
$
  with
\bea
\sum_{s=0}^n\frac{K(\PP_n^s,n,d,m,p)} {{n \choose s}}{\rm Pr}_\VV\{ S={s}  \}= 1-2\suml_{s=0}^n \left({n-d \choose s} - {n-2d+p \choose s}\right) \kappa^s(1-\kappa)^{n-s} \,.\\
\eea

\item Let $\T=\PP_n^{\leq n}$,  $\QQ$ be a distribution satisfying the constraint \eqref{Q_restrictions} and  suppose $\UU$ is the $Bin(n,q)$ distribution on $ \{0,1,\ldots n \}$. Then from \eqref{eq:binomial_sample1} we obtain
\bea
\label{eq:binomial_sample}
\!\!\gamma^*(\QQ, \RR,N) \!=\!\sum_{\dd =0}^n {n\choose \dd }q^\dd (1\!-\!q)^{n-\dd }\min \left\{
1,
\frac{1}{ {{n}\choose {\dd }}} \sum_{m=0}^{n}\, \sum_{p=0}^{\min\{\dd ,m\}}\!
{\textstyle
 {{n}\choose {p\;m-p \; \dd -p\; n-\dd -m+p}}
}
 q^N_{\X,n,\dd ,m,p,\kappa}
\right\} ,\;\;\;\;
\eea
where
\bea
q_{\X,n,\dd ,m,p,\kappa}=1{-}\suml_{s=0}^n \left({{{n{-}\dd }\choose {s}}
+{{n{-}m}\choose{s}}{-}2{{n{-}\dd {-}m{+}p}\choose{s}}} \right)\kappa^s(1-\kappa)^{n-s} \,.
\eea
\end{enumerate}
\end{corollary}

}

 In Table~\ref{compare_binom1}, using the results of part one of Corollary~\ref{simple_corollary} we tabulate the value of $N_\gamma$ supposing $\T=\PP_n^d$ with $d=3$ for different values of $s$ and $n$. This table considers more choices for $s$ when compared to Table~\ref{exactly_3}. In Table~\ref{compare_binom2}, we tabulate the value of $N_\gamma$ obtained via part one of Corollary~\ref{cor:N*1} supposing $\VV$ is the $Bin(n,\lceil\lambda n \rceil /n)$ distribution. The probability parameter has been set to $\lceil\lambda n \rceil /n$ such that each $X_i$ in $\D_N=\{X_1,\ldots,X_N\}$ will have $\lceil\lambda n \rceil$  elements on average to compare with the results of Table~\ref{compare_binom1}. The results of these tables indicate it is preferable to have a design with constant-row-weight rather than including each item in a test with some fixed probability (at least for choices of $s$ of interest).

%(note that using the convention
%$K(\X,n,d,d,d)=0$
%of Theorem \ref{th:exist_Kch},
%we have $q_{\X,n,d,d,d}= 0 $
%for all $d=0,\ldots,n$).
%

\begin{table}[!h]
\begin{center}
\begin{tabular}{|P{7mm}|}
  \multicolumn{1}{c}{  } \\
 \hline
  \multicolumn{1}{|c|}{ $\lambda$ } \\ \hline
  0.05 \\
   0.10 \\
  0.15 \\
  0.20   \\
  0.25 \\
  0.30  \\
  0.35   \\
  0.40  \\
  0.45 \\
        0.50 \\
  \hline
\end{tabular}
\begin{tabular}{ |P{12mm} |P{12mm}|P{12mm}|P{14mm}| }
  \hline
  \multicolumn{4}{|c|}{ $\gamma=0.01$ } \\ \hline
$n=10$ & $n=20$& $n=50$  & $n=100$ \\ \hline
 35&82 &86 &112\\
35 &47 &58 &64\\
  25& 33&43  &48\\
 25 &33  &40 &44\\
 27 &32 &39 &43\\
 27 &34 &40 &44\\
37   &43  &45  &48\\
 37 & 43 &50 &54\\
62 & 52&62 &64\\
 62 & 66& 73&79\\
  \hline
\end{tabular}
\begin{tabular}{ |P{12mm} |P{12mm}|P{12mm}|P{14mm}| }
  \hline
  \multicolumn{4}{|c|}{ $\gamma=0.05$ } \\ \hline
$n=10$ & $n=20$& $n=50$  & $n=100$ \\ \hline
28 & 66& 72& 94\\
 28& 38& 48&54\\
  20&27 & 36 & 41\\
  20& 27 &33 &37\\
  22& 26& 33&36\\
  22&28 &34 &38\\
   29& 35 & 38 &41\\
  29& 35 &42 &46\\
 51&43 & 52&55\\
  51&55 & 63&69\\
  \hline\end{tabular}
\end{center}
\caption{Values of $N_\gamma$ for binary model with $d=3$, $s =\lceil \lambda n \rceil$ for
various  $n$ and $\lambda$.}
\label{compare_binom1}
\end{table}

\begin{table}[!h]
\begin{center}
\begin{tabular}{|P{7mm}|}
  \multicolumn{1}{c}{  } \\
 \hline
  \multicolumn{1}{|c|}{ $\lambda$ } \\ \hline
  0.05 \\
   0.10 \\
  0.15 \\
  0.20   \\
  0.25 \\
  0.30  \\
  0.35   \\
  0.40  \\
  0.45 \\
        0.50 \\
  \hline
\end{tabular}
\begin{tabular}{ |P{12mm} |P{12mm}|P{12mm}|P{14mm}| }
  \hline
  \multicolumn{4}{|c|}{ $\gamma=0.01$ } \\ \hline
$n=10$ & $n=20$& $n=50$  & $n=100$ \\ \hline
 49& 96 &92 &115\\
49 & 55& 61&66\\
  34& 38&  46&49\\
  34&  38& 42&45\\
  34& 37& 41&44\\
  34& 38& 42&45\\
   41& 46 &  46&49\\
  41&  46& 51&55\\
 58&53 & 62&64\\
  58 &65 &73 &79\\
  \hline
\end{tabular}
\begin{tabular}{ |P{12mm} |P{12mm}|P{12mm}|P{14mm}| }
  \hline
  \multicolumn{4}{|c|}{ $\gamma=0.05$ } \\ \hline
$n=10$ & $n=20$& $n=50$  & $n=100$ \\ \hline
 39& 78& 76& 97\\
 39&44 & 51&56\\
  27&31 & 38 &42\\
  27& 31 & 35&38 \\
  27& 30&34 &37\\
 27 & 31&35 &38\\
   33& 37 & 39 &41\\
  33& 37 & 43&47\\
 47&44 & 52&55\\
  47&54 & 62&69\\
  \hline\end{tabular}
\end{center}
\caption{Values of $N_\gamma$ for binary model with $\VV$ the $Bin(n,\lceil \lambda n \rceil /n)$ distribution for
various  $n$ and $\lambda$.}
\label{compare_binom2}
\end{table}

\subsection{Simulation study}
\label{sec:4.2}

 In Tables~\ref{table1}--\ref{table3}, for a given value of $1-\gamma^*$ we tabulate the value of $1-\gamma$ for the binary group testing model, where $\T=\PP_n^d$, $\X=\PP_n^s$ and $\QQ$ and $\RR$ are uniform on $\T$ and $\X$ respectively. Similarly to Tables~\ref{table1_add}--\ref{table2_add}, we also include the explicit upper bounds $N_\gamma$ and the value of $N_{\gamma^*}$ obtained via Monte Carlo methods with $50,000$ trials for different values of $n, s$ and $\gamma^*$. We see once again, that for small values of $\gamma^*$, the union bound used in Theorem~\ref{th:sec2.2} appears very sharp.

\begin{figure}[!h]
\centering
\begin{minipage}{.5\textwidth}
 \centering
 \includegraphics[width=1\textwidth]{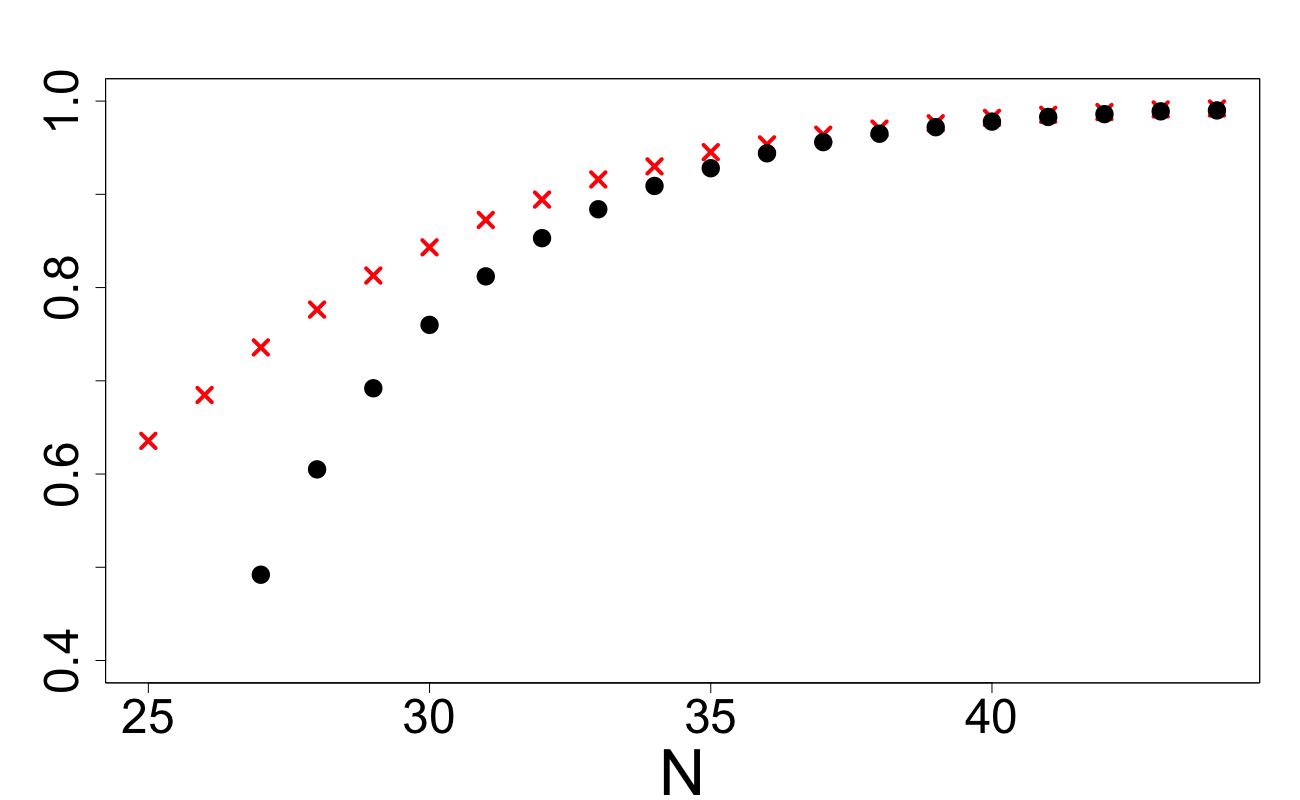}
   \captionsetup{width=.9\linewidth}
\caption{ Binary model: $\gamma$ vs $\gamma^*$ for $n=100$ and $d=3$.}
\label{Bonf_n_50}
\end{minipage}%
\begin{minipage}{.5\textwidth}
 \centering
\includegraphics[width=1\textwidth]{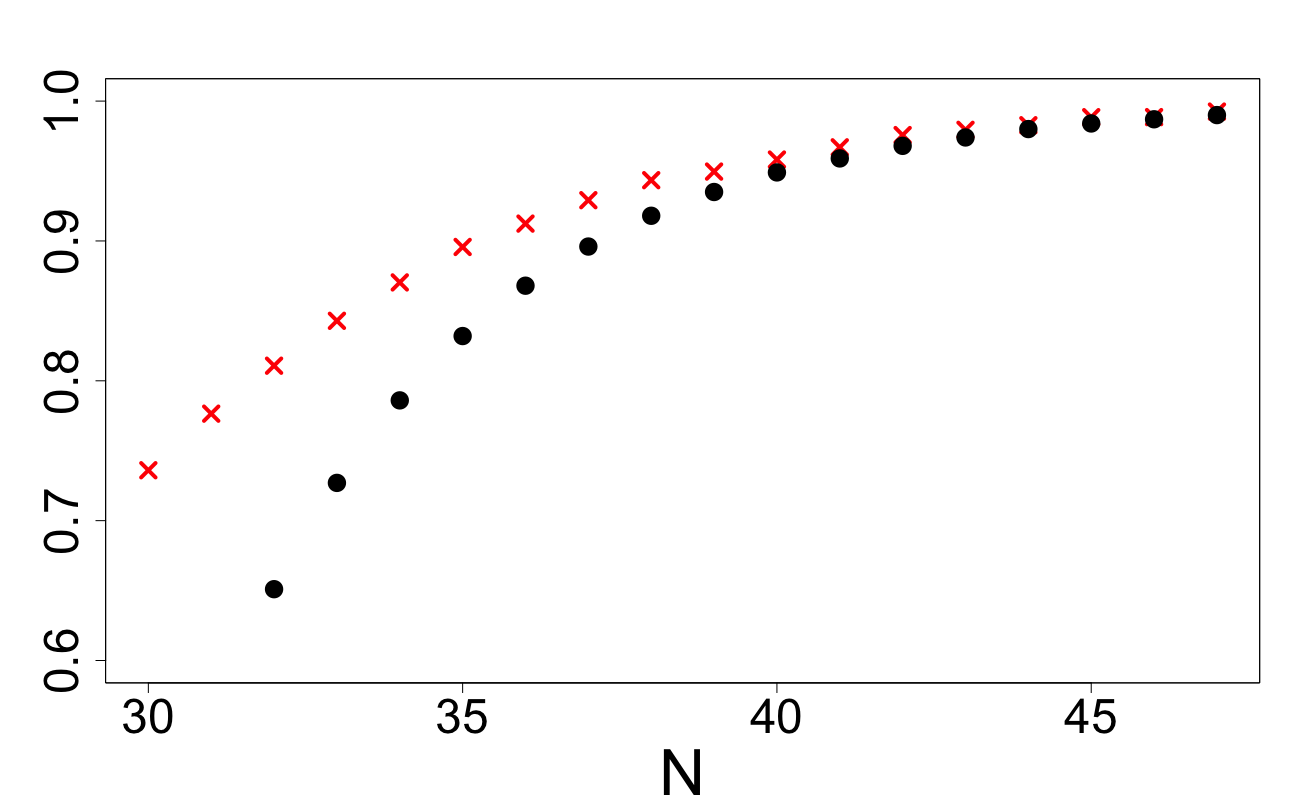}
  \captionsetup{width=.9\linewidth}
\caption{ Binary model: $\gamma$ vs $\gamma^*$ for $n=200$ and $d=3$.}
\end{minipage}
\end{figure}

\begin{figure}[!h]
\centering
\begin{minipage}{.5\textwidth}
 \centering
 \includegraphics[width=1\textwidth]{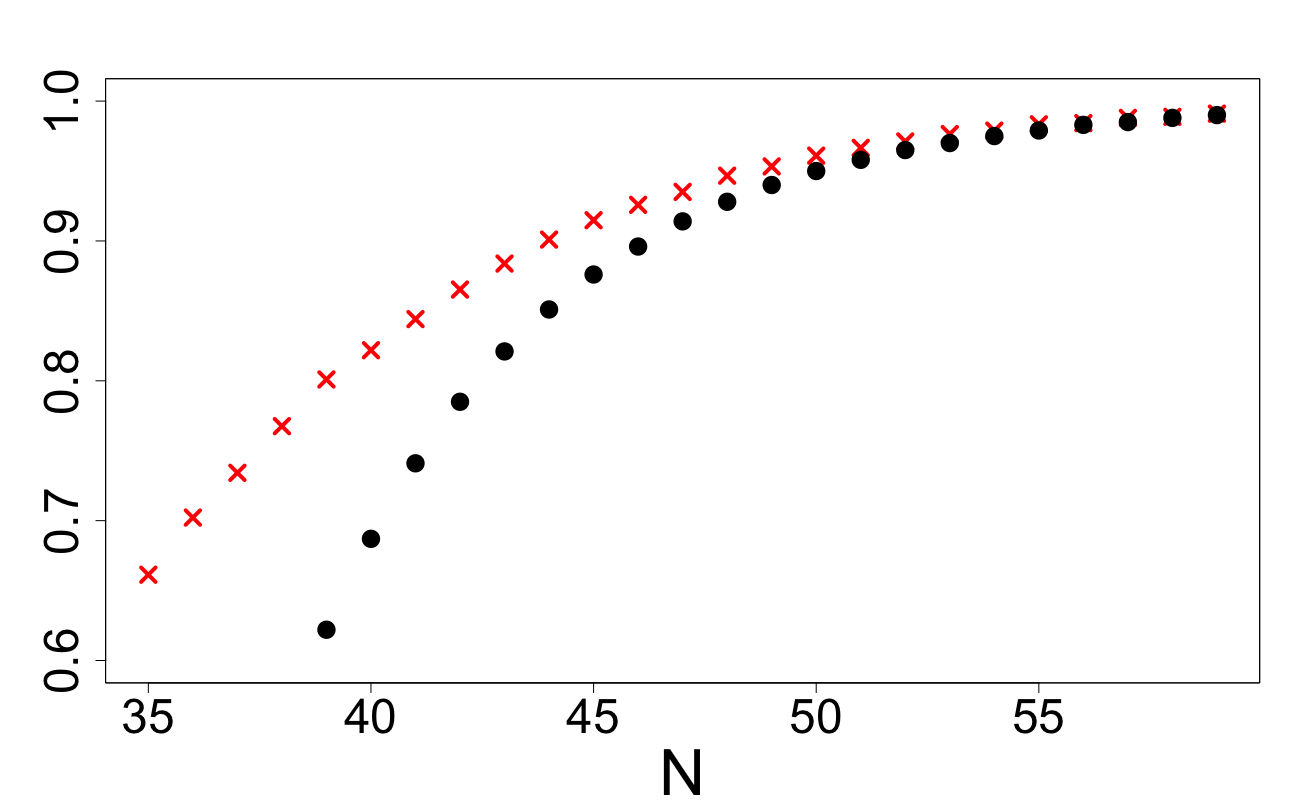}
   \captionsetup{width=.9\linewidth}
\caption{ Binary model: $\gamma$ vs $\gamma^*$ for $n=100$ and $d=4$.}

\end{minipage}%
\begin{minipage}{.5\textwidth}
 \centering
\includegraphics[width=1\textwidth]{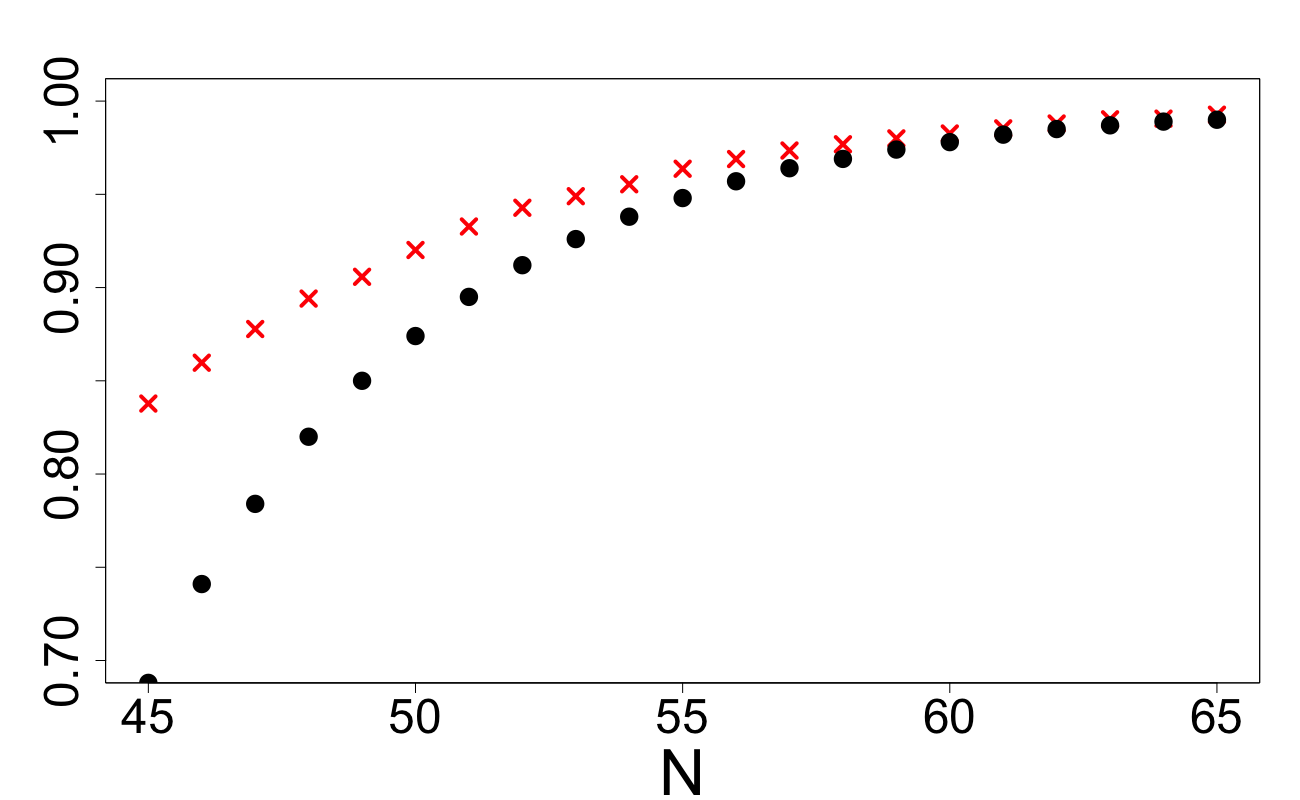}
  \captionsetup{width=.9\linewidth}
\caption{ Binary model: $\gamma$ vs $\gamma^*$ for $n=200$ and $d=4$.}
\label{Bonf_n_200}
\end{minipage}
\end{figure}

In Figures~\ref{Bonf_n_50}--\ref{Bonf_n_200}, using red crosses we depict the probability ${{\rm Pr}_{\QQ,\RR}\{ T \textrm{ is separated by } \D_{N}  \} }$ as a function of $N$ obtained with $50,000$ Monte Carlo simulations. With the black dots we plot the value of $1-\gamma$ as a function of $N_\gamma$. For these figures, we have chosen $s=\lfloor (1- 2^{-1/d} )n \rfloor $ based on the asymptotic considerations discussed in the beginning of Section~\ref{sec:as_b_weak}.

\begin{table}[!h]

\begin{center}
\begin{tabular}{|c||c|c|c|c|c|c|c|c|c|c|c|c|}
\cline{2-13}
\multicolumn{1}{c|}{} &\multicolumn{3}{c|}{$n=20$}&\multicolumn{3}{c|}{$n=50$}&\multicolumn{3}{c|}{$n=100$}&\multicolumn{3}{c|}{$n=200$}  \\
\cline{1-13}
$\lambda$&$N_{\gamma^*} $ & $N_\gamma$ & $1-\gamma $ &$N_{\gamma^*}$ & $N_\gamma$  & $1-\gamma$&$N_{\gamma^*}$ & $N_\gamma$  & $1-\gamma$& $N_{\gamma^*}$ &$N_\gamma$  & $1-\gamma$ \\
\hline
\hline
0.10  &  36  &  38  & 0.96  & 44 &48  &0.96 & 49 & 54 & 0.96& 55&59&0.97  \\
0.20 &  25 &   27  &  0.96 &   30 & 33 & 0.96  &33  & 37&0.96 & 38& 41& 0.96 \\
 0.30&  29 &   31 &  0.96 & 33 &35  & 0.96  & 35 & 38 & 0.97& 39& 41 &0.96 \\
0.40 &33   &   35&  0.96 &  37& 42 & 0.97  & 42 & 46& 0.97& 47& 51&0.96  \\
0.50 &  52 &  55& 0.97  & 57& 63 & 0.97    &  62& 69& 0.98 &68 &76 & 0.97 \\
\hline
% etc. ...
\end{tabular}
\end{center}
\caption{Binary model with $\gamma^*=0.05$, $d=3$
$s =\lceil \lambda n \rceil$,
various  $n$ and $\lambda$.}
\label{table1}
\end{table}

\begin{table}[h]

\begin{center}
\begin{tabular}{|c||c|c|c|c|c|c|c|c|c|c|c|c|}
\cline{2-13}
\multicolumn{1}{c|}{} &\multicolumn{3}{c|}{$n=20$}&\multicolumn{3}{c|}{$n=50$}&\multicolumn{3}{c|}{$n=100$}&\multicolumn{3}{c|}{$n=200$}  \\
\cline{1-13}
$\lambda$&$N_{\gamma^*} $ & $N_\gamma$ & $1-\gamma $ &$N_{\gamma^*}$ & $N_\gamma$  & $1-\gamma$&$N_{\gamma^*}$ & $N_\gamma$  & $1-\gamma$& $N_{\gamma^*}$ &$N_\gamma$  & $1-\gamma$ \\
\hline
\hline
0.10  &  32 & 34  & 0.92 & 40 &  44&0.93 & 45 & 50 &0.93 & 51& 55&0.94  \\
0.20 & 22 &   24 &  0.92 & 28   &  31&0.93   & 32 &34  &0.92 & 36&38 & 0.93 \\
 0.30& 26 &   28 &  0.93 & 30 & 32 &  0.93 &33  &35  &0.93 &36 &  38& 0.93\\
0.40 &  29&  32 &  0.93 & 35 & 39 & 0.94  & 39 &43 &0.94 &44 & 47& 0.94 \\
0.50 &  46 & 51 &  0.95 &52 &  59&     0.96& 58 &65  & 0.96 &64 & 72& 0.95 \\
\hline
% etc. ...
\end{tabular}
\end{center}
\caption{Binary model with $\gamma^*=0.1$, $d=3$
$s =\lceil \lambda n \rceil$,
various  $n$ and $\lambda$.}
\end{table}

\begin{table}[h]

\begin{center}
\begin{tabular}{|c||c|c|c|c|c|c|c|c|c|c|c|c|}
\cline{2-13}
\multicolumn{1}{c|}{} &\multicolumn{3}{c|}{$n=20$}&\multicolumn{3}{c|}{$n=50$}&\multicolumn{3}{c|}{$n=100$}&\multicolumn{3}{c|}{$n=200$}  \\
\cline{1-13}
$\lambda$&$N_{\gamma^*} $ & $N_\gamma$ & $1-\gamma $ &$N_{\gamma^*}$ & $N_\gamma$  & $1-\gamma$&$N_{\gamma^*}$ & $N_\gamma$  & $1-\gamma$& $N_{\gamma^*}$ &$N_\gamma$  & $1-\gamma$ \\
\hline
\hline
0.10  &  24 &  29 & 0.84 & 34 &39  &0.87 & 38 & 44 &0.87 & 43&49& 0.88 \\
0.20 & 18 &   21 & 0.84  &  24  & 27 &0.85   & 26 &31  &0.86 &29 &34 &0.86  \\
 0.30& 21 &   24 &  0.85 & 25 & 28 & 0.85  & 27 &  31&0.85 & 30& 35 &0.87 \\
0.40 & 24 &  28 &  0.86 & 30 & 35 &  0.87 &  34& 39&0.88 & 37& 43&  0.88\\
0.50 &  37 & 45 &0.90   & 44 &  53&   0.89  & 50 &  60& 0.92 &53 &67 & 0.95 \\
\hline
% etc. ...
\end{tabular}
\end{center}
\caption{Binary model with $\gamma^*=0.25$, $d=3$
$s =\lceil \lambda n \rceil$,
various  $n$ and $\lambda$.}
\label{table3}
\end{table}

From Tables~\ref{table1}--\ref{table3} and Figures~\ref{Bonf_n_50}--\ref{Bonf_n_200} we can draw the following conclusions. For small values of $\gamma$, the value of $\gamma^*$ is very close to $\gamma$ (equivalently $N_\gamma$ is very close to $N_\gamma$).  For larger values of $\gamma$, we see that $\gamma^*$ is often very conservative with the true $\gamma$ being significantly smaller.

\AZ {We  use the following decoding technique for random designs and improved random designs of Section~\ref{sec:improve}.
We start with the COMP procedure described in the beginning of Section~\ref{sec:disj} to eliminate uniquely defined non-defective items. 
Then, in the case where the defective factors are unknown,  we perform several additional individual tests to exactly locate the defective items (such tests are very easy to design). 
In  simulation studies we do not need this as the group $T=T_i$ consisting of  defective items is known and we only need to establish whether there is another group $T^\prime=T_j$ giving exactly the same test results. In one random test, the probability that the results coincide is  $p_{ij}$ defined in \eqref{eq:pij}.
As follows from formula \eqref{eq:k_ijb0},  this probability is high only if $|T_i\setminus T_j|=1$; this is used explicitly in the proof of Theorem~\ref{th:5.1} and noticed in the beginning of Section~\ref{sec:as_b_weak}.  In $N$ tests, such probability becomes $p_{ij}^N$ and if $N$ is not very small, $p_{ij}^N$ becomes negligible when $|T_i\setminus T_j|>1$.
The probability $\tilde{p}_{ij}$ that both results are 1 are also small when $|T_i\setminus T_j|>1$. 
Therefore, for checking whether $T$ is not the unique group of items consistent with 
all the test results, 
it is enough to only check item groups $T^\prime$ with  $|T\setminus T^\prime|=1$.
The same considerations can be used for the additive and other group testing models.
}

\subsection{{Improving on random designs in group testing problems}}
\label{sec:improve}

{

Any $N$-point design $\D_N= \{X_1,\ldots, X_N\}$ has an equivalent  matrix representation as an $N\times n$-matrix $\M(\D_N)$ where columns relate to items and rows to test groups.
Let $a_{i,j}=1$ if item $a_j$  $(j=1,\ldots,n)$ is included into the test group $X_i$ $(i=1,\ldots,N)$; otherwise $a_{i,j}=0$. Then the test matrix corresponding to design $\D_N$ is
$ \M(\D_N):=(a_{i,j})_{i,j=1}^{N,n}\, $. We shall denote the rows of $\M(\D_N)$ by ${\cal X}_i := (a_{i,1},\ldots,a_{i,n})$ for $i=1,\ldots,N$.
A design is called {\it constant-column-weight design} if  all columns of $ \M(\D_N)$ have the same number of ones whereas for a {\it constant-row-weight design}   all rows of $\M(\D_N)$ have the same number of ones. The designs which are both constant-row-weight and constant-column-weight designs are referred to as doubly regular designs, see Section~1.3 in~\cite{Johnson}.
If, for a given design, one of the constancy assumptions is approximately true, we shall use the prefix `near-constant'.

In the most important  case $\X=\PP_n^s$, all designs (including random designs and the designs constructed in this section) are automatically constant-row-weight designs. To improve on the separability properties of random designs, we will construct  near-constant-column weight designs and hence our designs will be nearly doubly regular designs.
Moreover,
we will impose  restrictions on the Hamming  distance between the tests (equivalently the rows of $ \M(\D_N)$). Summarizing, the designs  of this section will have near-constant-column weights, constant-row-weights  and have an additional restriction on the Hamming  distance between the rows of  $\M(\D_N)$. \AZ{ Notice that the fact that keeping large Hamming distances between columns of the test matrix  $\M(\cdot)$ tend to improve  separability properties of the design has been noted in group testing literature, see e.g. \cite{aldridge2016improved}. Moreover, the main idea behind the $d$-disjunct designs of Macula \cite{macula1996simple} is maximization of the minimal Hamming distance between these columns.
} 

%In group testing problems, the test groups $X_i$ ($i=1, \ldots, N$) can be represented as row-vectors
%$X_i=(a_{i,1}, \ldots , a_{i,n})$ of size $n$, where all $a_{i,j} \in \{0,1\}$.

Here we shall describe the algorithm of  construction of the nested designs  we propose; a formal description as a  pseudo-code for the algorithm can be found in Appendix B. We start with a one-element design $\D_1= \{X_1\}$, where $X_1$ is  a random group. At $k$-th step we have a design $\D_{k-1}= \{X_1, \ldots, X_{k-1}\}$ and we are looking for a new test group $X_{k}$ to be added to the design $\D_{k-1}$. To do this, we generate 100 candidate test groups $U_k = \{ X_{k,1},\ldots, X_{k,100} \}$ with  $X_{k,i}\in \PP_n^s$ according to the following procedure. For 75 of the candidate tests, repeat the following. Check the frequency of occurrence of each item and locate the items with the smallest number of occurrences. If there are greater than $s$ of these items, return a random sample of size~$s$. If there are fewer than $s$, say $s'$, such lowest-frequency items, return all $s'$ items and supplement the remaining $s-s'$ items with a random sample from the group containing items that have not appeared the fewest. This describes  Algorithm 1 in the Appendix B. To form the remaining 25 candidate tests, we simply sample them randomly from $\X=\PP_n^s$.
The 100 candidate tests chosen in this manner  encourage nearly equal column weights of the constructed designs $\D_k$ for all $k$.
Of the 100 candidates of the set $U_k$, we select a single test group as $X_{k}$  by maximizing the smallest Hamming distance to all previous points in the design  $\D_{k-1}$.
Specifically, we locate any test group (or groups)  $X'\in U_k$ such that $\min_{1\leq j \leq k-1}d_H(X,X_j) \to \max_{X \in U_k}$. This may result in more than one such  $X'$. If this occurs, we select the group $X'\in U_k$ such that $\sum_{i=1}^{N}d_H(X',X_i)$ is largest.
 This whole process is described as Algorithm 2 in Appendix B.
}

For the random design $\D_N=\{X_1, \ldots, X_{N}\}$ with each $X_i \in \D_N$ chosen independently and uniformly in $\PP_n^s$, the distribution of the Hamming distance between any two rows of $\M(\D_N)$ can be computed. Without loss of generality, we only need to consider the first and second rows of $\M(\D_N)$, that is ${\cal X}_1$ and ${\cal X}_2$. The random variable of interest is $d_H({\cal X}_1,{\cal X}_2)$.
Assume $s\leq n/2$. Then for $x=0,1,\ldots s$ we clearly have
\bea
{\rm Pr}\{ d_H({\cal X}_1,{\cal X}_2) = 2x \} = {{s\choose s-x}{n-s \choose x}}/{{n \choose s}} \,.
\eea

\begin{figure}[h]
\centering
\begin{minipage}{.49\textwidth}
  \centering
  \includegraphics[width=1\textwidth]{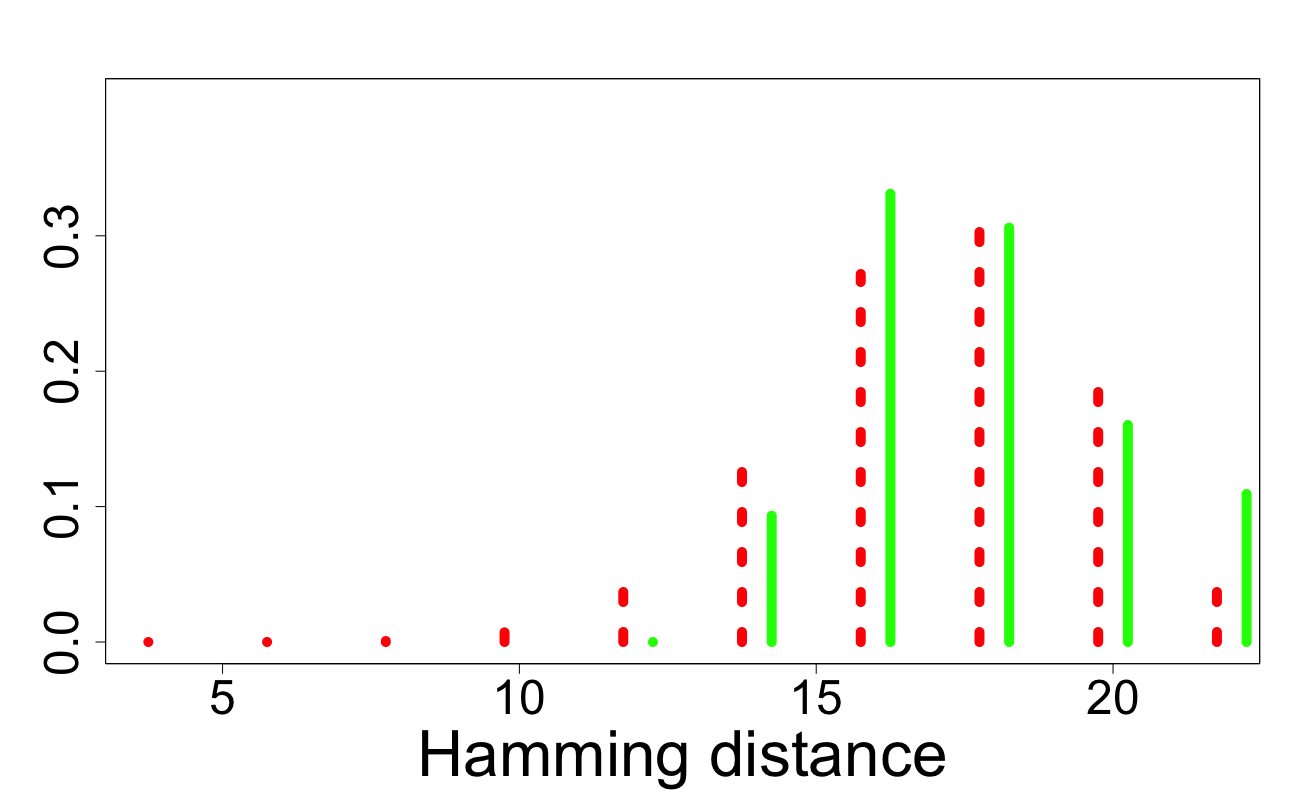}
   \captionsetup{width=.9\linewidth}
\caption{Distribution of inter-point Hamming  distances for random (red) and after the application of Alg. 1 (green); $n=50$ and $s=11$. }
\label{interpoint1}
\end{minipage}%
\begin{minipage}{.49\textwidth}
  \centering
 \includegraphics[width=1\textwidth]{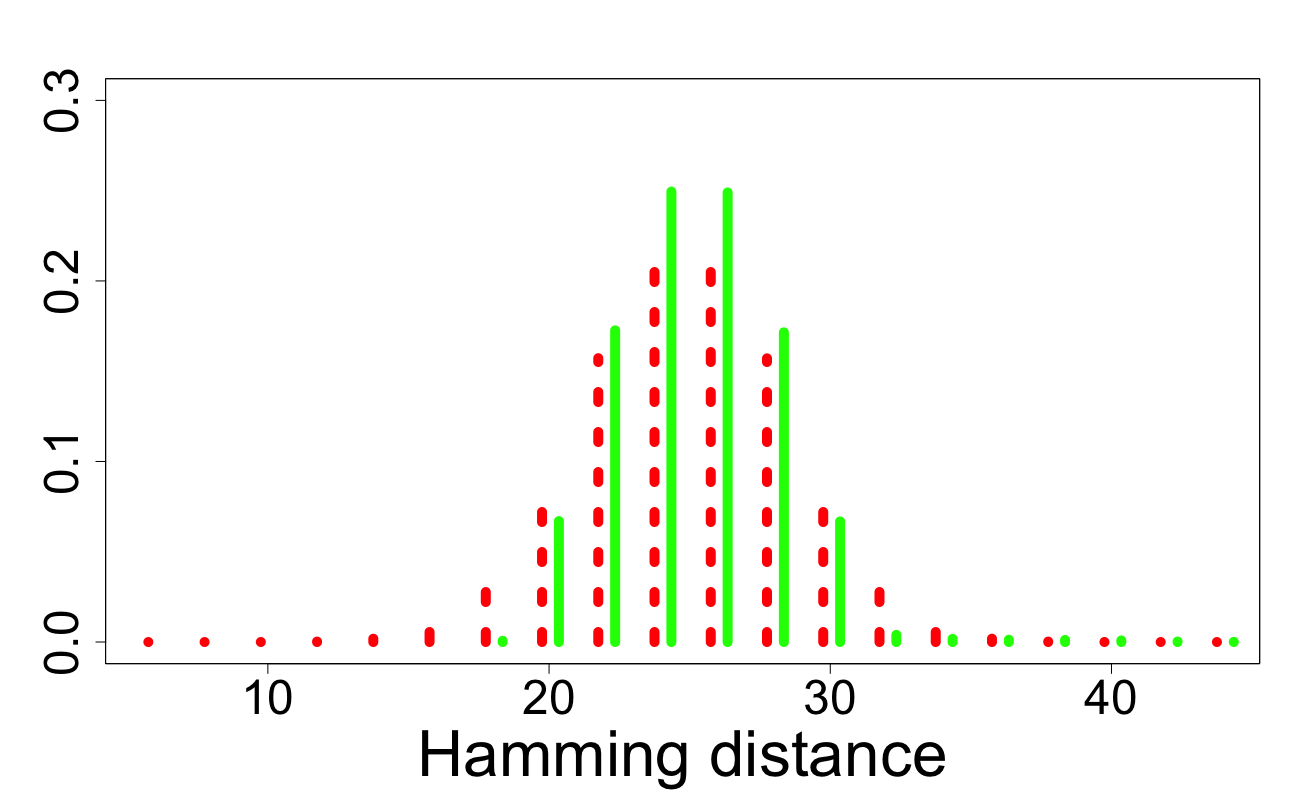}
  \captionsetup{width=.9\linewidth}
\caption{Distribution of inter-point Hamming  distances for random (red) and after the application of Alg. 1 (green); $n=50$ and $s=25$. }
\end{minipage}
\end{figure}

\begin{figure}[h]
\centering
\begin{minipage}{.49\textwidth}
  \centering
  \includegraphics[width=1\textwidth]{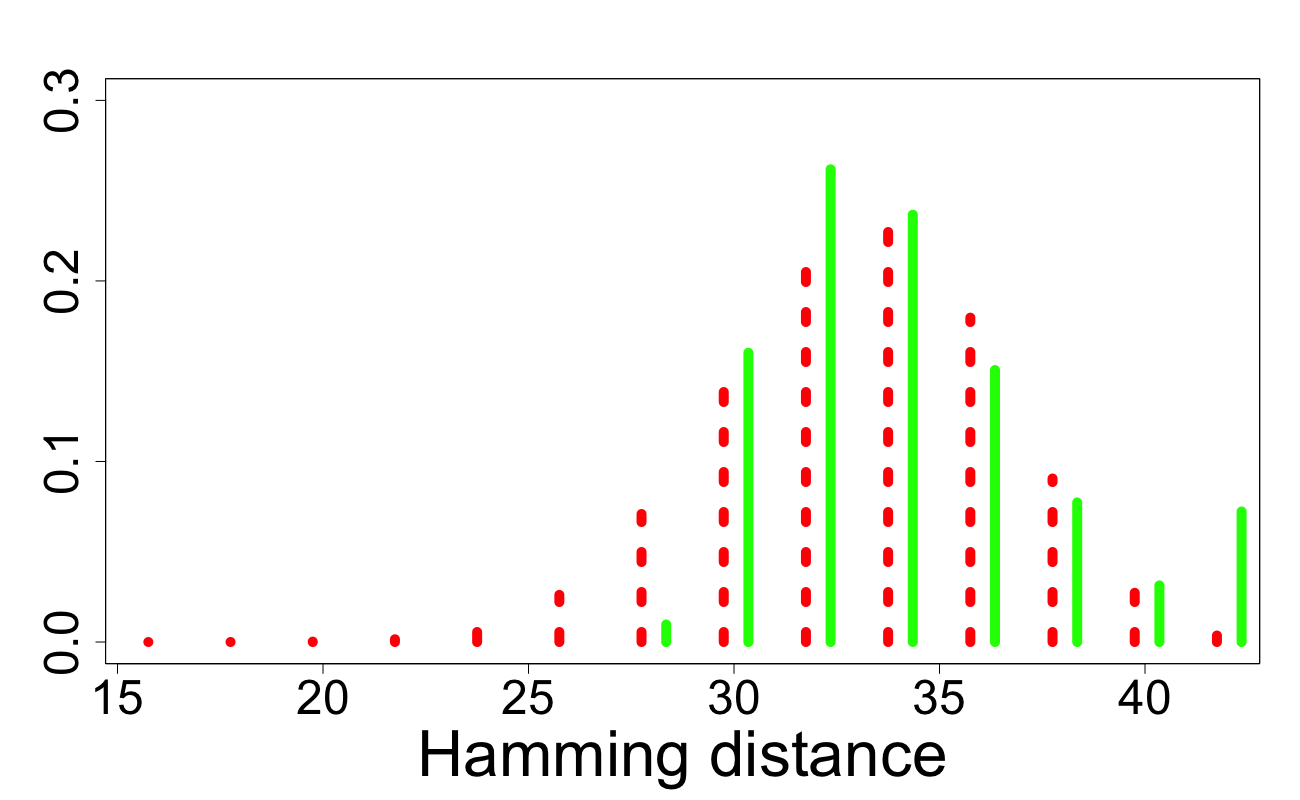}
   \captionsetup{width=.9\linewidth}
\caption{Distribution of inter-point Hamming  distances for random (red) and after the application of Alg. 1 (green); $n=100$ and $s=21$. }
\end{minipage} %
\begin{minipage}{.49\textwidth}
  \centering
 \includegraphics[width=1\textwidth]{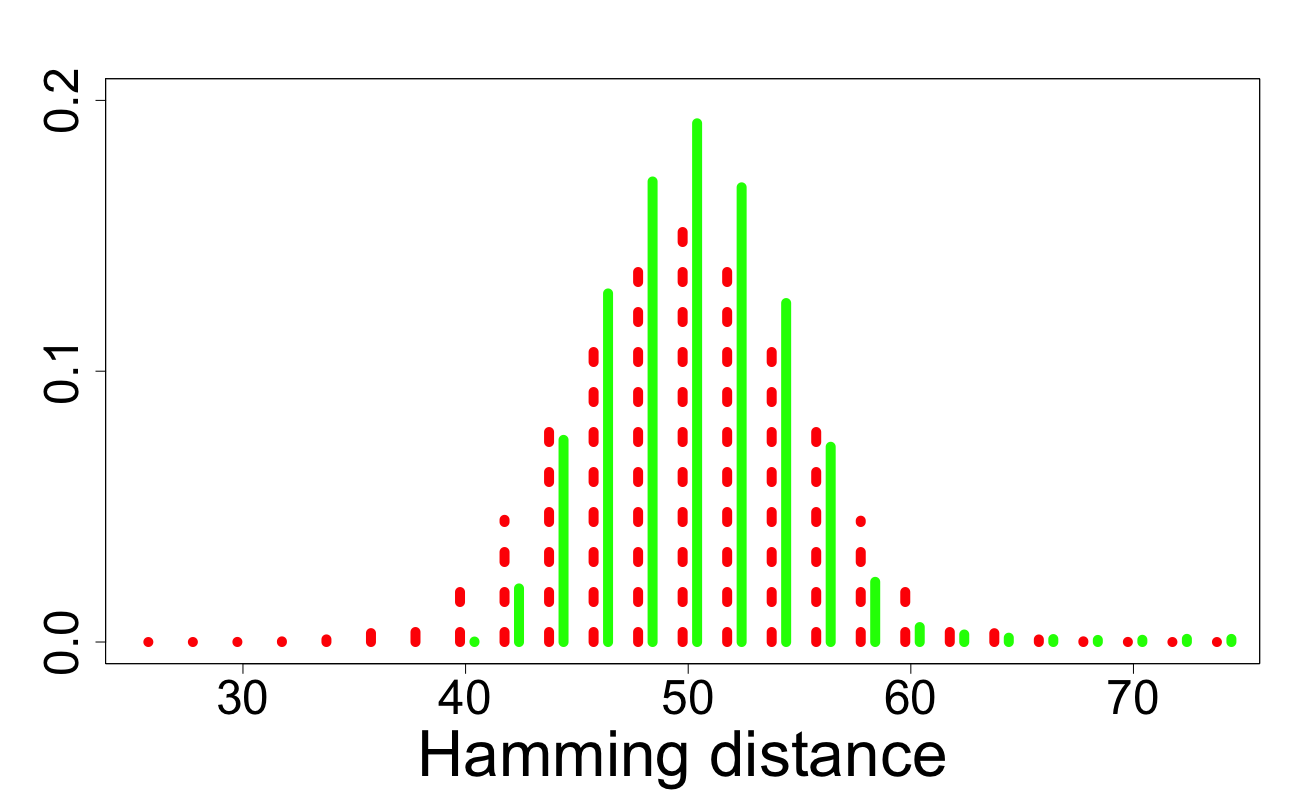}
  \captionsetup{width=.9\linewidth}
\caption{Distribution of inter-point Hamming  distances for random (red) and after the application of Alg. 1 (green); $n=100$ and $s=50$. }
\label{interpoint2}
\end{minipage}
\end{figure}

In Figures~\ref{interpoint1}--\ref{interpoint2}, we plot the distribution of inter-row distances of $\M(\D_N)$ in dotted red and $\M(\D_N')$ in solid green, where $\D_N'$ is a  design obtained by Algorithm 2. The truncation of the lower tail of the distribution in red demonstrates that Algorithm 2 performs very well at preventing small Hamming distances and encouraging large ones.

%
%\AZ{To do entropy pics}\\
%Assume $|X_1|=|X_2|=s$, $|X_1 \cap X_2|=s-v$, $|T|=d$, $T$ is uniformly distributed. Consider the binary model. Consider random vector $F=(f(X_1,T), f(X_2,T)$.
%Then
%$$
%{\rm Pr} \{F=(0,0)\}= \a{n-s-v}{d}/\a{n}{d}\, ,
%$$
%$$
%{\rm Pr} \{F=(0,1)\}= {\rm Pr} \{F=(1,0)\}=\left[\a{n-s}{d}- \a{n-s-v}{d}\right]/\a{n}{d}
%$$
% $$
%{\rm Pr} \{F=(1,1)\}=1- {\rm Pr} \{F=(0,0)\}-2{\rm Pr} \{F=(1,0)\}
%$$

\subsection{Simulation study for quasi-random designs }\label{quasi_sim}

\begin{figure}[h]
\centering
\begin{minipage}{.5\textwidth}
  \centering
  \includegraphics[width=1\textwidth]{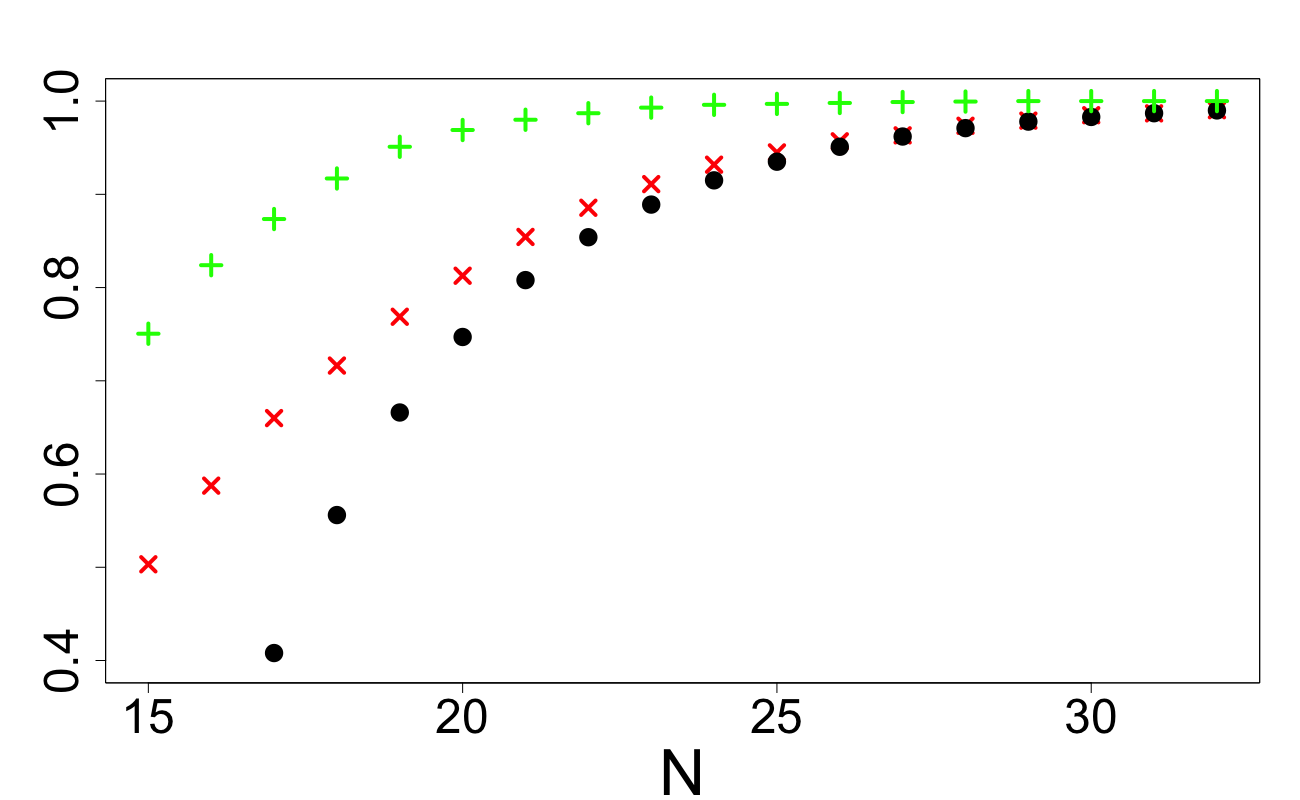}
    \captionsetup{width=.9\linewidth}
\caption{Binary model with $n=20, s=5$;\\ random (red) vs improved random (green). }
\label{alg_pic1}
\end{minipage}%
\begin{minipage}{.5\textwidth}
  \centering
 \includegraphics[width=1\textwidth]{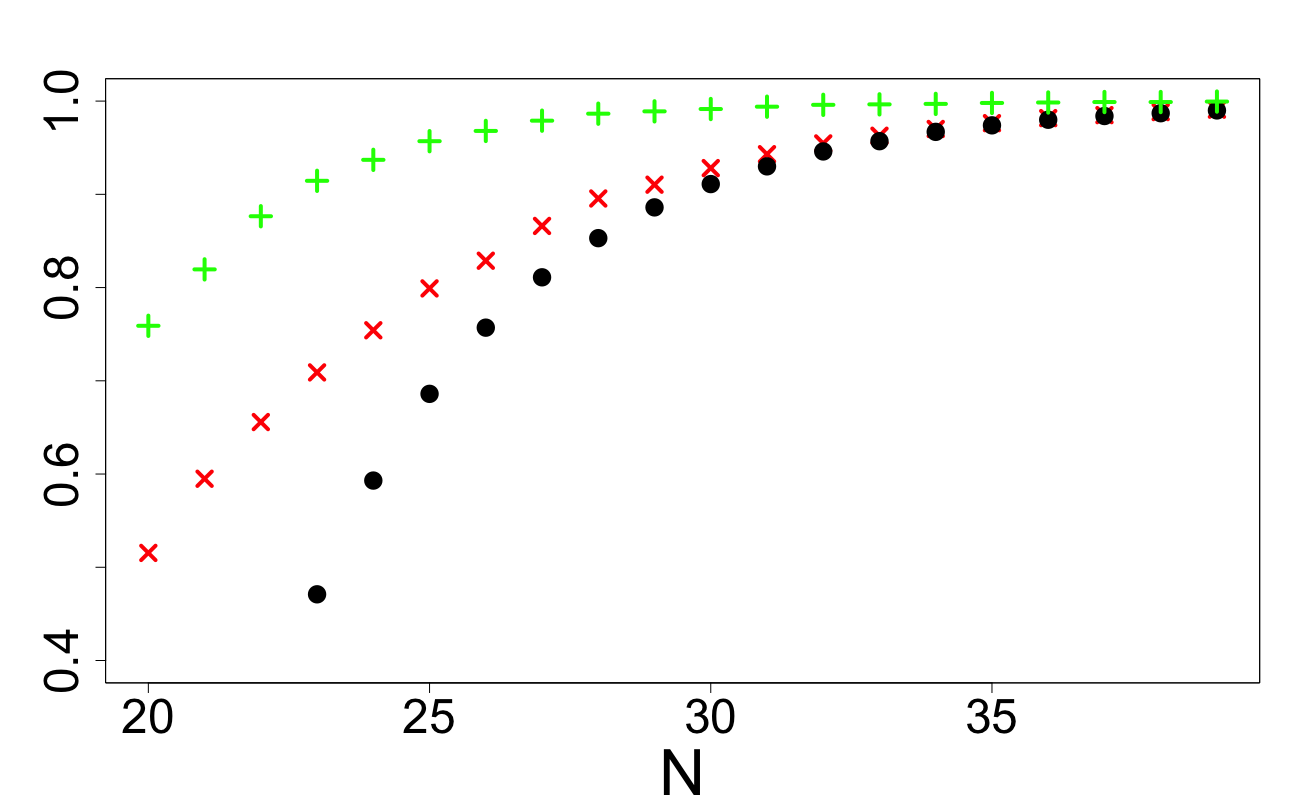}
   \captionsetup{width=.9\linewidth}
\caption{Binary model with $n=50, s=11$;\\ random (red) vs improved random (green). }
\end{minipage}
\end{figure}

\begin{figure}[h]
\centering
\begin{minipage}{.5\textwidth}
  \centering
  \includegraphics[width=1\textwidth]{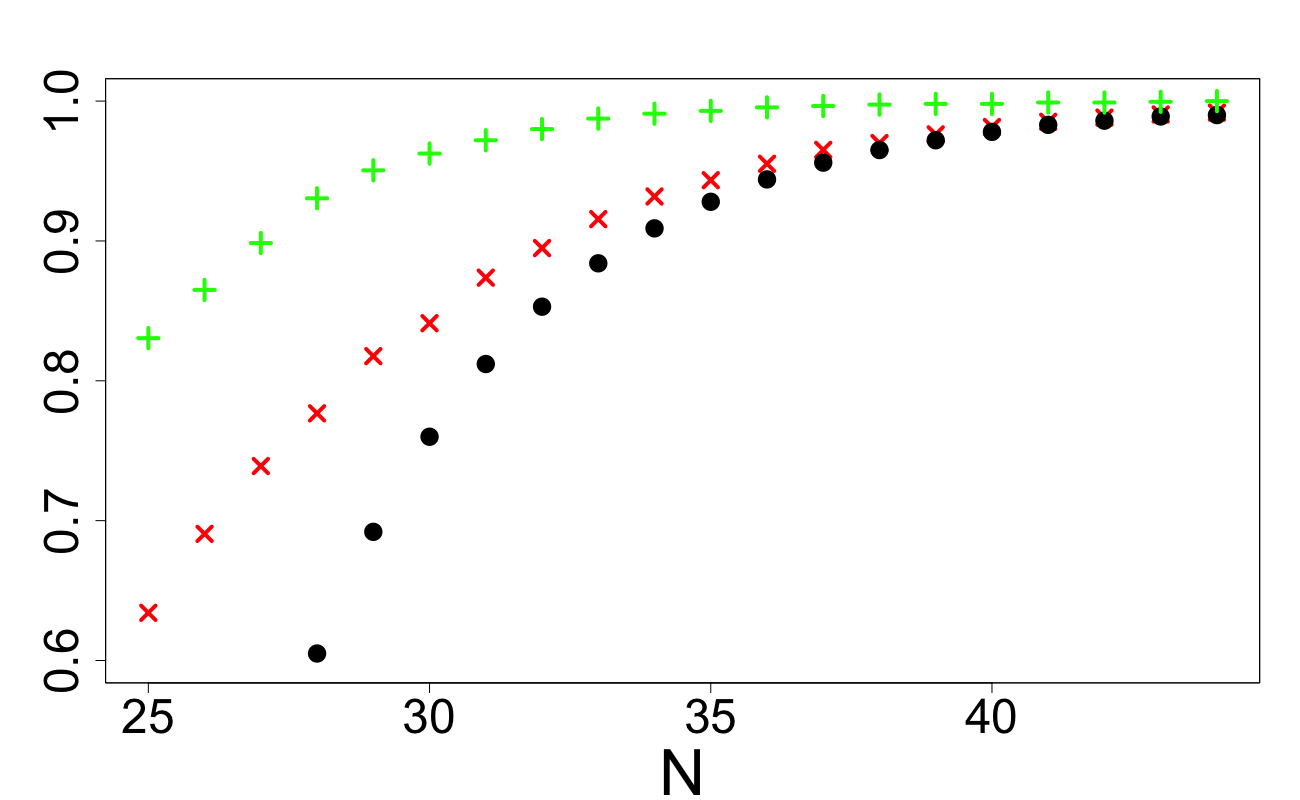}
    \captionsetup{width=.9\linewidth}
\caption{Binary model with $n=100, s=21$;\\ random (red) vs improved random (green). }
\end{minipage}%
\begin{minipage}{.5\textwidth}
  \centering
 \includegraphics[width=1\textwidth]{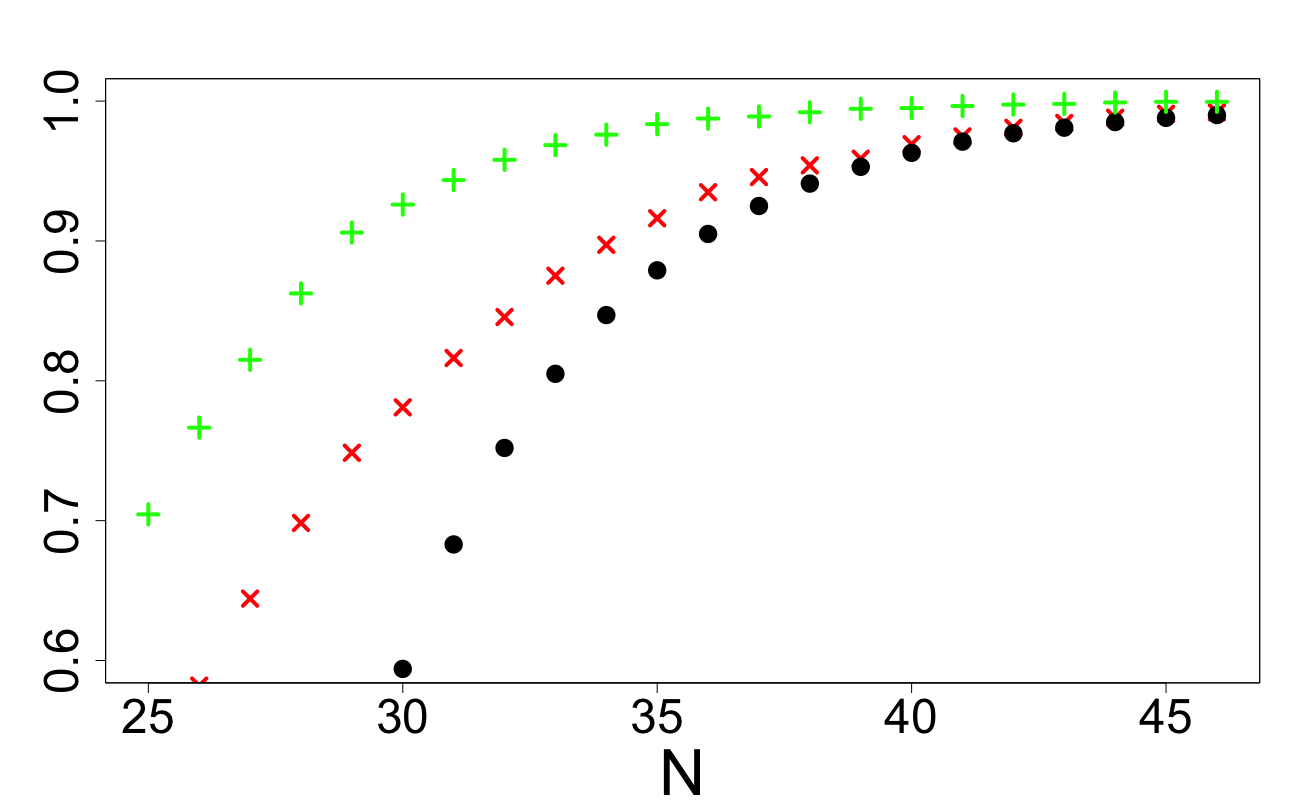}
   \captionsetup{width=.9\linewidth}
\caption{Binary model with $n=150, s=31$;\\ random (red) vs improved random (green). }
\label{alg_pic222}
\end{minipage}
\end{figure}

In Figures~\ref{alg_pic1}--\ref{alg_pic222}, we demonstrate the effect Algorithm 2 has on the probability of separation for the binary group testing problem. Using the red crosses we depict the probability ${\rm Pr}_{\QQ,\RR}\{ T \,$ is separated by $ \D_{N}  \} $ as a function of $N$. With the black dots we plot the value of $1-\gamma^*$ as a function of $N_\gamma$. With green plusses we depict the probability of separation when the design  $\D_{N}'$ is obtained by Algorithm 2. For these figures we have set $d=3$ and $s=s(n) =\lambda_d n$ with $\lambda_d$ chosen asymptotically optimally as $\lambda_d=1- 2^{-1/d}$ (see Section~\ref{sec:as_b_weak}). From these figures we can see Algorithm 2 significantly increases the probability of separation for the binary testing problem. This is particularly evident for smaller values of $N$.

%
%\begin{figure}[h]
%\centering
%\begin{minipage}{.5\textwidth}
%  \centering
%  \includegraphics[width=1\textwidth]{n_20_comparison_additive.png}
%    \captionsetup{width=.9\linewidth}
%\caption{Additive testing problem; $n=20, s=10$.}
%\end{minipage}%
%\begin{minipage}{.5\textwidth}
%  \centering
% \includegraphics[width=1\textwidth]{n_50_comparison.png}
%   \captionsetup{width=.9\linewidth}
%\caption{Additive testing problem; $n=50, s=25$.}
%\end{minipage}
%\end{figure}
%

\subsection{Comparison with designs constructed from the disjunct matrices}
\label{sec:disj}

{
Given a test matrix $ \M(\D_N):=(a_{i,j})_{i,j=1}^{N,n}\, $, let ${\cal S}(a_j):=\{i:a_{i,j}=1   \}$ denote set of tests in which item $a_j$ is included. For a subset ${\cal L}\subseteq {\cal A}$, let ${\cal{S}}({\cal L}) =   \cup_{a_j \in {\cal L}} {\cal S} (a_j)$. Then a test matrix $\M = \M(\D_N)$ is called $d$-disjunct if for any subset ${\cal L}\subseteq {\cal A}$ satisfying $|{\cal L}|=d$ and any $a_j\notin {\cal L}$, we never have ${\cal S}(a_j)\subseteq {\cal{S}}({\cal L})$. A $d$-disjunct matrix can be used to uniquely identify  $d$ or less defective items
 and has the following simple decoding procedure to identify the true defective set: all items in a negative test are identified as non-defective whereas all remaining items are identified as (potentially) defective. This simple procedure is called the combinatorial orthogonal matching pursuit (COMP) algorithm, see \cite[p. 37]{Johnson}.
Consider the following construction of $d$-disjunct matrices $\M$.}

Let $[m] := \{1,2,...,m \}$ be a set of integers.  Then each of the $n$
columns is labeled by a (distinct) $k$ subset of $[m]$. The numbers $m$ and $k$ must satisfy $n \leq {m \choose k}$. Set $\M$ to have ${m \choose d}$ rows with each row labeled by a (distinct) $d$-subset of $[m]$, where $d < k < m$; $a_{i,j}$ = 1 if and only if the label of row $i$ is contained in the label of column $j$.
It was proved in \cite{macula1996simple}, that this procedure makes $\M$  $d$-disjunct. {
The number of rows in $\M$, and hence the number of tests performed, is $N={m \choose d}$ which can be very large and can make identification of the defective set expensive. To avoid a large number of tests, it was recommended in \cite{macula1998probabilistic}  to set $d=2$ regardless of the true $d$; we will call such a matrix 2-disjunct. Whilst the 2-disjunct matrix will no longer guarantee the identification of the defective set if the true $d>2$, it was claimed in \cite{macula1998probabilistic}, see also \cite{d2005construction}, that with high probability the defective set will be identified.}

In Tables~\ref{disjunct_table} and~\ref{disjunct_table2}, we investigate the probability the defective set $T$ is identified when $\T=\PP_n^{3}$ and $\T=\PP_n^{4}$ for designs constructed by the following three procedures: (a) the design corresponding to the 2-disjunct matrix $\M$ with the full decoding; (b) the design corresponding to the 2-disjunct matrix $\M$ with only the COMP procedure used for decoding; (c) $\D_N= \{X_1,\ldots, X_N\} $ with each $X_i \in \D_N $ chosen independently and $\RR$-distributed on $\X=\PP_n^{s}$ where $s$ is chosen according to its asymptotically optimal value (see Section~\ref{sec:as_b_weak}); (d) the design is an improved random design constructed from Algorithm~1.  For different values of $n$, when constructing the  2-disjunct matrix $\M$ we have chosen $m$ and $k$ such that $n \leq {m \choose k}$, $2 < k < m$ and  $N={m \choose 2}$ is as small as possible. For $n=50,100,200$ and $300$, this results in choosing $m=8$ and $k=3$, $m=9$ and $k=4$, $m=10$ and $k=4$ and $m=11$ and $k=4$ respectively. We have then set the random and improved random  designs (constructed from Algorithm 2) (c) and (d) to have the same value of $N$. In these tables, the letter next to $1-\gamma$ corresponds to the procedure used. Within Tables~\ref{disjunct_table} and~\ref{disjunct_table2}, results have been obtained from Monte Carlo simulations with $100,000$ repetitions.

We can make the following conclusions  from the results presented in Tables~\ref{disjunct_table} and~\ref{disjunct_table2}: (i) random designs  are slightly inferior to the designs obtained from 2-disjunct matrices (note, however, that random designs are nested and can be constructed for any $N$), (ii) the COMP decoding procedure alone is insufficient and makes the pair [design, decoding procedure] poor,  and (iii) improved random designs constructed by applying Algorithm~1 have much better separability than both random designs and
the designs obtained from 2-disjunct matrices.

\begin{table}[h]
{
\begin{center}
\begin{tabular}{|c|c|c|c|c|c|c|c|c|c|c|c|}
\cline{1-6}
$n$ & $N$ &$1-\gamma$ (a) &$1-\gamma$ (b)& $1-\gamma$ (c) & $1-\gamma$ (d) \\
\hline
\hline
50 & 28 & 0.99 & 0.82  &0.89 &  0.96 \\
100  & 36 & 0.95 &  0.67 & 0.95&0.97  \\
200 &  45 & 0.98 & 0.70&0.98  & 0.98 \\
300 &  55 & 0.98 &0.77 &0.98  & 0.99 \\
\hline
% etc. ...
\end{tabular}
\end{center}
\caption{Separability comparison for 2-disjunct, random and improved random designs: $\T=\PP_n^{3}$.}
\label{disjunct_table}
}
\end{table}

\begin{table}[h]
{
\begin{center}
\begin{tabular}{|c|c|c|c|c|c|c|c|c|c|c|c|}
\cline{1-6}
$n$ & $N$ &$1-\gamma$ (a) &$1-\gamma$ (b)& $1-\gamma$ (c) & $1-\gamma$ (d) \\
\hline
\hline
50 & 28 & 0.90 &0.51  & 0.53 & 0.86 \\
100  & 36 & 0.76  & 0.26 & 0.70& 0.92 \\
200 &  45 & 0.86 &  0.29& 0.84  & 0.96 \\
300 &  55 & 0.92 & 0.38 & 0.94  & 0.99 \\
\hline
% etc. ...
\end{tabular}
\end{center}
\caption{Separability comparison for  2-disjunct, random and improved random designs: $\T=\PP_n^{4}$.}
\label{disjunct_table2}
}
\end{table}

\subsection{Efficiency of the COMP decoding procedure for random designs} \label{sec:4.8}

\begin{figure}[!h]
\centering
\begin{minipage}{.5\textwidth}
  \centering
  \includegraphics[width=1\textwidth]{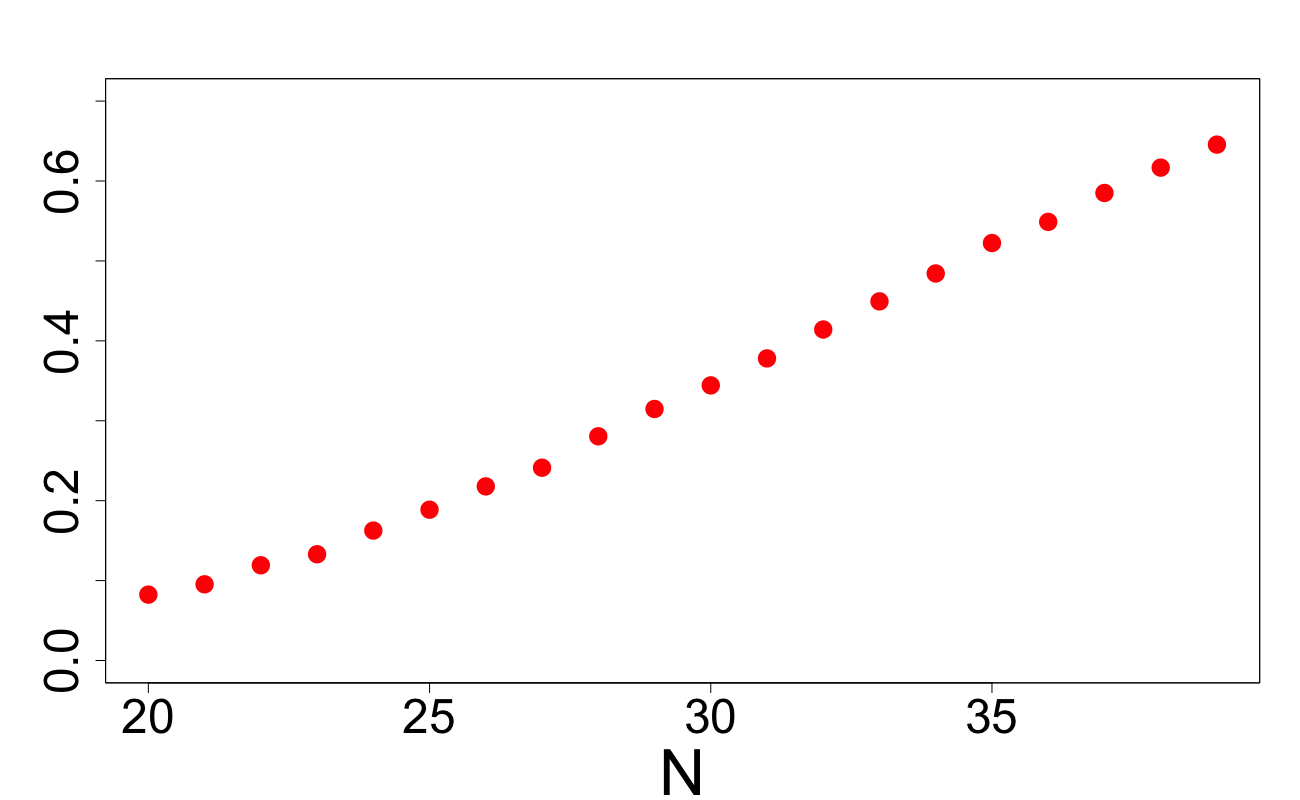}
\caption{Binary model with $n=50, s=11$. }
\label{COMP1}
\end{minipage}%
\begin{minipage}{.5\textwidth}
  \centering
 \includegraphics[width=1\textwidth]{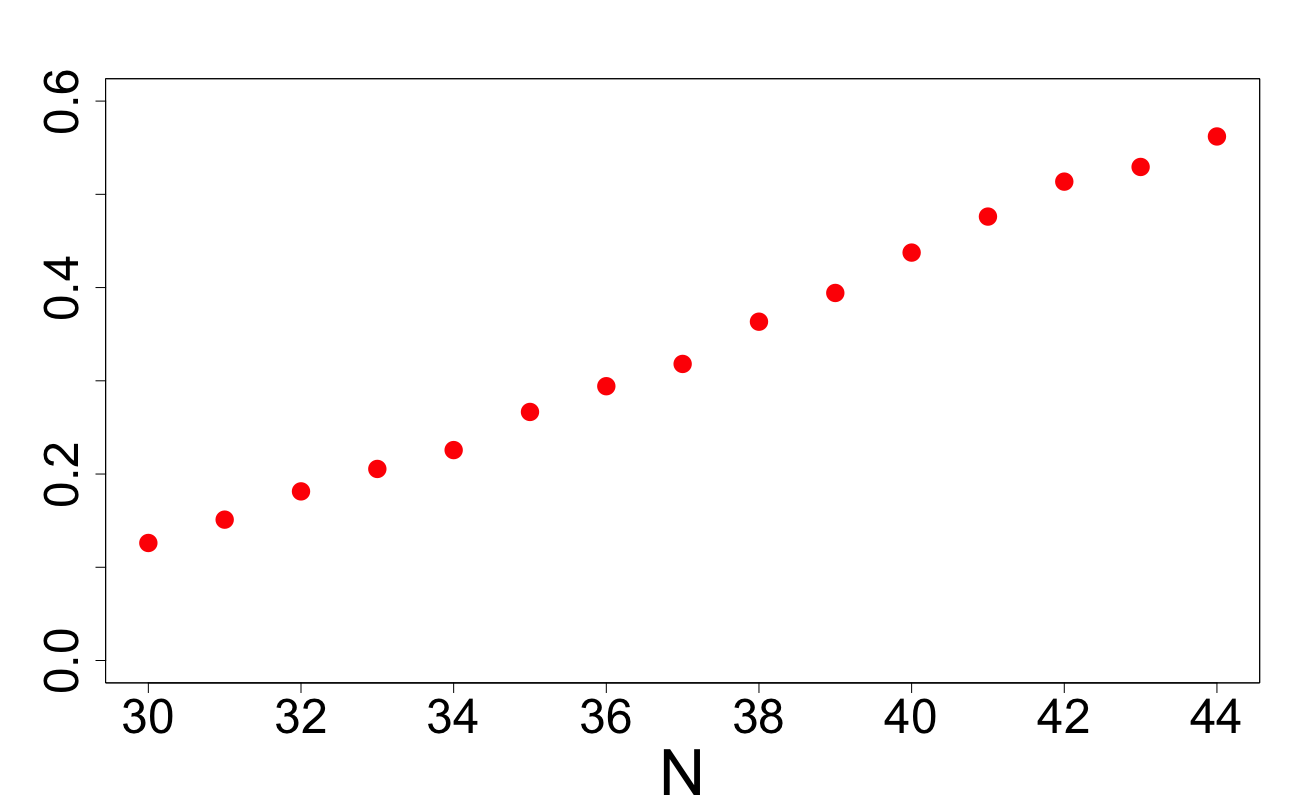}
\caption{Binary model with $n=100, s=21$.  }
\label{COMP2}
\label{alg_pic2}
\end{minipage}
\end{figure}

For a disjunct test matrix $\M$, the COMP decoding procedure described in Section~\ref{sec:disj} is guaranteed to find the defective set and can do so very efficiently (possibly defective items become definitely defective). When the design is not disjunct, say $\D_N$ is constructed randomly, there is no guarantee the COMP procedure will identify the true defective set. Instead, the procedure will provide a set containing the true defective set possibly mixed in with some non-defectives. In \cite[p.37]{Johnson}, the set returned by the COMP algorithm is referred to as the largest satisfying set. For situations when the COMP procedure does not return a uniquely defined $T$, further analysis (based  on the tests with positive results) must be performed to reduce the number of possible target groups of items  $T$ consistent with all available test results.
In Figures~\ref{COMP1}--\ref{COMP2}, we  investigate the  efficiency of COMP expressed as the ratio
\bea
{ \mbox{Pr}_{\QQ,\RR}\{ \text{COMP decoding returns exactly $T$ for  design $\D_N$}  \}  }/{{{\rm Pr}_{\QQ,\RR}\{ T \textrm{ is separated by } \D_{N}  \} }} \,
\eea
for the designs $\D_N$ is constructed randomly. The values in these figures have been obtained from Monte Carlo methods with $50,000$ repetitions.
From these figures we observe that despite for larger  $N$ the COMP procedure has a higher efficiency, this efficiency is still very  low.
{We thus conclude, also taking into account the second conclusion at the end of Section~\ref{sec:disj},  for random designs $\D_N$ the COMP procedure alone will not guarantee identification of the target set frequently enough and must be supplemented by further analysis of positive results.}

\subsection{Binary group testing with lies} \label{sec:4.6}

As discussed in Section~\ref{lies_section}, the results of this paper can be extended to the case where several lies are allowed by introducing the final sum on the right hand side of \eqref{eq:exist_lies0}.
As an example, we shall provide a generalisation of part one of Corollary~\ref{simple_corollary}.

\begin{corollary}\label{simple_corollary_lies}
Let
the test function be defined by
{\rm (\ref{eq:f(X,T)b})}. Let $\T=\PP_n^d$ and $\X=\PP_n^s$, where $n\geq 2$, $1\leq d<n$, $1\leq s<n$ and suppose at most $L$ lies are allowed. Let $\QQ$ and $\RR$ be uniform distributions on $\T$ and $\X$ respectively. For a fixed $N\geq1$, let $\D_N= \{X_1,\ldots, X_N\} $ be a random $N$-point design $\D_N$ with each $X_i \in \D_N $ chosen independently and $\RR$-distributed. Then $\gamma^*(\QQ, \RR,N) $ for the $L$-lie problem can be obtained from \eqref{main_cor1}
by replacing
\bea
\frac{K(\PP_n^s,n,d,d,p)} {{n \choose s}}=
1-2\cdot\frac{{{n-d}\choose {s}}-
{{n-2d+p} \choose {s}}}{{{n}\choose{s}}} \,
\eea
\end{corollary}
with
\bea
\sum_{l=0}^{2L} {{N}\choose {l}}
\left(1-2\cdot\frac{{{n-d}\choose {s}}-
{{n-2d+p} \choose {s}}}{{{n}\choose{s}}} \right)^{N-l}
\left(2\cdot\frac{{{n-d}\choose {s}}-
{{n-2d+p} \choose {s}}}{{{n}\choose{s}}} \right)^l  \, .
\eea

\begin{figure}[!h]
\centering
\begin{minipage}{.5\textwidth}
  \centering
  \includegraphics[width=1\textwidth]{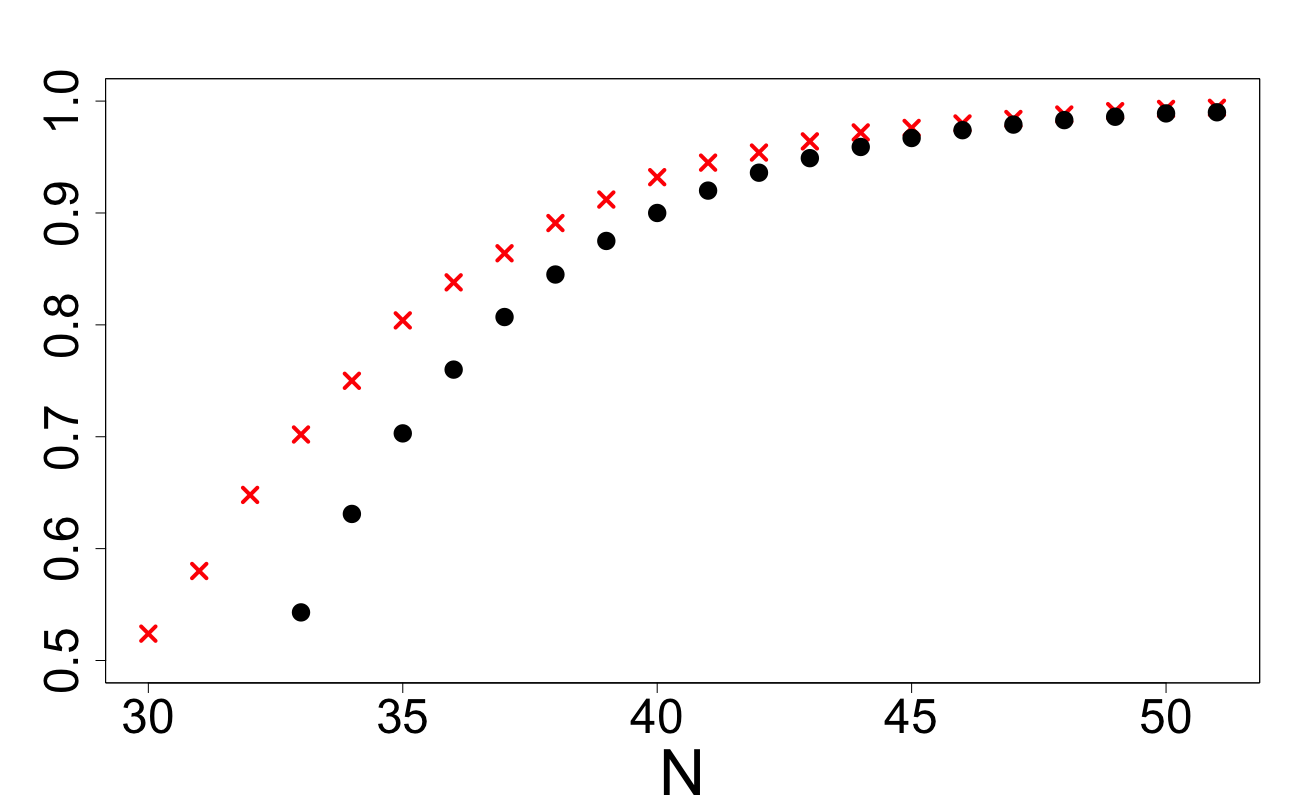}
     \captionsetup{width=.9\linewidth}
\caption{Lies; binary model with $n=20,L=1, s=5$.}
\label{Lies_figure1}
\end{minipage}%
\begin{minipage}{.5\textwidth}
  \centering
 \includegraphics[width=1\textwidth]{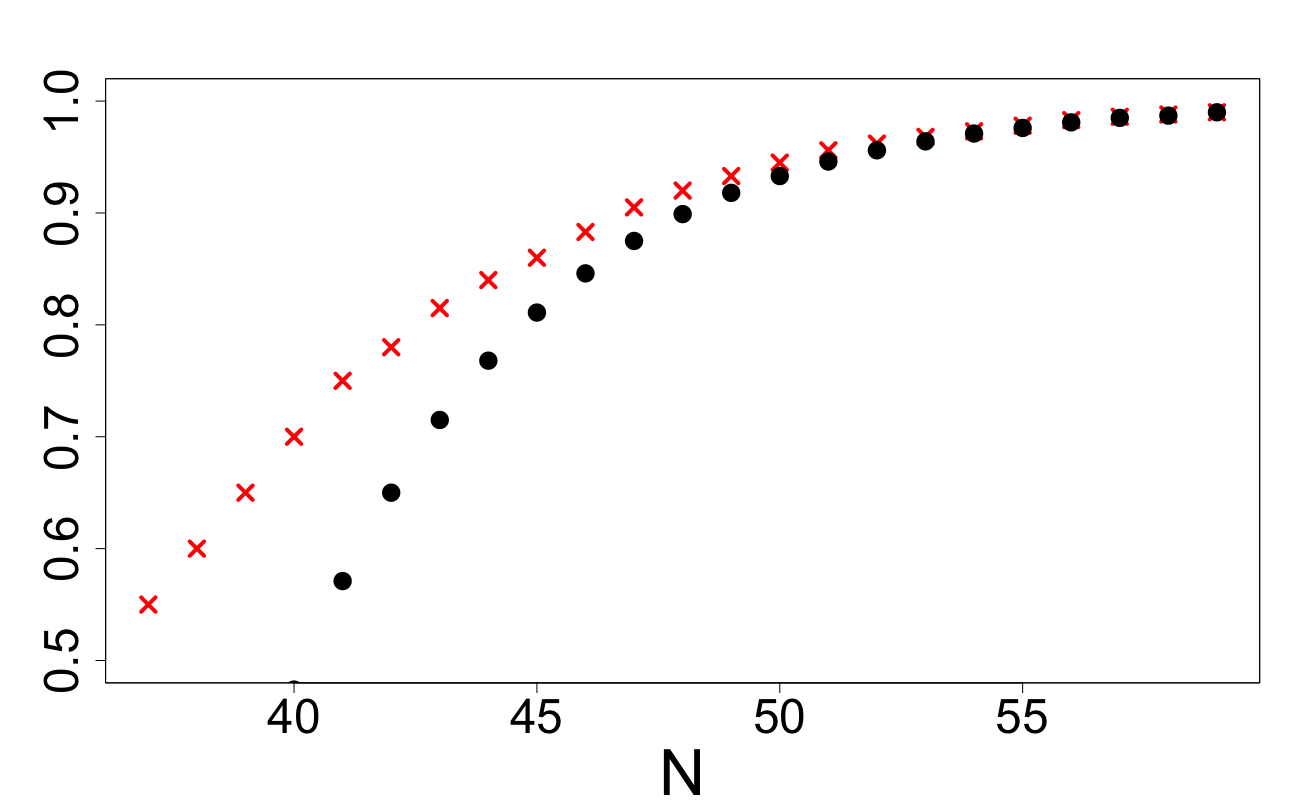}
    \captionsetup{width=.9\linewidth}
\caption{Lies; binary model with $n=50,L=1, s=13$.}
\label{Lies_figure2}
\end{minipage}
\end{figure}

In Table~\ref{Lies_table1} and Table~\ref{Lies_table2}, we document the values of $N^*_\gamma$ obtained from Corollary~\ref{simple_corollary_lies} for $L=1$ and $L=2$ respectively, for several choices of $s$ and $n$. When comparing these tables with Table~\ref{compare_binom1}, we see the significant increase in tests needed when lies are present.
In Figures~\ref{Lies_figure1}--\ref{Lies_figure2}, using red crosses we depict ${\rm Pr}_{\QQ,\RR}\{ T \textrm{ can be uniquely identified by } \D_{N}  \textrm{ with at most $1$ lies} \}$ as a function of $N$. This has been obtain from Monte Carlo methods with $50,000$ repetitions. With the black dots we plot the value of $1-\gamma^*$ as a function of $N_\gamma$ obtained via Corollary~\ref{simple_corollary_lies}. In these figures we have set $s=n/4$ on the basis of Table~\ref{Lies_table1}.
We see once again for small values of $\gamma$, the value of $\gamma^*$ is very close to $\gamma$ (equivalently $N_\gamma$ is very close to $N_\gamma$). For larger values of $\gamma$, we see that $\gamma^*$ is very conservative.

\begin{table}[H]
\begin{center}
\begin{tabular}{|P{7mm}|}
  \multicolumn{1}{c}{  } \\
 \hline
  \multicolumn{1}{|c|}{ $\lambda$ } \\ \hline
  0.05 \\
   0.10 \\
  0.15 \\
  0.20   \\
  0.25 \\
  0.30  \\
  0.35   \\
  0.40  \\
  0.45 \\
        0.50 \\
  \hline
\end{tabular}
\begin{tabular}{ |P{12mm} |P{12mm}|P{12mm}|P{14mm}| }
  \hline
  \multicolumn{4}{|c|}{ $\gamma=0.01$ } \\ \hline
$n=10$ & $n=20$& $n=50$  & $n=100$ \\ \hline
56&126   &130  &166 \\
56&73   & 87 & 95\\
41&  52 & 66 & 72\\
41& 52  & 61 & 66\\
44&  51 &  59& 64\\
59&  53 & 61 & 66\\
59&  67 & 68 & 71\\
59& 67  & 75 & 81\\
98&  81 &  92& 94\\
98&  101 & 109 &115 \\
  \hline
\end{tabular}
\begin{tabular}{ |P{12mm} |P{12mm}|P{12mm}|P{14mm}| }
  \hline
  \multicolumn{4}{|c|}{ $\gamma=0.05$ } \\ \hline
$n=10$ & $n=20$& $n=50$  & $n=100$ \\ \hline
47& 108  &113  &145 \\
47&  63 &  76& 83\\
34&   44&  58& 63\\
34&   44&  53& 58\\
37&   44&  52& 56\\
37&   46&  53& 58\\
50&   58&  59& 63\\
50&   58&  66& 71\\
83&   69&  81& 83\\
83&   87&  96& 102\\
  \hline\end{tabular}
\end{center}
\caption{Values of $N_\gamma$  for binary model with $d=3$, $L=1$, $s =\lceil \lambda n \rceil$,
various  $n$ and $\lambda$.}
\label{Lies_table1}
\end{table}

\begin{table}[!h]
\begin{center}
\begin{tabular}{|P{7mm}|}
  \multicolumn{1}{c}{  } \\
 \hline
  \multicolumn{1}{|c|}{ $\lambda$ } \\ \hline
  0.05 \\
   0.10 \\
  0.15 \\
  0.20   \\
  0.25 \\
  0.30  \\
  0.35   \\
  0.40  \\
  0.45 \\
        0.50 \\
  \hline
\end{tabular}
\begin{tabular}{ |P{12mm} |P{12mm}|P{12mm}|P{14mm}| }
  \hline
  \multicolumn{4}{|c|}{ $\gamma=0.01$ } \\ \hline
$n=10$ & $n=20$& $n=50$  & $n=100$ \\ \hline
 73 &  163 & 166  &210 \\
  73&  94 &  111 & 120\\
  53&  67 &  84 &91 \\
  53&  67 &  78 & 84\\
  57&  66 &  76 & 81\\
  57&   69&  79 & 84\\
  77&   87&   87& 91\\
  77&   87&  96 & 102\\
  127&  104 & 118  &120 \\
  127&   131&  139 & 146\\
  \hline
\end{tabular}
\begin{tabular}{ |P{12mm} |P{12mm}|P{12mm}|P{14mm}| }
  \hline
  \multicolumn{4}{|c|}{ $\gamma=0.05$ } \\ \hline
$n=10$ & $n=20$& $n=50$  & $n=100$ \\ \hline
64 & 143  & 147  &188 \\
 64&  83 &  99 & 108\\
 46&   59&  75 & 82\\
 46&  59 &  69 & 75\\
 50&   58&  68 & 73\\
 50&   61&  70 & 75\\
 67&   77&  78 & 81\\
 67&   77&  86 & 92\\
 111&  92 &  105 & 107\\
 111&   115&   124& 131\\
  \hline\end{tabular}
\end{center}
\caption{Values of $N_\gamma$ for binary model with $d=3$, $L=2$, $s =\lceil \lambda n \rceil$,
various  $n$ and $\lambda$.}
\label{Lies_table2}
\end{table}

{

\section{Asymptotic results}
\label{sec:5}

To start this section, let us make a general comment about the asymptotic expansions in group testing. In most of the known expansions (usually based on the use of Bernoulli designs) the authors are interested in the main asymptotic term only. The authors believe that this is not enough if the asymptotic expansions are intended for the use as (even rough) approximations; see, for example,  a discussion in Section~\ref{sec:as_b_weak} on the asymptotic existence bound in the case of  weak recovery in the binary model. All our expansions in the case of very sparse regime (that is, for fixed $d$) are accurate up to the constant term which we have confirmed by numerous numerical studies. As a result, all our sparse-regime asymptotic expansions can be used as rather accurate approximations already for moderate values of $n$ such as $n=1000$. Typically, this is not so if only the leading term in the expansions is kept. The situation in the sparse regime (when $d \to \infty$ but $d/n \to 0$ as $n \to \infty$) is different and depends on the rate of increase of $d$. If $d$ increases as $\log n$ then once again our expansions are rather accurate up to the constant term. However, if $d=n^\beta +o(1)$ as $n \to \infty$ with some $0<\beta<1$ then we usually can guarantee only the leading term in the expansions and hence the expansions become pretty useless if one wants to use them for deriving approximations. Moreover, our technique completely fails in the case when $d$ grows like ${\rm const}\! \cdot \! n$ as $n \to \infty$.

\subsection{Technical results}
\label{sec:5.1}

The main technical result used for derivation of asymptotic upper bounds in the error-free  environment (no lies) for both  exact and weak recoveries is Theorem 5.1 in \cite{zhigljavsky2003probabilistic}, which we formulate below as Theorem~\ref{th:5.1}. This theorem is especially useful in the case  when ${\X}  = \PP_n^{s}$ with $s=s(n)=\lambda n + o(1)$ (here $0<\lambda<1$ and $n\to \infty$) and ${\T}  $ is either $ \PP_n^{d}$ or $ \PP_n^{\leq d}$ with $d$ fixed (that is, for a very sparse regime). As we show below some results can be extended to a sparse regime when $d \to \infty$ but $d/n \to 0$ as $n \to \infty$. However, unless $d$ tends to infinity very slowly (like $\log n$, for example), we lose the very attractive feature of the expansions, which is the correct constant term.

The authors are not confident that Theorem~\ref{th:5.1} can be applied to the problem of binomial group testing. Also, there are some extra technical difficulties in applying this theorem for  Bernoulli designs. At least, we cannot get the constant term $c$ in \eqref{eq:51d} for   Bernoulli designs (for these designs, the main term $C\log n$ is the same as for our main case ${\X}  = \PP_n^{s}$ with $s=\lambda n + o(1)$
and suitable  $\lambda$).

\begin{theorem}
\label{th:5.1}
Let $I$  be some integer, $c_i,r_i,\alpha_i$
$(i\!=\!1,\ldots,I)$
be some real numbers, $c_i\!>\!0,$ $ 0\!<\!r_i\!<\!1$,
at least one of $\alpha_i$ be positive,
 $\{q_{i,n}\}$, $\{r_{i,n}\}$ be
families of positive numbers $(i\!=\!1,\dots,I)$ such
that
 $0 \!<\! r_{i,n}\! <\!1 $ for all $i$ and
\be
\label{eq:q}
q_{i,n} = c_i n^{\alpha_i}(1+o(1)), \quad r_{i,n}=r_i+o\left(
\frac{1}{\log n}\right)\;\; \;\; {\rm as}\;\; n \ra \infty\,.
\ee
Define $M(n)$ as the solution (with respect to $M$) of the equation
$
%\label{eq:51e}
 \sum^I_{i=1} q_{i,n} r_{i,n}^M \,= 1\,
$
 and set
\be
\label{eq:51a}
N(n) = \min\,\left\{ k=1,2,\ldots\;\mbox{\rm such that } \;\;
 \sum^I_{i=1} q_{i,n} r_{i,n}^k \,< 1 \right\}\,,
\ee
\be
\label{eq:C}
C=\max_{i=1,\ldots,I}\; \frac{\alpha_i}{-\log r_i}\,.
\ee
Finally, let $c$ be the solution of the equation
$ \sum_{j \in
{\cal J}} c_j r_{j}^c \,= 1\, ,
$
where
${\cal J}$
is the subset of the set $\{1,\ldots,I\}$ at which the maximum in
{\rm (\ref{eq:C})} is attained.
Then $ N(n)= \lfloor M(n) \rfloor+1 $
and
\be
\label{eq:51d}
M(n)=C \log n + c +o(1) \quad {\rm as } \quad n\ra\infty\,.
\ee
\end{theorem}

Note that $C$ and $c$ in \eqref{eq:51d} are  constants in the sense that they  do not depend on $n$. Extensive numerical results for   exact and weak recoveries in the binary, additive and multichannel models show that the resulting asymptotic formula \eqref{eq:51d} (in cases
${\X}  = \PP_n^{s}$  and ${\T}  =\PP_n^{d}$ or $ {\T}  =\PP_n^{\leq d}$) is very accurate even for moderate values of $n$. In fact, in all these cases the difference $N(n)-[C \log n + c]$ tends to zero very fast (as $ n \to \infty$)  as long as $d$ is not too large (here $N(n)$ is the upper bound in any of the existence theorems and is defined in \eqref{eq:51a}). In the sparse regime, when $d \to \infty$ (but $d/n \to 0$), the approximation $N(n)\simeq C \log n + c$ is still  accurate but $n$ has to be significantly larger for this approximation to have close to zero accuracy. To distinguish the cases of exact recovery ($\gamma=0$) and weak recovery ($\gamma>0$) we shall write $M_0(n)$ for the upper bounds \eqref{eq:51d} in case of exact recovery  and
$M_\gamma(n)$ in case of weak recovery.

As follows from Theorem~\ref{th:sec6.2} and Corollary~\ref{cor:2} of Section~\ref{sec:3.1} for weak recovery (similar considerations are true for exact recovery), in cases
${\X}  = \PP_n^{s}$  and either ${\T}  =\PP_n^{d}$ or $ {\T}  =\PP_n^{\leq d}$, the existence bounds have the form \eqref{eq:51a}.
Establishment of the asymptotic relations \eqref{eq:q}, from which everything else follows, is usually a straightforward application of the following two simple asymptotic formulas (see Lemmas~5.1 and 5.2 in \cite{zhigljavsky2003probabilistic}).
\begin{itemize}
  \item[(a)] Let $n \ra \infty$,   $u$ and $w$ be positive integers and
$s\!=\!\lambda n\! +\!O(1)$ as $n\ra\infty$
($0\!<\!\lambda\!<\!1$). Then
\bea
{{\a{n-w}{s-u}}}\big/{{\a{n}{s}}} = \lambda^u (1-\lambda)^{w-u}
+ O \left({1}/{n}\right)\;\; {\rm as} \;\; n\ra\infty\,.
\eea
  \item[(b)] Let
  $Q(n,l,m,p)$
be as in \eqref{eq:Q},
$p$, $m$, $l$ be fixed and $n\ra\infty$. Then
  \bea
%\label{eq:Q1}
&&Q(n,l,m,p) = c_{l,m,p}\cdot n^{l+m-p}
\left(1+O\left({1}/{n}\right)\right), \quad n\ra\infty\,,\\
{\rm with}\;\;\;\;&&c_{l,m,p}=\left\{\begin{array}{ll}
                  {1}/\left[{p!(m\!-\!p)!(l\!-\!p)}\right] & {\rm if} \;\; m \neq l\,, \\
                  {1}//\left[{2p!((m\!-\!p)!)^2}\right] & {\rm if} \;\; m=l\,.
                 \end{array}\right.
\eea
\end{itemize}

The set ${\cal J}$ of Theorem~\ref{th:5.1} determines the set (or sets) $\T(n,l,m,p)$ (see \eqref{eq:T}) of pairs of target groups $(T,T')$ which are most difficult to separate by the random design. Theorem~\ref{th:5.1} establishes that by the time the pairs from these set/s $\T(n,l,m,p)$ will be separated (in the case of weak recovery, with probability $1-\gamma$), the pairs $(T,T')$  from all other sets $\T(n,l,m,p)$ will be automatically separated with much higher probability which is infinitely close to 1.
In most cases, the set ${\cal J}$ defined in Theorem~\ref{th:5.1} contains just one number and hence computation of the constant $c$ in \eqref{eq:51d} is immediate. Even if this is not the case, as in \eqref{eq:c_t} below, a very accurate approximation to the exact value of $c$ can be easily found.

\subsection{Additive model}
\label{sec:add_exact}

For the additive model, the case $\T=\PP_n^{\leq d}$  is not very interesting (the same applies to the Binomial testing) as we can make an initial test with all items included into the test group and hence determine the total number of defectives. Therefore, we only consider the case $\X\!=\!\PP_n^s$, $\T\!=\!\PP_n^d$.
Assume
$n\ra\infty$,
$s=s(n) =\lambda n +O(1)$ when
$n\ra\infty$,
 $0 < \gamma < 1.$ The optimal value of $\lambda$ is $1/2$, both for weak and exact recovery. For $\lambda=1/2$, ${\cal J}$ consists of the single index corresponding to $l=m=d$ and $p=0$. This gives for exact  and  weak recovery respectively:
\be
\label{eq:Nasympta}
M_0(n)&=&
  (d+1) \log_2 n\! -\!\log_2 (d-1)!\! -\!1 +o(1) \quad {\rm as } \quad n\ra\infty\,,
\\
\label{eq:51}
M_\gamma(n)&=&
\frac{d\log_2 n \!-\!\log_2 (d!\gamma)}{2d\!-\!\log_2((2d)!)\!+\!2\log_2(d!)} +o(1) \;\;  {\rm as}
\; n \ra \infty\, .
\ee

The asymptotic expressions \eqref{eq:Nasympta} and \eqref{eq:51} have first appeared as
\cite[Corollary 5.1]{ZhigljavskyZ95}.
% but the derivation of the constant terms was only sketched.
Let us make some observations from analyzing  formulas \eqref{eq:Nasympta} and \eqref{eq:51}.

 First, the denominator $F(d)={2d\!-\!\log_2((2d)!)\!+\!2\log_2(d!)}$ in  \eqref{eq:51}  is monotonically increasing with $d$ from $F(2)= 3-\log_2 3 \simeq 1.415$ to $\infty$. This implies that the problem of exact recovery is much more complicated than the problem of  weak recovery and ratio of leading coefficients in \eqref{eq:Nasympta} and \eqref{eq:51} tends to infinity as $d$ increases. This also shows the diminishing role of $\gamma$ in \eqref{eq:51} and the possibility to allow $\gamma$ to slowly decrease as $d$ increases.

 Second, the asymptotic expansion of $F(d) $ at $d=\infty$ is
  $F(d)= \frac12 \log_2 (\pi d )  +O \left( 1/{d} \right)$ with the respective approximation $F(d)\simeq \frac12 \log_2 (\pi d )$ being very accurate for all $d$.
  Stirling formula also gives $\log_2(d!) = d \log_2(d/e) + \frac 12 \log_2(2\pi d)+ O \left( 1/{d} \right)$ as $ d \to \infty$.
This allows us to write the following asymptotic version of \eqref{eq:51} in the sparse regime with $ d= n^\beta +O(1)$ and  $0<\beta<1$ as
\bea
%\label{eq:51as}
\!\!\!M_\gamma(n)=
\frac{ n^\beta (1+2 (1-\beta)\log n)}
{ \log (\pi n^\beta )}  \!
+O(1) \; \; {\rm as}
\; n \!\ra\! \infty\,.
\eea
The sparse-regime version of \eqref{eq:Nasympta} is very clear and need only the expansion
$\log_2((d-1)!) = d \log_2(d/e) + \frac 12 \log_2(2\pi /d)+ O \left( 1/{d} \right)$ as $ d \to \infty$. Thus, for   $d=\lfloor n^\beta \rfloor$ with $0<\beta<1$ we obtain
\bea
%\label{eq:51as1}
M_0(n)=
(\lfloor n^\beta \rfloor+1+\beta/2) \log_2 n\!
 +O((1) \; \;\;
 {\rm as}
\; n \!\ra\! \infty\,  .
\eea
%We believe that even the constant terms in \eqref{eq:51as1} are correct.

\subsection{Binary model, exact recovery}
\label{sec:bin_exact}

Consider first the case of exact recovery in the binary model with $\T=\PP_n^{d}$, $\X=\PP_n^{s}$ and $s=s(n) = \lambda n+O(1)$.
From Corollary 5.2 in \cite{zhigljavsky2003probabilistic} we obtain the following:
the optimal value of $\lambda$ is $\lambda=1/(d+1)$ for which the set ${\cal J}$ of Theorem~\ref{th:5.1} consists of one index corresponding to $l=m=d$ and $p=d-1$; this gives
\be
\label{eq:Nas01}
M_0(n)=
\frac{(d+1)\log_2 n-\log_2(d-1)!-1}
{-\log_2\left(1-{2d^d}/{(d+1)^{d+1}}\right)}\, +o(1) \quad {\rm as } \quad n\ra\infty\,.
\ee
The numerator in \eqref{eq:Nas01} coincides with the rhs in \eqref{eq:Nasympta}. The denominator in the rhs of \eqref{eq:Nas01},  $G(d):=-\log_2 [(1-{2d^d}/(d+1)^{d+1}]
$, provides the coefficient characterizing the complexity of  the binary model with respect to  the additive one.
Function $G(d)$ monotonically decreases from $
G(2)\simeq 0.507$ to 0 with
$G(d)=
2/[ d e \log 2 ]+O \left( {d}^{-2}
 \right)$ for large $d$.
This gives  us  the following  sparse-regime version of \eqref{eq:Nas01} ($  d=\lfloor n^\beta \rfloor,\; 0<\beta<1/2$):
\be
\label{eq:Nas01a}
\!\!\!\!\!\!M_0(n)= \lfloor n^\beta \rfloor e \log \sqrt{2}  \left[
(\lfloor n^\beta \rfloor \!+\!1\!+\!\beta/2) \log_2 n\!
\right]
 \!+\!O(1) \; \;
 {\rm as}
\; n \!\ra\! \infty\,     .\;\;\;\;
\ee

Consider now the case of exact recovery in the binary model with $\T=\PP_n^{\leq d}$, $d>2$, $\X=\PP_n^{s}$ and $s=s(n) = \lambda n+0(1)$.
From Corollary 5.3 in \cite{zhigljavsky2003probabilistic} we obtain the following:
the optimal value of $\lambda$ is $\lambda=1/d$ for which the set ${\cal J}$ of Theorem~\ref{th:5.1} consists of one index corresponding to $l=d$ and $m=p=d-1$; this gives
\be
\label{eq:Nas05}
M_0(n)=
\frac{d\log_2 n-\log_2(d-1)!}
{-\log_2\left(1-{(d-1)^{d-1}}/{d^d}\right)}\, +o(1) \quad {\rm as } \quad n\ra\infty\,.
\ee
The denominator  $H(d):=-\log_2 [\left(1-{(d-1)^{d-1}}/{d^d}\right)]
$  in the rhs of \eqref{eq:Nas05} is noticeably  smaller than the denominator $G(d)$ in the rhs of \eqref{eq:Nas01}.  For large $d$, we have
$H(d)=
1/[ (d-1) e \log 2 ]+O \left( {d}^{-2}
 \right)$.
This gives  us  the following  sparse-regime version of \eqref{eq:Nas05} for $\T=\PP_n^{\leq d}$ and $d=\lfloor n^\beta \rfloor$ with $0<\beta<1/2$:
\be
\label{eq:Nas05a}
\!\!\!\!\!\!M_0(n)= \lfloor n^\beta -1\rfloor e \log {2}  \left[
(\lfloor n^\beta \rfloor +\!\beta/2) \log_2 n
\right]
 \!+\!O(1) \; \;
 {\rm as}
\; n \!\ra\! \infty\,     .\;\;\;\;
\ee
Comparing \eqref{eq:Nas01a}  with \eqref{eq:Nas05a} we can conclude that in the sparse regime with  $d \to \infty$,  the problem of exact recovery in the binary model with $\T=\PP_n^{\leq d}$  is approximately twice harder than in the case of $\T=\PP_n^{ d}$ in the sense that
it requires approximately twice more tests needed to guarantee the exact recovery of all defectives.

\subsection{Binary model, weak recovery} \label{sec:as_b_weak}
Consider now the case of weak recovery; the non-asymptotic version is considered in Corollary~\ref{simple_corollary}.
Assume that
$\T$ is either $\PP_n^{d}$ or
$\T=\PP_n^{\leq d}$,
$d\geq 2$,
$0 < \gamma < 1$,
$\X=\PP_n^s$,
$s=s(n) =\lambda n +O(1)$ when
$n\rightarrow\infty$. Then  the optimal value of $\lambda$ is $\lambda=1- 2^{-1/d}$; for this value of $\lambda$
the set ${\cal J}$ of Theorem~\ref{th:5.1} consists of $d$ indices corresponding to $l=m=d$ and $p=0,1, \ldots,d-1$;
\be
\label{eq:Nasympt2}
  M_{\gamma}(n) =
d\log_2 n +c  +o(1) \quad {\rm as } \quad n\ra\infty\,,
\ee
where $c=c(\gamma,d)$ is the solution of the equation
\be
\label{eq:c_t}
\sum_{p=0}^{d-1}{2^{-c(d-p)/d}
{\displaystyle\frac{d!}{p!(d-p)!^2}}}=\gamma \, .
%\sum_{p=0}^{t-1}{2^{-ct/(t-p)}
%{\displaystyle\frac{t!}{p!(t-p)!^2}}}=1 \, .
\ee
Numerical results show that the asymptotic expansion \eqref{eq:Nasympt2} provides
an approximation $  N_{\gamma}(n) \simeq
d\log_2 n +c $ which is extremely accurate for even moderate values of $n$ such as  $n=10^3$.

By comparing \eqref{eq:c_t} with \eqref{eq:Nas01} and \eqref{eq:Nas05} we conclude that in the case of binary model,  weak recovery (for any $0 < \gamma<1$) is
a much simpler problem than exact recovery.

Since the set ${\cal J}$ of Theorem~\ref{th:5.1} consists of $d$ indices rather than one, the constant $c$ is a solution of the equation containing $d$ summands, see \eqref{eq:c_t}. Despite formally we cannot neglect any of the terms in \eqref{eq:c_t}, keeping just one term, with $p=t-1$, provides an easily computable but rather accurate lower bound for $c$:
$
c\geq c_{\ast}=d  \log_2 (d / \gamma)\, .
$
Table \ref{tableC} shows that the loss of precision in \eqref{eq:Nasympt2} due to  the substitution of $c$ by $c_{\ast}=
d  \log_2 (d / \gamma)$ in \eqref{eq:c_t} is minimal. As a by-product,
Table~\ref{tableC} shows that neglecting the constant term in the asymptotic expressions like \eqref{eq:Nasympt2} would make such asymptotic formulas totally impractical  as in practice $n$ is rarely astronomically large.

\begin{table}[h]

\begin{center}
\begin{tabular}{|c|c|c|c|c|c|c|c|c|c|c|c|}
%\cline{1-9}
\hline
$d$ & 2 & 3 & 5 & 10 & 20 & 30 &40 &50\\
\hline %\hline
$c$ &13.295 & 21.701 &39.858 & 89.722& 199.45& 316.73& 438.91& 564.74\\
\hline
$c_{\ast}$ & 13.288 & 21.686  &39.829 & 89.657& 199.31& 316.53& 438.64 &564.38\\
\hline
\end{tabular}
\end{center}
\caption{Values of $c$ defined as the solution of  \eqref{eq:c_t} and $c_{\ast}=d  \log_2 (d / \gamma)$ for $\gamma=0.02$ and different values of $d$
.}
\label{tableC}
\end{table}

As perhaps the main conclusion of this section, we offer the following approximation for $N_\gamma$
in the case of binary model with $\X=\PP_n^{s}$,
$\T=\PP_n^{d}$ and $\T=\PP_n^{\leq d}$ and $s$ chosen asymptotically optimally by $s=\lfloor n(1- 2^{-1/d})\rfloor$:
\be
\label{eq:b_ex}
N_\gamma(n)\simeq d \log_2 n + d  \log_2 (d / \gamma)\, .
\ee
If we use this formula and  express $\gamma$ through $N_\gamma(n)$, then we get an approximation \be
\label{eq:b_ex1}
\gamma^*(\QQ, \RR,N) \simeq 2^{-N/d} nd\ee
for the value $\gamma^*(\QQ, \RR,N) $ of part one of Corollary~\ref{simple_corollary}. Formulas \eqref{eq:b_ex} and \eqref{eq:b_ex1} connect all major parameters of interest, $n$, $d$, $N$ and $\gamma$, into one simple approximate relation. This relation can clearly show, in particular, allowed rates of increase of $d$ as a function of $n$ guaranteeing the same or even decreasing $\gamma$.

The approximation \eqref{eq:b_ex} is extremely accurate already for very moderate $n$ (say, $n \geq 200$) and not very large $d$. Rather surprisingly, the approximation \eqref{eq:b_ex1} becomes reasonably accurate for moderate $n$ too, as long as the r.h.s. in \eqref{eq:b_ex1} gets small enough. A very simple MAPLE code can provide such a comparison (with almost arbitrary computational precision) for values of $n$ up to $10^6$ and $d$ up to 20 or more. Actually, what is important for formula  \eqref{eq:b_ex1} getting high levels of accuracy is the value of $N$ which has to be large enough; this is consistent with
very high level of accuracy of \eqref{eq:b_ex} for large values of $N_\gamma(n)$.

\subsection{Extensions to noisy testing}
\label{sec:noisy_as}
In \cite{zhigljavsky2010nonadaptive}  a technique is developed of transforming the asymptotic upper bounds \eqref{eq:51d}, obtained from the non-asymptotic expression \eqref{eq:51a}, for an upper bounds for $N$ in the same model  when up to $L$ lies are allowed.
 Theorems 2 and 3 of \cite{zhigljavsky2010nonadaptive}   imply that any asymptotic  bound of the form \eqref{eq:51d} can be rewritten in the form
\be
\label{eq:51dn}
N(n) = C \log n + c_1 \log \log n  +c_0(n) \,,
\ee
where the constant $C$ is exactly the same as in  \eqref{eq:51d} and the constant $c_1$ is computable from the considerations very similar to
indicated in Theorem~\ref{th:5.1}. The main difficulty in using the asymptotic expansion   \eqref{eq:51dn} as an approximation for finite $n$ is related to a rather difficult structure of the function $c_0(n)$, which is bounded (with a computable upper bound) but not monotonic in $n$.
The first term in \eqref{eq:51dn} dominates the asymptotical behaviour of  $N(n)$. However, the constant $c_1$ is always larger than $C$ and,  depending on the allowed number of lies $L$, could be very large. This makes the second term in \eqref{eq:51dn}  significantly more influential than the first term (assuming, for example, $L=5$). Moreover, for small or moderate values of $n$, the values of $c_0(n)$ could also be larger than the main asymptotic term $C \log n$.
}

\section*{Appendix A: Proofs}

\subsubsection*{Proof of Theorem~\ref{th:22}}

 We are interested in computing the value of $\gamma^*$ which satisfies the following.
\bea
&&{\rm Pr}_{\QQ,\RR}\{ T \textrm{ can be uniquely identified by } \D_{N}  \textrm{ with at most $L$ lies} \}= 1-\gamma \\
&=& \sum_{i=1}^{|\T|}{\rm Pr}_{\RR}\{ T_i \textrm{ can be uniquely identified by} \D_{N}  \textrm{ with at most $L$ lies} \}{\rm Pr}_{\QQ}\{ T=T_i \}\\
&=& \sum_{i=1}^{|\T|}{\rm Pr}_{\RR}\{d_H(F_{T_i},F_{T_j})\geq 2L+1  \text{ for all } j\neq i \}{\rm Pr}_{\QQ}\{ T=T_i \}\\
&=& 1- \sum_{i=1}^{|\T|}{\rm Pr}_{\RR}\{d_H(F_{T_i},F_{T_j})\leq 2L  \text{ for at least one } j\neq i \}{\rm Pr}_{\QQ}\{ T=T_i \}\\
&\geq &1- \sum_{i=1}^{|\T|}{\rm Pr}_{\QQ}\{ T=T_i \}\sum_{j\neq i}{\rm Pr}_{\RR}\{d_H(F_{T_i},F_{T_j})\leq 2L  \} = 1-\gamma^* \,.
\eea

For a given design
 $\D_N=\{X_{1},\dots,X_{N}\}$, consider the matrix
$
 \|f(X_{i},T_j)\|_{i,j=1}^{N,|\T|}\,
$
 whose rows correspond to
 the test sets $X_i$ and the columns correspond to the targets $T_j$. Denote the columns of this matrix by $A_j$ ($j=1, \ldots,  |\T|$).

Let $(X_1,X_2,\dots,X_N)$ be a random sample from $\X$.
Then for any fixed pair $(i,j)$ such that $i\neq j$
 $\,(i,j=1,\dots,|\T|)$ and any integer $l$ $\, (0\leq l \leq N)$ we have
\bea \Pr\{d_H(A_i,A_j)=l\} = {{N}\choose {l}} \left(p_{ij}\right)^{N-l}
\left(1-p_{ij}\right)^l
\eea
and therefore
\bea
\Pr\{d_H(A_i,A_j)\leq 2L\} = \sum_{l=0}^{2L} {{N}\choose {l}}\left(p_{ij}\right)^{N-l} \left(1-p_{ij}\right)^l \,. \hspace{50mm} \Box \eea

\subsubsection*{Proof of Theorem~\ref{th:k_ij}}

Let $(T_i, T_j)\in\T(n,l,m,p)$
and $a$ be some integer. Introduce the sets
\bea
\label{eq:X^01}
\nonumber
\X^{a,a}
=\{X\in \X :\,|X\cap T_i|=a, \, |X\cap T_j|=a\}\, , \\
\nonumber
\X^{a,> a}
=\{X \in \X :\,|X\cap T_i|=a, \,|X\cap T_j| > a\}\, , \\
\nonumber \X^{> a,a}= \{X \in \X :\,|X\cap T_i| > a, \,|X\cap T_j|
= a\}\, .
\eea
Remind that $k_{ij}=
\left|\{X\in\X:\,f(X,T_i)=f(X,T_j)\}\right|$ and $f(X,T){=}
\min\{\KK,|X\cap T|\}$.

We have the equality
$f(X,T_i)=f(X,T_j)\,$
if and only if one of the three following cases occurs:
(i)
$X\in \X^{a,a}$ for some $a\geq 0$;
(ii)
$X\in
\X^{a, >a}$ for some $a \geq \KK$;
(iii)
$X\in
\X^{> a, a}$ for some $a \geq \KK$.
Therefore,
\be
\label{eq:k_ij-3}
k_{ij}=
\sum_{a\geq 0} |\X^{a,a}|+
  \sum_{a\geq \KK} |\X^{a,>a }|+
  \sum_{a\geq \KK} |\X^{>a,a}|.
\ee

%By definition,  $X\in\X_{uvr}$
%for some $u,v$ and $r$ such that the sets
%$\X_{uvr}
%=\X_{uvr}(T_i,T_j)$ are defined, see
%(\ref{eq:X-uvr}).
The set of  integers
$n$, $m$, $l$, $p$, $u$, $v$ and $r$ satisfy then the constraints
(\ref{ae:param_g}).
Using these constraints
and the definition of the coefficients
$R(\cdot)$, see
(\ref{def_R}), we can re-express
 the sums in the right-hand side of
(\ref{eq:k_ij-3}) as follows:
$$
\sum_{a\geq 0} |\X^{a,a}|=
 \sum_{r{=}0}^p \sum_{u{=}0}^{m{-}p}
R(n,l,m,p,u,u,r)\, ,
$$
$$
 \sum_{a\geq \KK} |\X^{a,>a }|=
\sum_{r{=}0}^p \sum_
{u{=}w}^{l{-}p}
\sum_{v{=}u{+}1}^{m{-}p}
R(n,l,m,p,u,v,r)\, ,
$$
where
$ w=\max\{0,\KK-r\} $,
and analogously
$$
 \sum_{a\geq \KK} |\X^{> a,a }|=
 \sum_{r{=}0}^p \sum_{v{=}w}^{m{-}p} \sum_{u{=}v{+}1}^{l{-}p}
R(n,l,m,p,u,v,r)\,  .
$$
By substituting this into
(\ref{eq:k_ij-3}) we get
(\ref{eq:k_ij-2}).
To finish the proof we just need  to mention that the above calculation
does not depend on the choice of the pair
$(T_i, T_j)\in\T(n,l,m,p)$ since $\X=\PP_n^{s}$ is balanced.
\hfill $\Box$

\subsubsection*{Proof of Theorem~\ref{th:sec6.2}}

Let $\D_N=\{X_1, \ldots, X_N\}$ be an $\RR$-distributed random design and let $T$ be $\QQ$-distributed. For some $0<\gamma<1$, we have
$
{\rm Pr}_{\QQ,\RR}\{ T \textrm{ is separated by } \D_N \} = 1-\gamma.
$

Let
$
\PP_N=
{\rm Pr} _{\QQ,\RR}\{ T \mbox{ is not separated
by } \D_N \}.$
Then
$
{\rm Pr}_{\QQ,\RR}\{ T \textrm{ is separated by } \D_N \} = 1-\PP_N.
$
By conditioning on $T \in \PP_n^{ \dd }$, for $0\leq \dd  \leq d$, and $\UU$-distributed random variable $\xi$ we have
\bea
\PP_N=
{\rm Pr} _{\QQ,\RR}\{ T \mbox{ is not separated by } \D_N \} =
\suml_{\dd =0}^d
 P_{N,n,\dd }(\X){\rm Pr}_\UU\{\xi =\dd   \} \, ,
\eea
where
$P_{N,n,\dd }(\X)$ is the probability
\bea
P_{N,n,\dd }(\X) =
{\rm Pr}_{\QQ,\RR} \{ T \mbox{ is not separated by } \D_N
| \, |T|=\dd
\} \, .
%} | \, |T|=t\} \, .
\eea
Since $\X$ is  balanced,  the probability
$P_{N,n,\dd }(\X) $
is correctly defined; that is, it does not depend on
the choice of a particular $T $ such that $|T|=\dd $.

For a pair $(T,T') \in \T \times \T$ of different targets, set
$P(N,T,T')$ to be the  probability of the event that $ T$ and $T'$ are
not separated after $N$ random tests. If $T=T_i$ and $T'=T_j$
then, in the notation of Section 2.1, $P(1,T,T')=
p_{ij}= k_{ij}/{{n\choose s}}$, where $k_{ij}$ are the R\'{e}nyi
coefficients and $P(N,T,T')= (P(1,T,T'))^N$.

For a fixed $T$, such that $|T|=\dd $,
the probability
$P_{N,n,\dd }(\X) $
 that
after $N$ random tests
$T$ is not
separated from all $T' \neq T$,
is less
than or equal to
$
P_{N,n,\dd }(\X)
\leq
Q_{N,n,\dd }(\X)
$
where
\bea
Q_{N,n,\dd }(\X)
=
\min \{ 1, \sum_{T' \neq T}
P(N,T,T')   \}
= \min \{1,  S_1+S_2+S_3 \}\, .
\eea
Here
$$
S_1= \sum_{T': |T'| <\dd  }
P(N,T,T'), \;\;
S_2= \sum_{T' \neq T, |T'|=\dd  }
P(N,T,T'), \;\;
S_3=\sum_{T': |T'| > \dd  }
P(N,T,T')\,.
$$

One can show that
$$
S_1=
\frac{1}{{{n}\choose{\dd }}}
 \sum_{m=0}^{\dd -1}\; \sum_{ p=0}^m\!
Q(n,\dd ,m,p)
 \left( \frac{K(\PP_n^s,n,\dd ,m,p)} {{n \choose s}}\right)^N \, ,
$$
$$
S_2=
\frac{2}{{{n}\choose {\dd }}}
 \sum_{ p=0}^{\dd -1}\!
Q(n,\dd ,\dd ,p)
 \left( \frac{K(\PP_n^s,n,\dd ,\dd ,p)} {{n \choose s}} \right)^N \,,
   \,
$$
and
$$
S_3=
\frac{1}{{{n}\choose {\dd }}}
 \sum_{m=\dd +1}^{d}\; \sum_{ p=0}^\dd \!
Q(n,\dd ,m,p)
 \left(  \frac{K(\PP_n^s,n,m,\dd ,p)} {{n \choose s}}\right)^N         \, .
$$

Using the definition of $q_{\X,n,d,m,p}$ we obtain
$$
  S_1+S_2+S_3 =
\frac{1}{{{n}\choose {\dd }}}
 \sum_{m=0}^{d}\; \sum_{p=0}^{\min\{\dd ,m\}}
{\textstyle
 {{n}\choose {p\;m-p \; \dd -p\; n-\dd -m+p}}
}
 q^N_{\X,n,\dd ,m,p}\, .
$$

From the inequality
$$
\PP_N= \sum_{\dd =0}^d
{\rm Pr}_\UU\{\xi =\dd   \} P_{N,n,\dd }(\X)
\leq
 \sum_{\dd =0}^d
{\rm Pr}_\UU\{\xi =\dd   \}
Q_{N,n,\dd }(\X)\, =  \sum_{\dd =0}^d
{\rm Pr}_\UU\{\xi =\dd   \}
 \min \{1,  S_1\!+\!S_2\!+\!S_3 \} \,,
$$
we obtain:
\bea
&&{\rm Pr}_{\QQ,\RR}\{ T \textrm{ is separated by } \D_N \} = 1-\gamma \geq1- \sum_{\dd =0}^d
{\rm Pr}_\UU\{\xi =\dd   \}
Q_{N,n,\dd }(\X)\, \\
& =& 1- \sum_{\dd =0}^d
{\rm Pr}_\UU\{\xi =\dd   \}
 \min \{1,  S_1\!+\!S_2\!+\!S_3 \} = 1-\gamma^* \,. \hspace{35mm} \Box
\eea

\subsubsection*{Proof of Theorem~\ref{th:k_ijb2}}

Rewriting
(\ref{eq:k_ij-2})
 for $\KK=1$ we obtain
$$
K(\X,n,l,m,p) =
 \sum_{r{=}0}^p \sum_{u{=}0}^{m{-}p}
R(n,l,m,p,u,u,r){+}
\sum_{r{=}1}^p
\sum_{u{=}0}^{l{-}p}
\sum_{v{=}u{+}1}^{m{-}p}
R(n,l,m,p,u,v,r) +
$$
$$
\sum_{r{=}1}^p \sum_{u{=}0}^{m{-}p} \sum_{v{=}u{+}1}^{l{-}p}
R(n,\!l,\!m,\!p,\!v,\!u,\!r)+
\sum_{u{=}1}^{l{-}p}
\sum_{v{=}u{+}1}^{m{-}p}
R(n,\!l,\!m,\!p,\!u,\!v,\!0)
+
\sum_{u{=}1}^{m{-}p}\sum_{v\!=\!u+\!1}^{l-p}R(n,\!l,\!m,\!p,\!v,\!u,\!0)
$$
$$
= \sum_{r{=}1}^p
\sum_{u{=}0}^{l{-}p}\sum_{v{=}0}^{m{-}p}R(n,l,m,p,u,v,r)
+\sum_{u{=}1}^{l{-}p}
\sum_{v{=}1}^{m{-}p}
R(n,l,m,p,u,v,0)+R(n,l,m,p,0,0,0).
$$

\noindent
By using
Lemma 3.1 in \cite{zhigljavsky2003probabilistic} the following identity holds
\bea
\a{n}{s}  =
\sum_{r{=}0}^p
\sum_{u{=}0}^{l{-}p}\sum_{v{=}0}^{m{-}p}R(n,l,m,p,u,v,r)\,,
\eea
which allows us to state
\bea
K(\X,n,l,m,p) =
\a{n}{s} -\left(\sum_{u=1}^{l-p}R(n,l,m,p,u,0,0)+\sum_{v=1}^{m-p}
R(n,l,m,p,0,v,0)\right)\, .
\eea
By then applying the expression for $R(\cdot)$ given in
(\ref{eq:RR}),
we obtain
\bea
K(\X,n,l,m,p) =
\a{n}{s} -
\sum_{u=1}^{l-p}
\a{l\!-\!p}{u} \a{n\!-\!l\!-\!m\!+\!p}{s-u}
-\sum_{v=1}^{m-p}
\a{m\!-\!p}{v} \a{n\!-\!l\!-\!m\!+\!p}{s-v}
\, .
\eea
Application of the Vandermonde convolution formula
then
provides (\ref{eq:k_ijb0}). \hfill $\Box$   \\

{
\section*{Appendix B: Pseudo-code for Algorithm 1 and Algorithm 2}\label{pseudocode}

\begin{algorithm}[!h]
{
\setcounter{algocf}{-1}
\SetAlgoLined
\KwInput{A design  $\D_N$.}
\KwResult{One test containing $s$ items to be used within Algorithm 2.}
$Output  = \{\}$;\\
For each item $1,\ldots, n$, determine the frequency it appears in $\D_N$; \\
\If{there are at least $s$ items with equal smallest frequency of occurrence}{
Append to $Output $ a sample of $s$ elements from these items;
}
  \Else{ Append to $Output$ all the items with the smallest frequency of occurrence, say $s'$ of these, and sample the remaining $s-s'$ items randomly from groups that have not appeared the fewest;\\
}

\Return{Output}
\caption{}
}
\end{algorithm}

\begin{algorithm}[!h]
{
\SetAlgoLined
\KwInput{$N$ and $N':=$ The number of candidate tests. }
\KwResult{A matrix $\M=\M(\D_N)$ or equivalent design $\D_N$.}
  Construct $\M(\D_N)$ with $\D_N=\{ X_1\}$, with $X_1$ $\RR$-distributed from $\X=\PP_n^{s}$.\\
\While{Number of rows in $\M(\D_N) < N$}{
Create the $N'$ candidate tests $C_{N'}=\{ X'_1,X'_2,\ldots X'_{N'}\}$ by: repeating Algorithm 1 on $\D_N$ a total of $0.75\times N'$ times; randomly sample without replacement from $\X=\PP_n^{s}$ a total of  $0.25  \times N'$ times;\\
Construct the test matrix $\M' := \M'(C_{N'})$; \\
 Determine the row $k$ in $\M'$ (that is $\M'_{k})$ that satisfies:  $\min_{1\leq  j\leq N} d_H(\M'_{k},\M_{j})=\max_{1\leq i\leq N'}\min_{1\leq  j\leq N} d_H(\M'_{i},\M_{j})$ - if ties occur, select the item such that that $\sum_{j=1}^{N}d_H(\M'_k,\M_j)$ is highest;\\
 Append $\M'_{k}$ to the rows of $\M = \M(\D_N)$.
}
 \Return{$\M(\D_N)$}
   \caption{}
   }
\end{algorithm}
}

\bibliographystyle{plainnat}

\bibliography{Group_testing}

\end{document}